\documentclass[11pt]{article}
\newcommand{\sectionbis}[1]{\section{#1}\setcounter{equation}{0}}
\renewcommand\appendix{
  \setcounter{section}{0}%
  \setcounter{subsection}{0}%
\renewcommand{\theequation}{A.\arabic{equation}} 
 \gdef\thesection{A}\section
}
\renewcommand{\theequation}{\arabic{section}.\arabic{equation}}
\usepackage{amsmath}
\usepackage{amssymb}
\usepackage{amscd}
\newcommand{\Gen}{{\rm Gen}}
\newcommand{\scheck}{{\bf\check s}}
\newcommand{\shat}{{\bf \hat s}}
\newcommand{\xx}{{\bf x}}
\newcommand{\XX}{{\bf X}}

\newcommand{\MP}{{\rm MP}}

\newcommand{\capacity}{{\rm cap}}

\newcommand{\OO}{{\rm O}}

\newcommand{\ZZ}{{\bf Z}}
\renewcommand{\ss}{{\bf s}}

\newcommand{\codim}{{\rm codim}}
\newcommand{\GG}{{\rm G}}

\newcommand{\HH}{{\rm H}}
\newcommand{\RR}{{\rm R}}

\newcommand{\cl}{{\rm cl}}

\newcommand{\Ball}{{\rm B}}

\renewcommand{\P}{\mathbb{P}}
\newcommand{\C}{\mathbb{C}}

\newcommand{\R}{\mathbb{R}}

\newcommand{\Z}{\mathbb{Z}}

\renewcommand{\H}{{\cal H}}

\newcommand{\ddc}{{\rm dd^c}}
\renewcommand{\d}{{\rm d}}
\newcommand{\DD}{{\rm D}}
\newcommand{\QQ}{{\rm Q}}
\newcommand{\cad}{{\it c.-\`a-d. }}
\newcommand{\ie}{{\it i.e. }}

\newcommand{\vol}{{\rm vol}}

\newcommand{\Lone}{{{\rm L}^1}}
\newcommand{\Ltwo}{{{\rm L}^2}}
\newcommand{\LL}{{\rm L}}
\newcommand{\Lp}{{{\rm L}^p}}

\newcommand{\Linfty}{{{\rm L}^\infty}}

\newcommand{\FS}{{\rm FS}}
\renewcommand{\S}{{\cal S}}

\newcommand{\supp}{{\rm supp}}
\newcommand{\E}{{\cal E}}

\newcommand{\voir}{{\it voir }}

\renewcommand{\O}{{\rm O}}
\renewcommand{\o}{{\rm o}}

\newcommand{\DSH}{{\rm DSH}}

\renewcommand{\AA}{{\rm A}}
\newtheorem{theoreme}{Th\'eor\`eme}[section]
\newtheorem{proposition}[theoreme]{Proposition}
\newtheorem{corollaire}[theoreme]{Corollaire}
\newtheorem{lemme}[theoreme]{Lemme}

\newtheorem{exemple}[theoreme]{Exemple}

\newtheorem{remarque}[theoreme]{Remarque}
\newtheorem{remarques}[theoreme]{Remarques}
\newenvironment{preuve}{\begin{trivlist}
\item[]{\bf D\'emonstration.}}
{\par\hfill $\square$\end{trivlist}}
\title{Distribution des valeurs de transformations m\'eromorphes et
applications} 
\author{Tien-Cuong Dinh et Nessim Sibony}
\begin{document}
\maketitle
\begin{abstract} A meromorphic transform between complex manifolds is 
a surjective mutivalued map with an analytic graph. 
\par
Let $F_n$ be a sequence of meromorphic transforms from a compact
K\"ahler manifold $(X,\omega)$ into compact K\"ahler manifolds
$(X_n,\omega_n)$.
Let $\sigma_n$ be an appropriate probability measure on $X_n$ and
$\sigma$ the product measure of $\sigma_n$, on $\XX:=\prod_{n\geq 1} X_n$.
We give conditions which imply that for $\sigma$-almost
every $\xx=(x_1,x_2,\ldots)\in\XX$ 
$$\frac{1}{d(F_n)}
[(F_n)^*(\delta_{x_n})-(F_n)^*(\delta_{x_n'})]\rightarrow 0.$$
Here $\delta_{x_n}$ is the Dirac
mass at $x_n$ and $d(F_n)$ the maximal intermediate
degree of $F_n$.
\par
Using this formalism,
we obtain sharp results on the limit distribution of zeros,
for random $l$ 
holomorphic sections of high powers
$L^n$ of a positive holomorphic line bundle 
$L$ over a projective manifold $X$.
\par
We consider also the equidistribution problem for random
iteration of correspondences.
Assume that $f_n:X_{n-1}\longrightarrow X_n$ are
correspondences between $k$-dimensional compact manifolds and 
let $F_n:=f_n\circ\cdots\circ f_1$. We give conditions implying
that $[d(F_n)^{-1}](F_n)^*\omega_n^k$ has a limit.
In particular, when $f$ is a
meromorphic self correspondence of topological degree $d_t$
of a compact K\"ahler manifold $X$,
under a hypothesis on the dynamical degrees, we show
that $d_t^{-n}(f^n)^*\omega^k$ 
converges to a probability measure $\mu$,
satisfying $f^*\mu=d_t\mu$.
Moreover, quasi-p.s.h. functions are $\mu$-integrable. 
Every projective manifold admits such correspondences.
When $f$ is a meromorphic map, the measure $\mu$ is 
exponentially mixing.
\end{abstract}
\sectionbis{Introduction}
Une vari\'et\'e projective de type g\'en\'eral
admet au plus un nombre fini\break 
d'endomorphismes m\'eromorphes
dominants \cite{Kobayashi}. C'est dire que la dynamique de ces applications
est tr\`es pauvre. Il est souvent d\'elicat de construire des
applications m\'eromorphes d'une vari\'et\'e compacte dans
elle-m\^eme. Il n'en est plus de m\^eme d\`es qu'on consid\`ere
les correspondances. En effet, si $X$ est une vari\'et\'e projective de dimension
$k$ et si $g$ et $h$ d\'esignent deux projections holomorphes surjectives de
$X$ sur $\P^k$, le sous-ensemble analytique
$$\Gamma:=\{(x,y)\in X\times X,\ g(x)=h(y)\}$$
d\'efinit une correspondance sur $X$, \cad une fonction multivalu\'ee
$f:=h^{-1}\circ g$.
On peut aussi consid\'erer $h^{-1}\circ
u\circ g$, o\`u $u$ est un endomorphisme holomorphe de $\P^k$.
En choisissant $u$ convenablement, on obtient une correspondance
dont la dynamique est tr\`es riche (\voir aussi
\cite{ClozelUllmo,Voisin,Dinh} pour les exemples et les applications).
\par
Dans un cadre plus large: celui des transformations
m\'eromorphes, bon nombre de questions dynamiques,
ou de comportements asymptotiques de pr\'eimages,
se ram\`enent
\`a l'\'etude du probl\`eme suivant.
\\
\ \par
{\it Soient
$(X,\omega)$, $(X_n,\omega_n)$ des vari\'et\'es k\"ahl\'eriennes
compactes de dimensions respectives $k$ et $k_n$. On consid\`ere
une suite $F_n:X\longrightarrow X_n$ d'applications
m\'eromorphes, de correspondances ou plus g\'en\'eralement de
transformations m\'eromorphes et on se pose la question de donner
des crit\`eres v\'erifiables sur les $F_n$ et les $X_n$ assurant
que les pr\'eimages $F_n^{-1}(x_n)$ par $F_n$ des points $x_n\in
X_n$ sont \'equidistribu\'ees  dans $X$.
}
\\
\ \par
Pr\'ecisons les probl\`emes.
Une transformation m\'eromorphe de codimension $l$ de $X$ dans
$X_n$ est la donn\'ee d'un sous-ensemble analytique $\Gamma^{(n)}$ de
dimension pure
$k_n+l$ de $X\times X_n$, $1\leq l\leq k-1$. On suppose
que les projections $\pi$ et $\pi_n$,
restreintes \`a chaque composante irr\'eductible de $\Gamma^{(n)}$,
sur $X$ et $X_n$,
sont surjectives.
Pour $x_n\in
X_n$ g\'en\'erique, la fibre
$F_n^{-1}(x_n):=\pi(\pi_{n|\Gamma^{(n)}})^{-1}(x_n)$
est de dimension $l$. Si
$\delta_{x_n}$ d\'esigne la masse de Dirac en $x_n$, on pose
$$(F_n)^*(\delta_{x_n}):=\pi_*(\pi_{n|\Gamma^{(n)}})^*\delta_{x_n}.$$
C'est un courant de bidimension $(l,l)$ port\'e par
$F_n^{-1}(x_n)$.
\par
Sans hypoth\`ese sur
les transformations $F_n$, on ne peut
s'attendre \`a trouver en g\'en\'eral une limite des
courants $(F_n)^*(\delta_{x_n})$, convenablement normalis\'es.
Mais, $\Gamma^{(n)}$, les graphes de $F_n$, \'etant
holomorphes, on peut esp\'erer que dans les cas ``int\'eressants'', 
g\'en\'eriquement
$F_n^*(\delta_{x_n})$ et $F_n^*(\delta_{x_n'})$, convenablement
normalis\'es, ont la m\^eme
r\'epartition asymptotique.
On se pose deux probl\`emes.
\\
\
\\
{\bf Probl\`eme d'\'equidistribution.}
Trouver des conditions v\'erifiables, pour que
g\'en\'eriquement sur les suites $(x_n)$ et $(x_n')$
$$\frac{1}{d(F_n)}\left[F_n^*(\delta_{x_n})
-F_n^*(\delta_{x_n'})\right] \longrightarrow 0$$
au sens des courants. Ici $d(F_n)$ d\'esigne la masse d'une fibre
g\'en\'erique de $F_n$, celle-ci est ind\'ependante de la fibre.
\\
\
\\
{\bf Probl\`eme de convergence.}
Dans certains cas, trouver la limite de la suite de
courants
$$\frac{1}{d(F_n)}F_n^*(\delta_{x_n}).$$
Un cas particulier du probl\`eme de convergence est celui o\`u on se
donne une suite
$f_n:X_{n-1}\longrightarrow X_n$ de correspondances,
\ie une suite de transformations m\'eromorphes
de codimension $0$ entre
des vari\'et\'es
$X_{n-1}$, $X_n$ de m\^eme dimension.
On veut sous des hypoth\`eses
convenables, trouver la limite de la suite de mesures
$$\frac{1}{d_1\ldots d_n} (f_n\circ\cdots \circ f_1)^*(\delta_x)$$
et en donner les propri\'et\'es. 
On est dans le cas o\`u $F_n:=f_n\circ\cdots\circ f_1$,
le nombre de
points d'une fibre g\'en\'erique de $f_n$ \'etant \'egal \`a
$d_n$.
\\
\
\par
Commen\c cons par une situation \'etudi\'ee r\'ecemment
par Shiffman et\break 
Zelditch
\cite{ShiffmanZelditch1} et sur laquelle
il existe par ailleurs une vaste
litt\'erature classique, ainsi que de nombreux articles r\'ecents
de physiciens. C'est celle des z\'eros de polyn\^omes al\'eatoires.
Nous renvoyons \`a \cite{ShiffmanZelditch1} pour quelques
\'el\'ements bibliographiques. 
\par
Soit $L$ un fibr\'e holomorphe ample sur une vari\'et\'e
projective $X$.
On munit $L$ d'une m\'etrique dont la forme de courbure $\omega$
est positive.
Notons $\HH^0(X,L^n)$ l'espace des sections
holomorphes de $L^n$ et notons $X_n:=\P\HH^0(X,L^n)$ l'espace
projectif associ\'e. Il est
facile de d\'efinir une transformation m\'eromorphe $F_n$ telle
que pour $s_n\in X_n$, $F_n^*(\delta_{s_n})$ soit le courant
d'int\'egration $[Z_{s_n}]$ sur l'ensemble des z\'eros de $s_n$.
On se donne des mesures de probabilit\'e $\sigma_n$ sur $X_n$ et
on consid\`ere la mesure $\sigma$ 
produit des $\sigma_n$ sur $\prod_{n\geq 1} X_n$.
\begin{theoreme}{\bf (Shiffman-Zelditch \cite{ShiffmanZelditch1})}
Lorsque
$\sigma_n$ est la mesure de Lebesgue sur $X_n$ pour tout $n\geq 1$,
on a
$$\frac{1}{n}[Z_{s_n}]\longrightarrow \omega$$
pour $\sigma$-presque tout $(s_n)\in\prod_{n\geq 1} X_n$.
\end{theoreme}
\par
Nous donnons une version abstraite de cet \'enonc\'e
(th\'eor\`eme 4.1), dans le 
cadre  des transformations m\'eromorphes.
Notre m\'ethode fournit des pr\'ecisions nouvelles.
\par
1. Si $\psi$ d\'esigne une forme test de classe ${\cal C}^2$
$$\sigma_n\big(|\langle n^{-1}[Z_{s_n}]
-\omega,\psi\rangle| \geq \epsilon \big)\leq
c\|\psi\|_{{\cal C}^2}n^{mk}\exp(-\epsilon\alpha n)$$
o\`u $c>0$, $m\geq 0$ et $\alpha>0$ sont des constantes ind\'ependantes de 
$n$ et de $\epsilon$.
\par
2. On obtient le m\^eme r\'esultat en prenant pour $\sigma_n$ les
mesures de Lebesgue normalis\'ees sur la partie r\'eelle de
$X_n$. Ce qui revient \`a \'etudier les z\'eros des "sections \`a
coefficients r\'eels", par exemple, les z\'eros complexes des
polyn\^omes homog\`enes \`a coefficients r\'eels pour le cas o\`u $X=\P^k$.
\par
3. On obtient un th\'eor\`eme d'\'equidistribution, avec
des estimations analogues, pour les z\'eros communs de $l$
sections holomorphes:
$$\frac{1}{n^l}[Z_{s_n^1}\cap\ldots\cap Z_{s_n^l}]\longrightarrow
\omega^l.$$
\par
Il s'agit dans ces r\'esultats de r\'esoudre le probl\`eme
d'\'equidistribution. En effet,
un th\'eor\`eme de
Tian-Zelditch \cite{Tian, Zelditch} (\voir le th\'eor\`eme 7.2)
implique que la moyenne des courants
$n^{-l}[Z_{s_n^1}\cap\ldots\cap Z_{s_n^l}]$
converge vers $\omega^l$ quand $n$ tend vers l'infini. 
\\
\
\par
Dans le cas g\'en\'eral, notre crit\`ere d'\'equidistribution
utilise deux notions: l'une li\'ee \`a la croissance des
transformations $F_n$ et l'autre \`a la g\'eom\'etrie des $X_n$.
En particulier si les $X_n$ sont \'egales, seule la croissance
des $F_n$ intervient.
\par
Les indicateurs de croissance sont les degr\'es interm\'ediaires
d'ordres $k_n$ et $k_n-1$,
classiques en th\'eorie de distribution des valeurs, associ\'es aux $F_n$:
$$d(F_n):=\int_X F_n^*(\omega_n^{k_n})\wedge \omega^{k-l} \ \ \
\mbox{ et }\ \ \
\delta(F_n):=\int_XF_n^*(\omega_n^{k_n-1})\wedge\omega^{k-l+1}.$$
C'est le comportement de la suite des $\delta(F_n)d(F_n)^{-1}$
qui joue un r\^ole. Observons que le calcul de $d(F_n)$ et
$\delta(F_n)$ est cohomologique, $d(F_n)$ est la masse d'une
fibre g\'en\'erique de $F_n$ et par exemple dans le cas de l'espace
projectif $\delta(F_n)$ est la masse de l'image r\'eciproque 
d'une droite g\'en\'erique.
\par
La g\'eom\'etrie des vari\'et\'es intervient par
l'interm\'ediaire des meilleures constantes pour r\'esoudre
$\ddc$ dans une classe de cohomologie donn\'ee.
Pour les 
extensions ci-dessus du th\'eor\`eme de
Shiffman-Zelditch, la propri\'et\'e
suivante est cruciale: 
{\it pour toute fonction $\varphi$ q.p.s.h. sur
$\P^k$ v\'erifiant $\ddc\varphi\geq -\omega_\FS$ et
$\int_{\P^k}\varphi\d m=0$ on a
$$\max_{\P^k}\varphi \leq c(1+\log k).$$}
Ici $m$ est la mesure invariante sur $\P^k$ ou sur $\R\P^k$,
$\omega_\FS$ \'etant la forme de Fubini-Study sur $\P^k$
normalis\'ee par $\int_{\P^k}\omega_\FS^k=1$ et $c>0$ est une
constante ind\'ependante de $k$. On
calcule dans cet exemple que $\delta(F_n)=1$ et $d(F_n)=\O(n)$. 
Les estimations des constantes g\'eom\'etriques 
pour les vari\'et\'es multiprojectives permettent
de d\'eduire la convergence pour les z\'eros communs de plusieurs
sections holomorphes. C'est bien s\^ur la d\'ependance
par rapport \`a la dimension  qui
est importante.
\\
\ \par
On peut souligner deux autres raisons pour s'int\'eresser au comportement de 
suites d'applications et pas seulement aux it\'er\'es d'une application. 
La premi\`ere est que du point de vue physique, on compose des applications voisines,
et un th\'eor\`eme dans ce cadre est signe de la robustesse du r\'esultat.
La seconde est que pour l'it\'eration d'applications birationnelles $f$, on est 
amen\'e \`a consid\'erer le comportement de la suite $(f^n,f^{-n})$. Dans ce cas,
on obtient des r\'esultats qui sont beaucoup moins \'evidents si on  
en consid\`ere s\'epar\'ement les suites $f^n$ et $f^{-n}$.
\\
\
\par
Les r\'esultats g\'en\'eraux que nous obtenons sont
d\'emontr\'es 
au paragraphe 4. Nous en donnons des applications dans les
paragraphes suivants. Pr\'ecisons un cas simple. 
\begin{theoreme} Soient $F_n:X\longrightarrow X'$ des
transformations m\'eromorphes de codimension $l$ entre vari\'et\'es
k\"ahl\'eriennes compactes de dimensions respectives $k$ et $k'$. 
Soient $\delta_n,d_n$ leurs degr\'es
interm\'ediaires d'ordre $k'-1$ et $k'$. Notons $\E$ l'ensemble (exceptionnel) des
$x'\in X$ tels que 
$$ \frac{1}{d_n}[(F_n)^*(\delta_{x'})-(F_n)^*({\omega'}^{k'})]
\ \ \mbox{ ne converge pas vers }0.$$
\begin{enumerate}
\item[{\rm (1)}] Si $\sum\delta_n d_n^{-1}<+\infty$ alors $\E$ est
pluripolaire. 
\item[{\rm (2)}] Si $\sum \exp(-\delta_n^{-1}d_n t)<+\infty$ pour tout $t>0$
alors $\E$ est de mesure de Lebesgue nulle
et de $\sigma$-mesure nulle
pour toute mesure mod\'er\'ee $\sigma$.
\end{enumerate}
\end{theoreme}
Les mesures de Lebesgue sur $X$ ou sur une
sous-vari\'et\'e analytique totalement r\'eelle de dimension maximale 
sont des exemples de mesures mod\'er\'ees.
\\
Dans ce contexte, le cas (1) avec des applications rationnelles
$F_n$ de $\P^k$ dans $\P^{k'}$ a \'et\'e \'etudi\'e par 
Russakovskii-Sodin et Russakovskii-Shiffman
\cite{RussakovskiiSodin, RussakovskiiShiffman}. 
Pour le cas des it\'er\'es d'une application holomorphe 
de $\P^k$ \voir \cite{FornaessSibony2, Sibony2, BriendDuval,
DinhSibony2}.
\\
\ \par
Les transformations $F_n$ \'etant al\'eatoires, la suite 
$d_n^{-1}(F_n)^*({\omega'}^{k'})$ n'a pas de limite en g\'en\'eral. 
Lorsque les transformations $F_n$ sont les it\'er\'es
$f^n$ d'une correspondance $f:X\longrightarrow X$, on obtient le
th\'eor\`eme suivant qui fournit une solution au probl\`eme de 
convergence dans ce cas (d'autres variantes de ce r\'esultat pour les 
it\'erations al\'eatoires sont donn\'ees au paragraphe 5).
\begin{theoreme} Soit $f:X\longrightarrow X$ une correspondance
m\'eromorphe de degr\'e topologique $d_t$
sur une vari\'et\'e k\"ahl\'erienne compacte
$(X,\omega)$ de dimension $k$. Supposons que le degr\'e dynamique
d'ordre $k-1$
$$d_{k-1}:=\limsup_{n\rightarrow \infty} \left(\int_X
(f^n)^*\omega^{k-1}\wedge \omega\right)^{1/n}$$ 
v\'erifie  $d_{k-1}<d_t$.
Alors la suite
$\mu_n:=d_t^{-n}(f^n)^*(h_n\omega_n^k)$ converge vers une mesure
$\mu$ 
$f^*$-invariante, \ie $f^*\mu=d_t\mu$. La convergence est
uniforme sur les $h_n$ positives v\'erifiant
$\int h_n\omega^k=1$ et
$\|h_n\|_{\Ltwo(X)}^{1/n} =\o(d_{k-1}^{-1}d_t)$. Lorsque $f$ 
est une application m\'eromorphe, la mesure
$\mu$ est m\'elangeante \`a vitesse exponentielle. De plus toute
fonction q.p.s.h. $\varphi$ est $\mu$-int\'egrable ($\mu$ est PLB) et 
$\langle\mu_n,\varphi\rangle \longrightarrow\langle\mu,\varphi\rangle$.
\end{theoreme} 
\par
En appliquant les m\'ethodes de Lyubich \cite{Lyubich},
Briend-Duval \cite{BriendDuval} et
\cite{DinhSibony2,Dinh}, on peut montrer que les points
p\'eriodiques r\'epulsifs sont denses dans le support de $\mu$ et
que l'ensemble exceptionel $\E$ est une r\'eunion finie ou d\'enombrable
d'ensembles analytiques. Lorque $f$ est une application
rationnelle sur une vari\'et\'e projective, V. Guedj \cite{Guedj}
a, dans un travail r\'ecent, construit la
mesure $\mu$ et prouv\'e les propri\'et\'es ci-dessus.
Dans \cite{DinhSibony3}, nous avons montr\'e que $\mu$ est
d'entropie maximale $\log d_t$.
Observons que toute vari\'et\'e projective admet des
correspondances satisfaisant l'hypoth\`ese du th\'eor\`eme 1.3.
\\
\
\par
Notre approche du probl\`eme d'\'equidistribution g\'en\'eral
reprend celle que nous avons utilis\'ee
pour l'\'etude des applications
\`a allure polynomiale dans une vari\'et\'e de Stein
\cite{DinhSibony2}. Ici, apr\`es avoir d\'efini convenablement les images r\'eciproques
$F_n^*$ et les images directes $(F_n)_*$ des courants, il faut \'evaluer
\begin{eqnarray}
\frac{1}{d(F_n)}\left\langle F_n^*(\delta_{x_n})
-F_n^*(\omega_n^{k_n}), \psi\right\rangle & = &
\frac{1}{d(F_n)}\left\langle \delta_{x_n}
-\omega_n^{k_n}, (F_n)_*\psi\right\rangle
\end{eqnarray}
pour une forme test lisse $\psi$ de bidegr\'e $(l,l)$, la
forme de K\"ahler $\omega_n$
\'etant normalis\'ee par $\int_{X_n}\omega_n^{k_n}=1$.
\par
L'id\'ee est de remplacer $(F_n)_*\psi$ par une autre solution
$\psi_n$ de l'\'equation $\ddc\psi_n=\ddc(F_n)_*\psi$, \cad de
retrancher \`a $(F_n)_*\psi$ une constante convenable. 
Les mesures $\delta_{x_n}$ et $\omega_n^{k_n}$ \'etant de m\^eme masse, 
le membre \`a droite de (1.1) ne change pas lorsqu'on remplace $(F_n)_*\psi$
par $\psi_n$. 
On estime
$\psi_n$ en fonction de $\ddc(F_n)_*\psi$, d'o\`u l'introduction
de $\delta(F_n)$ qui essentiellement mesure la masse de $\ddc(F_n)_*\psi$. 
\par
Les propri\'et\'es de convergence de la suite
$\delta(F_n)[d(F_n)]^{-1}$ ont pour cons\'equence des
propri\'et\'es de convergence de (1.1) vers $0$ sauf sur des
ensembles dont nous pouvons majorer la mesure.
Lorsque la s\'erie $\sum\delta_n d_n^{-1}$
n'est pas convergente (c'est le cas pour la distribution des
z\'eros de sections holomorphes al\'eatoires de fibr\'es positifs) 
on est amen\'e \`a
utiliser la notion de mesure
mod\'er\'ee. 
On obtient gr\^ace 
\`a cela des r\'esultats
d'\'equidistribution presque partout comme au th\'eor\`eme 1.2
en supposant seulement que 
$\delta_nd_n^{-1}=\o(1/\log n)$. 
\par
Cette m\'ethode de dualit\'e permet d'obtenir des estim\'ees exponentielles 
des volumes de l'ensemble des ``mauvaises sections'' dans le th\'eor\`eme de 
Shiffman-Zelditch, ou encore des estim\'ees de
la vitesse de m\'elange dans le th\'eor\`eme 1.3.
Bien s\^ur, dans les applications, que nous donnons, 
il faut calculer les constantes g\'eom\'etriques
ainsi que les degr\'es dynamiques.
\\
\ \par

Revenons sur la solution du probl\`eme de convergence pour les correspondances:
th\'eor\`eme 1.3.
On montre 
que pour toute fonction q.p.s.h. $\varphi$, $\langle d_t^{-n} 
(f^n)^*\omega^k,\varphi\rangle$ est convergente. 
Pour cela, \'etant donn\'e une fonction q.p.s.h. $\varphi$ sur $X$. On pose 
$$b_0:=\int\varphi\omega^k,\ \ \ \ \ \ \ \varphi_0=\varphi-b_0$$
et
$$b_n:=\int f_*\varphi_{n-1}\omega^k,\ \ \ \ \ \ \
\varphi_n:=f_*\varphi_{n-1}-b_n.$$
On a par it\'eration
$$\left\langle \frac{(f^n)^*\omega^k}{d_t^n},\varphi\right\rangle
= b_0+\frac{b_1}{d_t}+
\cdots + \frac{b_n}{d_t^n} +\left\langle \omega^k,
\frac{\varphi_n}{d_t^n}\right\rangle.$$
Puisque les $\varphi_n$ sont convenablement normalis\'ees, on peut estimer 
$\varphi_n$ \`a l'aide de $\ddc\varphi_n$ et l'hypoth\`ese $d_{k-1}<d_t$ 
implique $\langle\omega^k,d_t^{-n}\varphi_n\rangle \rightarrow 0$. 
La ``continuit\'e'' de $f_*$ implique  
que $|b_n|$ est de l'ordre $d_{k-1}^n$. Ce qui d\'emontre 
la convergence.
\\
\ 
\par
Lorsque $f:X\longrightarrow X$
est une application birationnelle on a toujours $d_t=1$ donc 
les hypoth\` eses du th\'eor\`eme 1.3 ne sont pas v\'erifi\'ees. Nous obtenons 
cependant un r\'esultat d'\'equidistribution en consid\'erant simultan\'ement 
le graphe de $f^n$ et $f^{-n}$. On peut alors prendre des
intersections de ces graphes  
avec des sous-vari\'et\'es et
appliquer notre th\'eor\`eme abstrait. Donnons un exemple:
soit $f$ un automorphisme polynomial r\'egulier de $\C^k$
\cite{Sibony2}, notons
$d_+$ le degr\'e alg\'ebrique de $f$ et $d_-$ celui de $f^{-1}$. 
On montre \cite{Sibony2} que le courant
$$\frac{1}{(d_+)^{nl}(d_-)^{ml'}}[f^{-n}(\omega_\FS^l)\cap f^m(\omega_\FS^{l'})]$$
tend faiblement vers un courant positif ferm\'e
$T_{l,l'}$ de bidegr\'e $(l+l',l+l')$ lorsque $m$ et $n$ tendent vers l'infini.
Les entiers $l$, $l'$ admissibles sont d\'efinis
par les ensembles d'ind\'etermination de $f$ et $f^{-1}$.
Ce courant $T_{l,l'}$  a des propri\'et\'es d'invariance par rapport \`a $f$ qui 
en font un objet dynamique int\'eressant. Il peut \^etre obtenu autrement.
\par
Notons $\GG_l$
la grassmannienne des plans projectifs 
de dimension $k-l$ dans $\P^k$. Il existe un ensemble 
pluripolaire $\E\subset \GG_l\times \GG_{l'}$, $l\geq 1$, $l'\geq 1$, tel que
pour $(x,x')\in\GG_l\times\GG_{l'}\setminus \E$ la suite de courants
$$\frac{1}{(d_+)^{-nl}(d_-)^{-ml'}}[f^{-n}(\P^{k-l}_x)\cap f^m(\P_{x'}^{k-l'})]$$
converge faiblement vers 
$T_{l,l'}$ (\voir th\'eor\`eme 6.3).
\par
Pour d\'emontrer ce r\'esultat, on r\'esout le probl\`eme d'\'equidistribution
pour la suite de transformations m\'eromorphes $F_{n,m}:\P^k\longrightarrow
\GG_l\times \GG_{l'}$ d\'efinies par les relations 
$F^{-1}_{n,m}(x,x')=f^{-n}(\P^{k-l}_x)\cap f^m(\P_{x'}^{k-l'})$.
\\
\ 
\par
L'article est organis\'e de la mani\`ere suivante. Au paragraphe 2, nous 
donnons quelques r\'esultats de compacit\'e pour les fonctions q.p.s.h. 
dans une vari\'et\'e k\"ahl\'erienne compacte. 
C'est la base de notre approche. On utilise les fonctions q.p.s.h. pour tester 
la convergence. Leurs propri\'et\'es de 
compacit\'e jouent un r\^ole cl\'e dans les estimations.
\par
Nous d\'efinissons les op\'erations
image directe et image r\'eciproque des courants dans les cas que nous utilisons.
\par
Le paragraphe 3 introduit les op\'erations de composition, de produit
et d'intersection, sur les transformations
m\'eromorphes. Nous donnons
les estimations fondamentales des degr\'es interm\'ediaires d'une compos\'ee
ou d'un produit de transformations m\'eromorphes sur les
vari\'et\'es projectives.
\par
Au paragraphe 4 nous d\'emontrons notre r\'esultat de convergence abstrait 
(th\'eor\`eme 4.1) sur l'\'equidistribution, pour des transformations  m\'eromorphes
g\'en\'erales. Nous l'appliquons ensuite aux diverses situations que nous avons
discut\'ees. Dans un appendice, nous avons rassembl\'e quelques 
propri\'et\'es des ensembles pluripolaires et les estimations des constantes 
g\'eom\'etriques qui nous sont n\'ecessaires. Nous introduisons
aussi une notion de capacit\'e dans les vari\'et\'es compactes
qui peut \^etre utile dans d'autres questions.
\\
\
\par
Dans tout l'article $(X,\omega)$, $(X',\omega')$ et $(X_n,\omega_n)$, $n\geq 0$ 
d\'esignent des vari\'et\'es k\"ahl\'eriennes compactes de dimensions respectives
$k$, $k'$ et $k_n$. On suppose que $\int_{X'}{\omega'}^{k'}=1$ et 
$\int_{X_n}\omega_n^{k_n}=1$ pour tout $n\geq 1$. Les constantes 
g\'eom\'etriques $r(X,\omega)$, $\RR^*_1(X,\omega)$, $\RR^*_2(X,\omega,p)$ sont 
d\'efinies au paragraphe 2.1,
$\RR_i(X,\omega,\sigma)$, $\Delta(X,\omega,\sigma,t)$ au
paragraphe 2.2,
les degr\'es interm\'ediaires $\lambda_l(F)$, les degr\'es
dynamiques 
$d_l(F)$, le degr\'e topologique $d_t(F)$ et $\AA(f)$ sont
d\'efinis  au paragraphe 3. 
L'espace projectif $\P^k$ est muni de la forme de Fubini-Study 
$\omega_\FS$ normalis\'ee par $\int\omega_\FS^k=1$, de la mesure de probabilit\'e
invariante $\Omega_\FS$ et de la mesure de probabilit\'e invariante $m_\FS$ 
sur sa partie 
r\'eelle $\R\P^k$. L'espace multiprojectif $\P^{k,l}=\P^k\times\cdots
\times \P^k$ ($l$ fois) est muni de la forme de K\"ahler $\omega_\MP$, de la mesure 
de probabilit\'e invariante naturelle $\Omega_\MP$ et de la mesure 
de probabilit\'e invariante naturelle $m_\MP$ sur sa partie r\'eelle
(\voir appendice pour les d\'etails). 
\sectionbis{Pr\'eliminaires}
Dans ce paragraphe, nous donnons quelques propri\'et\'es 
des fonctions q.p.s.h. sur une vari\'et\'e
k\"ahl\'erienne compacte.
Nous d\'efinissons des
constantes\break 
g\'eom\'etriques relatives \`a la r\'esolution de
$\ddc$. Nous introduisons les op\'erateurs 
image directe et image
r\'eciproque d'un courant par une application holomorphe.
\par
Si $S$ est un courant r\'eel ferm\'e de bidegr\'e
$(r,r)$ de $X$, notons $\cl(S)$ sa classe
dans le groupe de cohomologie de
Dolbeault
$$\H^{r,r}(X,\R):=\H^{r,r}(X,\C)\cap\H^{2r}(X,\R).$$
On dira que $\cl(S)\leq \cl(S')$ si la classe
de $S'-S$ peut \^etre repr\'esent\'ee
par un courant positif ferm\'e. Si $S$, $S'$
sont des courants positifs ferm\'es de bidegr\'e $(r,r)$
et si $\cl(S)\leq \cl(S')$, leurs masses v\'erifient 
$\|S\|\leq \|S'\|$, o\`u on a pos\'e $\|S\|:=\int_X S\wedge \omega^{k-r}$.
Lorsque $S$ est positif, sa masse ne d\'epend que
de la classe $\cl(S)$. 
Les espaces $\Lp(X)$ sont d\'efinis par rapport \`a la forme volume
$\omega^k$.
\par\ 
\\
{\bf 2.1}. {\it Fonctions q.p.s.h. et fonctions d.s.h.}
\\
\
\par
Une fonction $\varphi$ sur $X$ est
{\it quasi-plurisousharmonique (q.p.s.h.)}
si elle
est int\'egrable, semi-continue sup\'erieurement (s.c.s.)
et v\'erifie  
$\ddc\varphi\geq -c\omega$ au sens des courants,
pour une constante $c\geq 0$. 
Une telle fonction $\varphi$ appartient \`a
$\Lp(X)$ pour tout $p\geq 1$. En effet, localement elle diff\`ere
d'une fonction p.s.h. par une fonction lisse. 
Pour toute suite $(\varphi_n)$ de
fonctions q.p.s.h. n\'egatives
v\'erifiant $\ddc\varphi_n\geq -\omega$,
on peut
extraire une sous-suite qui, ou bien, converge dans tout $\Lp(X)$, 
$p\geq 1$, vers une fonction q.p.s.h.
$\varphi$ v\'erifiant $\ddc\varphi\geq
-\omega$, ou bien, converge uniform\'ement vers $-\infty$,
\cite[p.94]{Hormander1}.
\begin{proposition} La famille des fonctions
$\psi$ q.p.s.h. v\'erifiant $\ddc\psi\geq
-\omega$ et l'une des deux conditions de normalisation
$$\max_X\varphi=0 \ \ \ \mbox{ ou } \ \ \ \int_X\varphi\omega^k=0$$
est compacte dans $\Lp(X)$ pour tout $p\geq 1$. De plus,
ces fonctions sont born\'ees sup\'erieurement par une m\^eme
constante. 
\end{proposition}
\begin{preuve} Soit $(\varphi_n)$ une suite de fonctions v\'erifiant
l'hypoth\`ese du lemme.
Notons $a_n:=\sup_X \varphi_n$. Aucune sous-suite de
$(\varphi_n-a_n)$ ne peut tendre uniform\'ement vers $-\infty$ car
$\sup_V (\varphi_n-a_n)=0$.
Par cons\'equent, la suite $(\varphi_n-a_n)$ est
born\'ee dans $\Lp(X)$ pour tout $p\geq 1$. Les conditions de
normalisation impliquent que
$$a_n=0\ \ \ \mbox{ ou } \ \ \
a_n\int_X\omega^k=-\int_X(\varphi_n-a_n)\omega^k.$$
Donc la suite $(a_n)$ est born\'ee. Ceci implique que
la suite $(\varphi_n)$ est
born\'ee dans $\Lp(X)$ et qu'on peut extraire des 
sous-suites convergentes.  
\end{preuve}
La proposition A.3 de l'appendice fournit des estimations sur les
int\'egrales $\int_X\exp(-\alpha\varphi)\omega^k$ avec $\alpha>0$,
pour $\varphi$ q.p.s.h.
\begin{proposition} Il existe une constante $r>0$ telle que, pour tout
courant positif ferm\'e $T$ de bidegr\'e $(1,1)$ et de masse $1$,
il existe une $(1,1)$-forme lisse $\alpha$, qui ne
d\'epend que de $\cl(T)$, et 
une fonction q.p.s.h. $\varphi$ v\'erifiant
$-r\omega\leq \alpha\leq r\omega$ et $\ddc\varphi -T=\alpha$.
\end{proposition}
\begin{preuve} On choisit des formes r\'eelles lisses
$\alpha_i$ de bidegr\'e $(1,1)$ avec $i=1,\ldots,m$ telles que
les classes $\cl(\alpha_i)$ engendrent le groupe de cohomologie
de Dolbeault $\H^{1,1}(X,\R)$. La famille des courants
$T$ de masse $1$
\'etant compacte, il existe une constante $c>0$ ind\'ependante
de $T$ et des nombres r\'eels $c_1,\ldots,c_m$ tels que
$\cl(T)=\sum c_i \cl(\alpha_i)$ avec $|c_i|\leq c$.
La derni\`ere relation entra\^{\i}ne l'existence
d'une fonction q.p.s.h.
$\varphi$
telle que $\ddc\varphi=T-\sum c_i \alpha_i$.
Soit $r>0$ une
constante telle que $-r\omega\leq\sum c_i\alpha_i\leq r\omega$
pour tous les $c_i$
v\'erifiant $|c_i|\leq c$. La constante $r$, la
forme $\alpha:=\sum c_i\alpha_i$ et la fonction $\varphi$
v\'erifient la proposition.
\end{preuve}
\par
On dit qu'un sous-ensemble de $X$ est {\it pluripolaire} 
s'il est contenu dans $\{\varphi=-\infty\}$ o\`u $\varphi$ 
est une fonction q.p.s.h. (\voir appendice).  
On appelle {\it fonction d.s.h.} toute fonction, d\'efinie hors d'un
sous-ensemble pluripolaire, qui s'\'ecrit comme
diff\'erence de deux
fonctions q.p.s.h. Deux fonctions d.s.h. sur $X$
sont {\it \'egales} si elles sont \'egales hors d'un
ensemble pluripolaire.  
Notons $\DSH(X)$ l'espace des fonctions d.s.h. sur $X$.
On v\'erifie facilement que, pour
une fonction d.s.h. $\psi$ sur $X$,
il existe deux courants $T^\pm$ 
positifs ferm\'es de bidegr\'e
$(1,1)$ tels que $\ddc\psi=T^+-T^-$. 
On a $\cl(T^+)=\cl(T^-)$ et
$\|T^+\|=\|T^-\|$. R\'eciproquement, d'apr\`es la proposition 2.2,
si $T^\pm$ sont deux
courants positifs ferm\'es de bidegr\'e $(1,1)$ tels que
$\cl(T^+)=\cl(T^-)$, alors il existe une fonction $\psi$ d.s.h. sur $X$
telle que $\ddc\psi=T^+-T^-$.
On peut choisir $\psi$ telle que $\int_X\psi\omega^k=0$.
\\
\ \par
Notons $r(X,\omega)$ la borne inf\'erieure des constantes $r$
qui v\'erifient la proposition 2.2. Posons
\begin{eqnarray}
\QQ(X,\omega) & := & \Big\{ \varphi\mbox{ q.p.s.h. sur } X,\ 
\ddc\varphi \geq -r(X,\omega)\omega\Big\}
\end{eqnarray}
Si $\dim\H^{1,1}(X,\R)=1$
on a $r(X,\omega)=1$. C'est le cas si $X$ est
l'espace projectif $\P^k$
muni de la forme de Fubini-Study $\omega_\FS$.
En g\'en\'eral,
on a $r(X,\omega)\geq 1$. Pour
l'espace multiprojectif
$\P^{k,l}$ muni de la forme de
K\"ahler $\omega_\MP$, on a
$r(\P^{k,l},\omega_\MP)\leq cl$ o\`u $c>0$ est une constante ind\'ependante de $k$
(\voir appendice A.11).
\begin{remarque}\rm
D'apr\`es Kodaira-Spencer \cite[p.73]{KodairaSpencer}, si
$(X_t)$ est une famille lisse de vari\'et\'es k\"ahl\'eriennes
compactes,
alors $\dim \H^{1,1}(X_t,\R)$ est localement constante et on peut
trouver des formes de K\"ahler $\omega_t$ sur $X_t$ qui
d\'ependent de fa\c con ${\cal C}^\infty$ du param\`etre $t$. On
en d\'eduit que les constantes $r(X_t,\omega_t)$ sont localement
major\'ees.
\end{remarque}
\par
Observons que deux fonctions $\psi_1$, $\psi_2$ dans $\Lone(X)$
diff\`erent par une
constante si et seulement si $\ddc\psi_1=\ddc\psi_2$. 
Nous d\'efinissons deux constantes positives
li\'ees \`a la r\'esolution de
$\ddc$ sur $(X,\omega)$ pour des solutions normalis\'ees.
Supposons que $\int_X\omega^k=1$. 
Posons pour tout $p\geq 1$
\begin{eqnarray}
\RR^*_1(X,\omega) & := &
\sup_\varphi \left\{\max_X\varphi,\ \varphi\in \QQ(X,\omega),
\int\varphi \omega^k=0\right\} \nonumber\\
& = & \sup_\varphi \left\{-\int\varphi \omega^k,\ \varphi\in \QQ(X,\omega),
\max_X\varphi =0\right\}\\
\RR^*_2(X,\omega,p) & := & \sup_\varphi
\left\{\|\varphi\|_{\Lp(X)},\ \varphi\in \QQ(X,\omega),
\int\varphi \omega^k=0\right\}
\end{eqnarray}
On verra \`a la proposition 2.5, que
$\RR^*_2(X,\omega,1)\leq 2\RR^*_1(X,\omega)$. 
\par
\
\\
{\bf 2.2}. {\it Mesures PLB et mesures mod\'er\'ees.}
\\
\ 
\par
Soit $\mu$ une mesure positive sur $X$. On
dira que $\mu$ est {\it PLB} si les fonctions
q.p.s.h. sont $\mu$-int\'egrables.
Dans le cas de dimension $1$, $\mu$ est PLB si et
seulement si elle admet localement un potentiel born\'e
\cite{DinhSibony2}. Il est clair que les mesures PLB ne chargent
pas les ensembles pluripolaires de $X$. 
On montre \`a la proposition A.1, qu'elles ne chargent
pas les sous-ensembles analytiques propres de $X$.
\\
\ \par
Soient $c>0$ et $\alpha>0$.
Nous dirons que $\mu$ est {\it
$(c,\alpha)$-mod\'er\'ee} si
$$\int_X \exp(-\alpha\varphi)
\d\mu\leq c$$
pour
toute $\varphi$ q.p.s.h. v\'erifiant $\ddc\varphi\geq -\omega$ et
$\max_X\varphi=0$. On d\'eduit d'un r\'esultat classique
\cite[p.105]{Hormander2} que 
la mesure $\omega^k$ est 
$(c,\alpha)$-mod\'er\'ee pour certaines constantes $c>0$ et $\alpha>0$
(\voir proposition A.3). On verra aussi que 
les mesures invariantes sur les sous-espaces
projectifs r\'eels $\R\P^k$ de $\P^k$ sont mod\'er\'ees
(\voir proposition A.9). Si une mesure $\mu$ de
$\P^1$ v\'erifie localement
$\int_\C |z-\xi|^{-\alpha}\d\mu(\xi)\leq A$ pour
une constante $A>0$ alors elle est $(c,\alpha)$-mod\'er\'ee pour une
constante $c>0$ convenable. On a not\'e $z,\xi$ des coordonn\'ees
affines sur une carte $\C\subset\P^1$. Il est clair que toute mesure mod\'er\'ee
est PLB. 
\begin{proposition}
Soit $\mu$ une mesure PLB sur $X$. La famille
des fonctions q.p.s.h. 
$\varphi$, v\'erifiant $\ddc\varphi\geq -\omega$, et
l'une des deux conditions
de normalisation 
$$\max_X\varphi=0 \ \  \mbox{ ou }\ \  \int\varphi\d\mu=0$$ 
est born\'ee dans $\Lone(\mu)$ et est born\'ee
sup\'erieurement. De plus, 
il existe $c>0$ ind\'ependant de $\varphi$ tel que
pour tout $t>0$ 
on ait $\mu(\varphi<-t)\leq ct^{-1}$.
\end{proposition}
\begin{preuve} Soit $(\varphi_n)$ une suite de fonctions
q.p.s.h. v\'erifiant $\ddc\varphi_n\geq -\omega$
et $\max_X\varphi_n=0$.
Montrons que $(\varphi_n)$ est born\'ee dans $\Lone(\mu)$. 
Sinon, quitte \`a
extraire une sous-suite, on peut supposer que
$\int\varphi_n\d\mu\leq -n^2$.
Posons
$\Phi:=\sum n^{-2}\varphi_n$. D'apr\`es la proposition 2.1,
$\Phi$ est une fonction q.p.s.h.
v\'erifiant $\ddc\Phi\geq -2\omega$. On a $\int_X
\Phi\d\mu=-\infty$. Cela contredit que $\Phi$ soit
$\mu$-int\'egrable. 
\\
\ \par
Soit maintenant $(\varphi_n)$ une suite de fonctions q.p.s.h.
v\'erifiant $\ddc\varphi_n\geq -\omega$ et
$\int\varphi_n\d\mu=0$. 
Posons $a_n:=\max_X\varphi_n$ et
$\widetilde\varphi_n:=\varphi_n-a_n$. On a $\sup_X
\widetilde\varphi_n=0$. D'apr\`es 
la partie pr\'ec\'edente, $(\widetilde\varphi_n)$
est born\'ee dans $\Lone(\mu)$. 
Or $a_n=-\int\widetilde\varphi_n\d\mu$, donc $(a_n)$ est
born\'ee et par suite $(\varphi_n)$ est born\'ee dans $\Lone(\mu)$.
On en d\'eduit aussi que $(\varphi_n)$ est born\'ee sup\'erieurement. 
\par
La famille de fonctions q.p.s.h. consid\'er\'ees
\'etant born\'ee dans $\Lone(\mu)$, il existe $c>0$
tel que $\|\varphi\|_{\Lone(\mu)}\leq c$ pour tout 
$\varphi$ dans cette famille. On a donc
$$\mu(\varphi<-t)\leq t^{-1}\|\varphi\|_{\Lone(\mu)}\leq ct^{-1}.$$ 
\end{preuve}
\par
Soit $\mu$ une mesure de probabilit\'e PLB. Il r\'esulte de la proposition 2.4
qu'on peut d\'efinir 
les meilleures constantes pour la r\'esolution de $\ddc$,
avec une normalisation associ\'ee
\`a $\mu$. Posons
\begin{eqnarray}
\RR_1(X,\omega,\mu) & := & \sup_\varphi
\left\{\max_X\varphi,\ \varphi\in\QQ(X,\omega),
\int\varphi\d\mu=0\right\}\nonumber\\
 & = & \sup_\varphi
\left\{-\int\varphi\d\mu,\ \varphi\in\QQ(X,\omega),
\max_X\varphi=0\right\}\\
\RR_2(X,\omega,\mu) & := & \sup_\varphi
\left\{\|\varphi\|_{\Lone(\mu)},
\  \varphi\in\QQ(X,\omega),\int\varphi\d\mu=0\right\}\\
\RR_3(X,\omega,\mu) & := & \sup_\varphi
\left\{\left|\int\varphi\omega^k\right|,
\  \varphi\in\QQ(X,\omega),\int\varphi\d\mu=0\right\}\nonumber\\
 & := & \sup_\varphi
\left\{\left|\int\varphi\d\mu\right|,
\  \varphi\in\QQ(X,\omega),\int\varphi\omega^k=0\right\}
\end{eqnarray}
Pour tout $t\in\R$, posons
\begin{eqnarray}
\Delta(X,\omega,\mu,t):=\sup_\varphi\left\{\mu(\varphi<-t),\
\varphi\in\QQ(X,\omega), \int\varphi\d\mu=0\right\}
\end{eqnarray}
\begin{proposition} Soit $\mu$ une mesure de probabilit\'e
PLB sur $(X,\omega)$. On a $\RR_2(X,\omega,\mu)\leq
2\RR_1(X,\omega,\mu)$
et $\RR_3(X,\omega,\mu)\leq \RR_1(X,\omega,\mu)+\RR^*_1(X,\omega)$. 
Si $\mu$ est $(c,\alpha)$-mod\'er\'ee alors
$\Delta(X,\omega,\mu,t)\leq c\exp(-\alpha r^{-1} t)$ o\`u
$r:=r(X,\omega)$.
\end{proposition}
\begin{preuve} Soit $\varphi$ une fonction q.p.s.h. telle que
$\ddc\varphi\geq -r\omega$ et 
$\int\varphi\d\mu=0$. Posons $m:=\max_X\varphi$. On a $m-\varphi
\geq 0$
et donc 
\begin{eqnarray*}
\int|\varphi|\d\mu & \leq & \int|\varphi-m|\d\mu+m=
\int(m-\varphi)\d\mu +m\\
& = & 2m-\int\varphi\d\mu=2m\leq 2\RR_1(X,\omega,\mu).
\end{eqnarray*}
Donc $\RR_2(X,\omega,\mu)\leq 2\RR_1(X,\omega,\mu)$. 
\\
\ \par
Puisque $\max_X\varphi-m=0$, d'apr\`es (2.2), on a
\begin{eqnarray*}
\left|\int\varphi\omega^k\right| & \leq & \int|\varphi-m|\omega^k+m=
\int(m-\varphi)\omega^k +m\\
& \leq & \RR^*_1(X,\omega)+\RR_1(X,\omega,\mu).
\end{eqnarray*}
Ceci implique que $\RR_3(X,\omega,\mu)\leq \RR_1(X,\omega,\mu)+ \RR^*_1(X,\omega)$.
\\
\ \par
Supposons que $\mu$ soit $(c,\alpha)$-mod\'er\'ee.
Posons $\psi:=r^{-1}(\varphi-m)$.
On a $\ddc \psi\geq -\omega$ et $\max_X\psi=0$. Puisque
$\int\varphi\d\mu=0$, on a
$m\geq 0$ et donc $\psi\leq r^{-1}\varphi$.
La mesure $\mu$ \'etant $(c,\alpha)$-mod\'er\'ee, on a
\begin{eqnarray*}
\mu(\varphi<-t) & = & \mu(r^{-1}\varphi<-r^{-1}t)\leq
\mu(\psi<-r^{-1}t)\\
& \leq & \exp(-\alpha
r^{-1}t) \int\exp(-\alpha \psi)\d\mu \\
& \leq & c\exp(-\alpha r^{-1}t).
\end{eqnarray*}
Donc $\Delta(X,\omega,\mu,t)\leq c\exp(-\alpha
r^{-1}t)$. 
\end{preuve}
\par\
\\
{\bf 2.3}. {\it Image directe d'un courant.}
\\
\
\par
Soit $\pi$ une application holomorphe surjective de
$X$ dans $X'$.
Si $S$ est un courant de bidimension $(r,r)$ de $X$, avec $0\leq
r\leq k, k'$, le courant $\pi_*(S)$ est
d\'efini par
$$\langle \pi_*(S),\psi \rangle := \langle S, \pi^*(\psi)\rangle$$ 
pour toute forme lisse $\psi$ de bidegr\'e $(r,r)$ de $X'$.
\par
Si $S$ est une forme \`a coefficients dans $\Lone(X)$,
$\pi_*(S)$ l'est aussi. En effet, on peut supposer que $S$ est
r\'eelle positive; une forme positive
est \`a coefficients dans $\Lone(X)$ si et
seulement si, elle d\'efinit un courant positif de masse finie
dont les coefficients dans une carte
sont des fonctions mesurables. Les
coefficients de $\pi_*(S)$ sont obtenus par int\'egration sur les
fibres qui sont presque partout de m\^eme dimension.
\par
Si les fibres g\'en\'eriques de $\pi$ sont
discr\`etes, en g\'en\'eral,
l'image par $\pi_*$ d'une fonction q.p.s.h. n'est
pas q.p.s.h., elle est diff\'erence de telles fonctions, d'o\`u
l'introduction des fonctions d.s.h.  
\begin{proposition} 
Soient $(X,\omega)$ et $(X',\omega')$ des vari\'et\'es k\"ahl\'eriennes
compactes de dimension $k$. Supposons que
pour un point g\'en\'erique $x'\in X'$
la fibre $\pi^{-1}(x')$ soit finie et non vide. Alors 
\begin{enumerate}
\item[{\rm (a)}] L'image de $\DSH(X)$ par $\pi_*$ est contenue dans
$\DSH(X')$.
\item[{\rm (b)}] L'image de la famille $\{\varphi \mbox{ q.p.s.h.
sur } X,\
\ddc\varphi\geq -\omega,\int_X\varphi\omega^{k}=0\}$ par
$\pi_*$ est relativement compacte dans $\Lp(X')$ pour tout $p\geq 1$.
\end{enumerate}
\end{proposition}
\begin{preuve}
(a) Soit $\psi$ une fonction q.p.s.h. dans $X$. 
Notons $I(\pi)$ l'ensemble des points $x'\in X'$ tels que la
fibre $\pi^{-1}(x')$ ne soit pas finie.
Il est clair que la fonction $\pi_*(\psi)$ est d\'efinie hors
de l'ensemble $I(\pi)$ qui est pluripolaire
(\voir proposition A.1).
Rappelons que pour $x'\in X'\setminus I(\pi)$ on a
$$\pi_*(\psi)(x')=\sum_{\pi(x_i)=x'}\psi(x_i).$$
\par
Soient $T^\pm$
des courants positifs ferm\'es de bidegr\'e $(1,1)$ tels que
$\ddc\psi=T^+-T^-$. On a
$\ddc\pi_*\psi=\pi_*(T^+)-\pi_*(T^-)$. Puisque $\pi_*(T^\pm)$ sont
des courants positifs ferm\'es de bidegr\'e $(1,1)$, la fonction
$\pi_*(\psi)$ est d.s.h. 
\\
\ \par
(b) 
Soit $(\psi_n)$ une suite de fonctions q.p.s.h. v\'erifiant
$\ddc\psi_n\geq -\omega$ et $\int_{X}\psi_n\omega^{k}=0$.
D'apr\`es la proposition 2.1, cette suite de fonctions
est born\'ee dans $\Lp(X)$ et est born\'ee sup\'erieurement.
Il faut montrer que la suite des $\pi_*(\psi_n)$ est born\'ee dans 
$\Lp(X')$. Quitte \`a extraire une sous-suite, on peut
supposer que $\psi_n$
converge ponctuellement vers une fonction $\psi$
hors d'un ensemble de mesure nulle.
\par
Les fonctions $\pi_*(\psi_n)$ sont dans $\Lp(X')$
et v\'erifient 
$\ddc \pi_*(\psi_n)\geq -\pi_*(\omega)$.
D'apr\`es la proposition 2.2, il existe une fonction $\varphi$
q.p.s.h. sur $X'$ v\'erifiant $\ddc \varphi\geq \pi_*(\omega)
-c\omega'$ o\`u $c\geq 0$ est une constante. Posons
$\varphi_n:=\pi_*(\psi_n)+\varphi$.
On a $\ddc\varphi_n\geq -c\omega'$. 
Montrons que
$(\varphi_n)$ est born\'ee dans $\Lp(X')$. 
La suite $(\varphi_n)$ \'etant born\'ee sup\'erieurement,
il suffit d'observer que
$\psi_n$ converge presque partout vers $\psi$ et donc,
aucune sous-suite de $(\varphi_n)$ ne converge
uniform\'ement vers $-\infty$. Il en r\'esulte que $\pi_*(\psi_n)$ est born\'ee
dans $\Lp(X')$. On peut en extraire des sous-suites convergentes.
\end{preuve}
\par\ 
\\
{\bf 2.4}. {\it Image r\'eciproque d'un courant.}
\\
\
\par
Soit $\pi$ une application holomorphe surjective de
$X$ dans $X'$.
L'image
r\'eciproque $\pi^*(\alpha)$ d'un courant $\alpha$ sur $X'$
est d\'efinie lorsque $\pi$ est une submersion.
Si $\pi$ n'est pas une submersion, on d\'efinit $\pi^*$ dans
les cas sp\'eciaux suivants.
Lorsque $\gamma$ est une
$(p,q)$ forme \`a coefficients dans $\Linfty(X')$, la forme
$\pi^*(\gamma)$
est bien d\'efinie. C'est une forme \`a coefficients dans
$\Linfty(X)$. Lorsque $\psi$ est une fonction d.s.h. sur $X'$, si
$\pi$ est surjective,
$\psi\circ \pi$ est aussi une fonction d.s.h. sur $X$;
en particulier, on a $\psi\circ \pi\in\Lp(X)$ pour tout $p\geq 1$.
On peut d\'efinir $\pi^*(\psi\gamma):=(\psi\circ\pi)
\pi^*(\gamma)$. 
\begin{proposition}
Soit $\pi:X\longrightarrow X'$ une application holomorphe
surjective. L'image
par $\pi^*$ de la famille
$$\left\{\varphi \mbox{ q.p.s.h. sur } X',\
\ddc\varphi\geq -\omega', \int_{X'}\varphi{\omega'}^{k'}=0\right\}$$
est relativement compacte dans $\Lp(X)$ pour tout $p\geq 1$.
\end{proposition}
\begin{preuve}
Soit $(\psi_n)$ une suite born\'ee dans $\Lp(X')$
telle que  $\ddc\psi_n\geq -\omega'$.
Il suffit de montrer que la suite $(\pi^*(\psi_n))$ est born\'ee
dans $\Lp(X)$.  
D'apr\`es la proposition 2.1, les fonctions $\psi_n$ sont born\'ees
sup\'erieurement par une m\^eme constante, il en est donc de
m\^eme   
des fonctions $\pi^*(\psi_n)$. De plus,
$\ddc \pi^*(\psi_n)\geq -\pi^*(\omega')\geq
-c\omega_1$ o\`u $c:=\|\pi\|_{{\cal C}^2}$. 
Quitte \`a extraire une
sous-suite, on peut supposer que $\pi^*(\psi_n)$ tend vers une
fonction $\varphi$ avec $\ddc\varphi\geq -c\omega$ ou
sinon $\pi^*(\psi_n)$ tend uniform\'ement vers $-\infty$.
Le deuxi\`eme
cas ne se produit pas car la suite $(\psi_n)$ est born\'ee dans
$\Lp(X')$. On en
d\'eduit que $\pi^*(\psi_n)$
tend vers $\varphi$ dans $\Lp(X)$. 
\end{preuve}
\par
D\'efinissons $\pi^*$ pour les courants de bidegr\'e $(1,1)$
lorsque 
 $\pi$ est surjective et $k\geq k'$.
Soit $S$ un courant positif ferm\'e de bidegr\'e $(1,1)$ dans
$X'$, et soit $\alpha$ une forme lisse de bidegr\'e $(1,1)$ de
$X'$ cohomologue \`a $S$. Il existe une fonction
q.p.s.h. $\psi$ dans $X'$ telle que $\ddc\psi=S-\alpha$.
Posons $\varphi:=\psi\circ \pi$. On a $\ddc\varphi\geq
-\pi^*(\alpha)$. C'est donc une fonction q.p.s.h. On pose
$\pi^*(S):=\ddc\varphi +\pi^*(\alpha)$. C'est
un courant positif ferm\'e de masse finie. Cette op\'eration est
continue et ind\'ependante du choix de $\alpha$ et de $\psi$
\cite{Meo}. 
On a aussi $\cl(\pi^*(S))=\cl(\pi^*(\alpha))$.
Si $\psi$ est une fonction d.s.h. sur $X'$ avec
$\ddc\psi=T^+-T^-$, on a 
$\ddc\pi^*(\psi)=\pi^*(T^+)-\pi^*(T^-)$.
\\
\ \par
Notons $I$ l'ensemble des points $x'\in X'$ tels que
$\dim\pi^{-1}(x')> k-k'$. C'est un sous-ensemble analytique
de codimension au moins $2$ de $X'$ car $\dim\pi^{-1}(I)$ est au
plus \'egale \`a $k-1$. Soit $H$ un sous-ensemble
analytique de dimension pure $l$ de $X'$. Supposons que $\dim
\pi^{-1}(H\cap I)<l+k-k'$. Dans ce cas, on peut d\'efinir
$\pi^*[H]$ comme le courant d'int\'egration sur $\pi^{-1}(H)$,
il est alors de m\^eme bidegr\'e que $[H]$.
\\
\ \par
Lorque $X'$ est une vari\'et\'e projective et $T$ est un
courant positif ferm\'e de bidegr\'e $(r,r)$, on peut d\'efinir
"la partie principale" de l'image r\'eciproque de $T$. On a le lemme
suivant.
\begin{lemme}{\bf \cite{DinhSibony3}} Soit $X'$ une vari\'et\'e
projective munie d'une forme de K\"ahler $\omega'$. 
Il existe une constante $c>0$, qui ne d\'epend que de
$(X',\omega')$, telle que pour tout courant positif ferm\'e $T$ de bidegr\'e
$(r,r)$ sur $X'$ on puisse trouver une suite de courants positifs
ferm\'es lisses $(T_m)_{m\geq 1}$, de bidegr\'e $(r,r)$,
v\'erifiant les propri\'et\'es suivantes
\begin{enumerate}
\item La suite $(T_m)$ converge vers un courant positif ferm\'e $T'$.
\item $T'\geq T$, c'est-\`a-dire que le courant $T'-T$ est positif.
\item On a $\|T_m\|\leq c\|T\|$ pour tout $m\geq 1$.
\end{enumerate}
\end{lemme}
\begin{remarques}\rm
Observons que l'ensemble des classes $\cl(S)$ des courants
$S$ positifs ferm\'es de bidegr\'e $(r,r)$
de masse 1, est born\'e dans $\H^{r,r}(X',\R)$. 
Il existe donc une
constante $\alpha_X>0$ ind\'ependante de $S$ telle que la classe de
$\alpha_X{\omega'}^r -S$ soit
repr\'esent\'ee par une forme lisse positive.
On dira que $S$ est {\it cohomologiquement domin\'e}
par $\alpha_X{\omega'}^r$. 
Les courants $T_m$ dans le lemme 2.8 sont cohomologiquement domin\'es
par $c_X\|T\|{\omega'}^r$ o\`u $c_X:=c\alpha_X$.
\par
Soit $X'$ une vari\'et\'e k\"ahl\'erienne compacte 
pour laquelle il existe une projection holomorphe surjective de $X$ sur une 
vari\'et\'e k\"ahl\'erienne compacte homog\`ene de m\^eme dimension.
Le lemme 2.8 reste alors valable pour tout courant $T$ qui ne charge pas 
les sous-ensembles analytiques de $X'$.
\end{remarques}
\par
Soit $\Omega\subset X$ l'ouvert, Zariski dense, o\`u $\pi$ est une
submersion locale. Le courant $(\pi_{|\Omega})^*(T)$ est bien
d\'efini, positif, ferm\'e sur $\Omega$. Le lemme 2.8 permet de
montrer que  $(\pi_{|\Omega})^*(T)$ est
de masse finie et donc, d'apr\`es Skoda \cite{Skoda},
son prolongement trivial $\widetilde{(\pi_{|\Omega})^*(T)}$ est
positif ferm\'e sur $X$.
\begin{proposition} {\bf \cite{DinhSibony3}}
Soit $T$ un courant positif ferm\'e sur une vari\'et\'e projective 
$X'$. Alors
$(\pi_{|\Omega})^*(T)$ est de masse finie et l'op\'erateur
$T\mapsto \widetilde{(\pi_{|\Omega})^*(T)}$ est semi-continu
inf\'erieurement. Plus pr\'ecis\'ement, si $T_n\rightarrow T$, tout courant
adh\'erent \`a la suite 
$\widetilde{(\pi_{|\Omega})^*T_n}$ est sup\'erieur ou \'egal \` a 
$\widetilde{(\pi_{|\Omega})^*T}$.  
\end{proposition}
\sectionbis{Transformations m\'eromorphes}
Dans ce paragraphe, nous d\'efinissons les
op\'erations: composition, produit, intersection, 
sur les transformations et les correspondances
m\'eromorphes. Nous \'etudions l'effet de ces
op\'erations sur les degr\'es
interm\'ediaires.
\par
\
\\
{\bf 3.1}. {\it D\'efinitions.}
\\
\ 
\par
Notons $\pi_1:X_1\times X_2\longrightarrow X_1$  et $\pi_2:X_1\times
X_2 \longrightarrow X_2$, les projections
canoniques de $X_1\times X_2$
sur $X_1$ et $X_2$.
On appelle {\it $m$-cha\^{\i}ne holomorphe positive}
de $X_1\times X_2$ toute
combinaison finie $\Gamma:=\sum \Gamma_j$ o\`u les
$\Gamma_j$ sont des
sous-ensembles analytiques irr\'eductibles
de dimension $m$ de $X_1\times X_2$.
Les $\Gamma_j$ ne sont pas n\'ecessairement distincts.
D'apr\`es un th\'eor\`eme de Lelong \cite{Lelong},
l'int\'egration sur la
partie lisse d'une $m$-cha\^{\i}ne
holomorphe positive $\Gamma$ d\'efinit
un courant positif ferm\'e $[\Gamma]$ de bidimension
$(m,m)$ sur $X_1\times X_2$. 
Notons 
$\overline \Gamma$ l'image de $\Gamma$ par l'application
$(x_1,x_2)\mapsto(x_2,x_1)$.
\par
Soit $l$ un entier naturel, $k_1-k_2\leq l<k_1$.
On appelle {\it transformation m\'eromorphe} $F$ de $X_1$ dans $X_2$
toute $(k_2+l)$-cha\^{\i}ne
holomorphe positive $\Gamma=\sum \Gamma_j$
de $X_1\times X_2$ telle que la restriction
de $\pi_i$ \`a chaque composante irr\'eductible $\Gamma_j$
soit surjective, $i=1,2$. 
On dira que $\Gamma$ est {\it le graphe} de $F$
et que $\codim(F):=l$
est {\it la codimension}
de $F$.
La transformation m\'eromorphe $\overline F$ de $X_2$ dans $X_1$
associ\'ee \`a la
$(k_2+l)$-cha\^{\i}ne
holomorphe $\overline\Gamma$ est appel\'ee {\it
transformation m\'eromorphe adjointe} de $F$, elle est de codimension $k_2-k_1+l$.
\par
Posons $F:=\pi_2\circ (\pi_{1|\Gamma})^{-1}$ et
$F^{-1}:=\pi_1\circ(\pi_{2|\Gamma})^{-1}$. Ces "applications"
sont d\'efinies sur les sous ensembles de $X_1$ et $X_2$.
{\it La fibre}
$F^{-1}(x_2)$ de $x_2\in X_2$ est g\'en\'eriquement un
sous-ensemble analytique de dimension $l$ de $X_1$.
\par
Notons $I_i(F):=\{x\in X_i,\ \dim\pi_i^{-1}(x)>k_2+l-k_i\}$. C'est
un sous-ensemble analytique de codimension au moins 2 de $X_i$.
En effet, s'il \'etait de codimension 1, le graphe contiendrait
un ouvert sur lequel $\pi_i$ ne serait pas une submersion.
On dira que $I_1(F)$ (resp. $I_2(F)$) est le {\it premier} (resp.
{\it deuxi\`eme}) {\it ensemble d'ind\'etermination} de $F$.
\\
\ \par
D\'efinissons les op\'erateurs $F^*$ et $F_*$.
Soit $S$ un courant de
bidegr\'e $(r,r)$ sur $X_2$, $k_2+l-k_1\leq r\leq k_2$. On
d\'efinit $F^*(S):=(\pi_1)_*(\pi_2^*(S)\wedge [\Gamma])$. 
C'est un courant de bidimension $(k_2+l-r,k_2+l-r)$
port\'e par $F^{-1}(\supp(S))$.
Cet
op\'erateur est d\'efini dans les deux espaces suivants:
\begin{enumerate}
\item  L'espace des formes lisses.
\item L'espace engendr\'e par les courants $[H]$
o\`u $H$ est un sous-ensemble
analytique de dimension pure $k_2-r$ de $X_2$ v\'erifiant
$\dim(\pi_{2|\Gamma})^{-1}(H\cap I_2(F))\leq k_2-r+l-1$. 
\end{enumerate}
\par
Consid\'erons une r\'esolution des singularit\'es
$\pi:\widetilde\Gamma\longrightarrow \Gamma$ \cite{Hironaka}. Posons
$\tau_i:=\pi_i\circ\pi$. On a 
$F=\tau_2\circ\tau_1^{-1}$ et
$F^{-1}=\tau_1\circ\tau_2^{-1}$. On a aussi
$F^*=(\tau_1)_*(\tau_2)^*$ lorsque cet op\'erateur agit sur
les courants d\'ecrits ci-dessus.
\par
Posons $F_*:=\overline F^*=(\tau_2)_*(\tau_1)^*$.
Observons que lorsque $\codim (F)=0$, l'op\'erateur $F_*$ 
agit aussi sur les fonctions q.p.s.h. et
donc sur les courants positifs ferm\'es
de bidegr\'e $(1,1)$. Dans ce cas, on peut aussi utiliser, \`a la place 
des applications $\tau_i$,
les fonctions q.p.s.h.
sur les sous-ensembles analytiques.
\par
Lorsque $S$ est lisse, le courant $F^*(S)$ est une forme \`a
coefficients dans $\Lone(X)$. Pour $x_2\in X_2\setminus I_2(F)$ le
courant $F^*(\delta_{x_2})$ est un courant d'int\'egration sur
une $l$-cha\^{\i}ne positive port\'ee par 
$F^{-1}(x_2)$.
\par
Pour tout $s$,
$k_2-k_1+l\leq s\leq k_2$, on appelle {\it
degr\'e interm\'ediaire d'ordre $s$} de $F$ le nombre
\begin{eqnarray}
\lambda_s(F) & := & \int_{X_1} F^*(\omega_2^s)
\wedge \omega_1^{k_2+l-s}=\int_{X_2}
\omega_2^s \wedge F_*(\omega_1^{k_2+l-s})\nonumber\\
& = & \int_\Gamma \pi_1^*(\omega_1^{k_2+l-s})\wedge\pi_2^*(\omega_2^s)
\end{eqnarray}
Par continuit\'e, la masse du courant $F^*(\delta_{x_2})$, 
qui se calcule cohomologiquement,
ne d\'epend pas de $x_2$,  pour
$x_2$ g\'en\'erique dans $X_2$. On en d\'eduit que cette masse est
\'egale au
dernier degr\'e interm\'ediaire $\lambda_{k_2}(F)$ de $F$.
\par
Finalement, on dira qu'un point $(x_1,x_2)\in\Gamma$ est
{\it g\'en\'erique} si la restriction de
$\pi_{i|\Gamma}$ \`a un voisinage de $(x_1,x_2)$ est une
submersion pour $i=1,2$. Notons $\Gen(\Gamma)$ l'ensemble de ces
points. C'est un ouvert de Zariski dense de $\Gamma$.  
\par
\
\\
{\bf 3.2}. {\it Composition de transformations m\'eromorphes.} 
\\
\
\par
Soit $F$ une transformation m\'eromorphe de codimension $l$
de $X_1$ dans $X_2$ comme ci-dessus. 
Soit $F'$ une autre transformation m\'eromorphe de $X_2$ dans $X_3$
associ\'ee \`a une $(k_3+l')$-cha\^{\i}ne holomorphe
$\Gamma'=\sum\Gamma_j'$ de
$X_2\times X_3$. Supposons que $l+l'<k_1$.
\par
Consid\'erons d'abord le cas o\`u $\Gamma$ et $\Gamma'$ sont
irr\'eductibles. D\'efinissons la compos\'ee $\Gamma'\circ\Gamma$ des graphes.
Notons $\pi_1$, $\pi_2$ les projections de $X_1\times X_2$
sur $X_1$ et $X_2$
et $\pi_2'$, $\pi_3'$ les projections de $X_2\times X_3$ sur $X_2$ et
$X_3$. Soient $x_2\in X_2$ un point g\'en\'erique et $(x_1,x_2)$, $(x_2,x_3)$ 
des points g\'en\'eriques de $\Gamma\cap\pi_2^{-1}(z_2)$ et 
$\Gamma'\cap {\pi_2'}^{-1}(x_2)$. On peut supposer que $(x_1,x_2)\in \Gen(\Gamma)$
et $(x_2,x_3)\in\Gen(\Gamma')$. Soient $U\subset \Gen(\Gamma)$ et 
$U'\subset\Gen(\Gamma')$ des petits voisinages de $(x_1,x_2)$ dans $\Gen(\Gamma)$ et
de $(x_2,x_3)$ dans $\Gen(\Gamma')$. Par d\'efinition de $\Gen(\Gamma)$ 
et $\Gen(\Gamma')$, 
on peut supposer que $U$ et $U'$ admettent des structures produit
$U\simeq W_1\times V_2$ et $U'\simeq 
V_2\times W_3$ o\`u $V_2$ d\'esigne un voisinage de
$x_2$ dans $X_2$. Les projections $\pi_2$, $\pi_2'$ 
de $U$ et $U'$ sur $X_2$ co\"{\i}ncident avec les 
projections des produits sur le facteur $V_2$. Les projections de $U$ sur $X_1$
et de $U'$ sur $X_3$ correspondent \`a des applications holomorphes 
$\tau:U\longrightarrow X_1$ et  
$\tau':U'\longrightarrow X_3$.
\par
Le mod\`ele local de $\Gamma'\circ\Gamma$ est l'image de $W_1\times V_2\times W_3$
dans $X_1\times X_3$ par l'application $(x_1,x_2,x_3)\mapsto 
(\tau(x_1,x_2),\tau'(x_2,x_3))$. On suppose que cette image est 
de dimension $k_3+l+l'$. On dira alors que $\Gamma$ et $\Gamma'$ {\it 
se composent correctement}.
Si cette hypoth\`ese est v\'erifi\'ee dans les petits ouverts, 
elle est v\'erifi\'ee en tout point g\'en\'erique. 
\par
Le graphe $\Gamma'\circ\Gamma$ de $F'\circ F$ est
alors l'adh\'erence de
l'ensemble des $(x_1,x_3)\in X_1\times X_3$ pour
lesquels il existe $x_2\in X_2$
avec $(x_1,x_2)\in\Gen(\Gamma)$ et
$(x_2,x_3)\in \Gen(\Gamma')$ tel qu'aux voisinages de ces points $\Gamma$ 
et $\Gamma'$ se composent correctement. Le point $(x_1,x_3)$ est compt\'e
avec la multiplicit\'e $m$ si $m$ est le nombre de $x_2$
pour lesquels $x_1$, $x_2$, $x_3$ v\'erifient
la propri\'et\'e ci-dessus. 
Puisque $\Gamma$ et $\Gamma'$ se composent correctement,
 $m$ est fini.
\par
Dans le cas o\`u $\Gamma$ et $\Gamma'$ ne sont pas
irr\'eductibles, on pose $\Gamma'\circ\Gamma:=\sum
\Gamma_j'\circ\Gamma_i$ en supposant que $\Gamma_i$ et $\Gamma_j'$ se composent 
correctement pour tout $i$, $j$.  
Observons que $\Gamma'\circ\Gamma$ est une
$(k_3+l+l')$-cha\^{\i}ne holomorphe
et qu'on a $\codim(F'\circ F)=\codim(F)+\codim(F')=l+l'$.
Observons aussi que
les transformations de codimension $0$ entre vari\'et\'es de
m\^eme dimension (\cad les correspondances)
se composent toujours correctement.
\begin{proposition} Soit $X_2$ une vari\'et\'e
projective de dimension $k_2$ munie d'une forme de K\"ahler
$\omega_2$.
Alors il existe une constante $c>0$, qui ne d\'epend que de
$(X_2,\omega_2)$,
telle que 
$$\lambda_s(F'\circ F)\leq
c\lambda_{k_2-k_3+s-l'}(F)\lambda_s(F')$$
pour tout $s$ avec
$k_3-k_1+l+l'\leq s\leq k_3$ et pour toutes les transformations m\'eromorphes
$F$ de $(X_1,\omega_1)$ dans $(X_2,\omega_2)$ et
$F'$ de $(X_2,\omega_2)$ dans $(X_3,\omega_3)$.
\end{proposition}
\begin{preuve}
Observons que dans (3.1) les formes \'etant lisses ou
\`a coefficients dans $\Lone$, les int\'egrales peuvent ne porter
que 
sur des ouverts de mesure totale. Posons $S:=(F')^*(\omega_3^s)$.
C'est un courant positif ferm\'e 
\`a coefficients dans $\Lone(X_2)$ de
bidegr\'e $(r,r)$ sur $X_2$ o\`u $r:=k_2-k_3+s-l'$.
D'apr\`es le lemme 2.8 et la remarque 2.9, il existe des courants
$S_m$ positifs ferm\'es lisses cohomologiquement domin\'es par
$c_{X_2} \|S\|\omega_2^r$ qui convergent vers un courant $S'$
v\'erifiant $S'\geq S$. 
On a 
\begin{eqnarray*}
\lambda_s(F'\circ F) & = & \int_{\Gen(\Gamma)}
(\pi_1)^*(\omega_1^{k_3+l+l'-s}) \wedge
(\pi_2)^*(S) \\
& \leq & \int_{\Gen(\Gamma)}
(\pi_1)^*(\omega_1^{k_3+l+l'-s}) \wedge
(\pi_2)^*(S') \\
& \leq & \lim_{m\rightarrow\infty} \int_\Gamma
(\pi_1)^*(\omega_1^{k_3+l+l'-s})\wedge (\pi_2)^*(S_m)\\
& \leq & c_{X_2}\lambda_s(F') \int_\Gamma
(\pi_1)^*(\omega_1^{k_3+l+l'-s})\wedge
(\pi_2)^*(\omega^r)\\
& = & c_{X_2}\lambda_{k_2-k_3+s-l'}(F)\lambda_s(F').
\end{eqnarray*}
La premi\`ere \'egalit\'e r\'esulte de la description locale de
$\Gamma'\circ\Gamma$, la lin\'earit\'e 
permet ensuite d'utiliser des partitions de l'unit\'e. 
Pour la deuxi\`eme in\'egalit\'e, on utilise une suite 
exhaustive de compacts de $\Gen(\Gamma)$. 
\end{preuve}
\
\\
{\bf 3.3}. {\it Produit et intersection de transformations m\'eromorphes.} 
\\
\
\par
Consid\'erons deux transformations m\'eromorphes
$F_1:X\longrightarrow X_1$ et $F_2:X\longrightarrow X_2$
de codimensions respectives $l_1$ et $l_2$.
On suppose que $l_1+l_2\geq k$ et qu'il existe des ouverts,
Zariski denses, $\Omega_1\subset X_1$ et $\Omega_2\subset X_2$
tels que $F_1^{-1}(x_1)\cap F_2^{-1}(x_2)$
soit de dimension pure $l_1+l_2-k$ pour tout
$(x_1,x_2)\in\Omega_1\times \Omega_2$.
\par
D\'efinissons {\it le produit } $F_1\bullet F_2$.
C'est une transformation m\'eromorphe de $X$ dans $X_1\times X_2$
de codimension $l_1+l_2-k$
dont nous allons d\'ecrire le graphe. Notons
$\Gamma^1=\sum\Gamma_i^1$ et $\Gamma^2=\sum\Gamma_j^2$ les
graphes de $F_1$ et $F_2$. Consid\'erons d'abord le cas o\`u
$\Gamma^1$ et $\Gamma^2$ sont irr\'eductibles. Le graphe
$\Gamma^1\bullet \Gamma^2$ de $F_1\bullet F_2$ est alors l'adh\'erence de
l'ensemble des $(x,x_1,x_2)\in X\times\Omega_1\times\Omega_2$
avec $x\in F_1^{-1}(x_1)\cap F_2^{-1}(x_2)$. Dans le cas
g\'en\'eral, on pose $\Gamma^1\bullet\Gamma^2:=\sum
\Gamma^1_i\bullet\Gamma^2_j$. 
\par
On munit $X_1\times X_2$ de la forme de K\"ahler
$\omega_{12}:=
c_{12}(\pi_1^*(\omega_1)+\pi_2^*(\omega_2))$ o\`u $\pi_1$, $\pi_2$
sont des projections sur $X_1$ et $X_2$ et
$c_{12}^{-k_1-k_2}:={k_1+k_2\choose k_1}$. Le choix de $c_{12}$
implique que $\int_{X_1\times
X_2}(\omega_{12})^{k_1+k_2}=1$. 
\begin{proposition} Soit $X$ une
vari\'et\'e projective munie d'une forme de K\"ahler $\omega$.
Il existe une constante $c>0$ qui
ne d\'epend que de $(X,\omega)$ telle que 
pour tout $s$ v\'erifiant
$k_1+k_2-2k+l_1+l_2\leq s\leq k_1+k_2$ on ait
$$\lambda_s(F_1\bullet F_2)\leq
cc_{12}^s\sum {s\choose s_1}
\lambda_{s_1}(F_1)\lambda_{s_2}(F_2)$$
avec $k_1-k+l_1\leq s_1\leq k_1$, $k_2-k+l_2\leq s_2\leq k_2$
et $s_1+s_2=s$.
\end{proposition}
\begin{preuve} Posons $F:=F_1\bullet F_2$. On a pour $r:=k_1+k_2+l_1+l_2-k-s$
\begin{eqnarray}
\lambda_s(F) & = & \int_X F^*(\omega_{12}^s)\wedge
\omega^r\nonumber \\
& = & c_{12}^s\sum_{s_1,s_2} {s \choose s_1}
\int_X (F_1)^*(\omega_1^{s_1}) \wedge
(F_2)^*(\omega_2^{s_2})\wedge \omega^r
\end{eqnarray}
La derni\`ere \'egalit\'e se v\'erifie sur des mod\`eles locaux \`a l'aide 
de partitions de l'unit\'e.
Estimons l'int\'egrale dans (3.2). Posons
$S_i:=(F_i)^*(\omega_i^{s_i})$. 
D'apr\`es le lemme 2.8,
il existe des courants lisses $S_{i,m}$ cohomologiquement
domin\'es par $c_X\|S_i\|\omega^{s_i}$ qui tendent vers un
courant $S_i'\geq S_i$. L'int\'egrale dans (3.2) est donc
major\'ee par
$$c_X^2\|S_1\|\|S_2\|\int_X\omega^k=
c_X^2\lambda_{s_1}(F_1)\lambda_{s_2}(F_2).$$  
\end{preuve}
\par
Soient $G_1:X_1\longrightarrow X$ et $G_2:X_2\longrightarrow X$
deux transformations m\'eromorphes. On d\'efinit {\it
l'intersection} $G_1\cap G_2$ de $G_1$ et $G_2$ comme l'adjoint
$\overline{\overline G_1\bullet \overline G_2}$ du produit
$\overline G_1\bullet \overline G_2$ lorsque ce produit est bien
d\'efini. C'est une transformation m\'eromorphe de $X_1\times
X_2$ dans $X$. 
Pour $(x_1,x_2)$ g\'en\'erique,
$(G_1\cap G_2)(x_1,x_2)$ est l'intersection de $G_1(x_1)$
et $G_2(x_2)$.
\begin{remarque} \rm
Pour la validit\'e des relations de la proposition 3.1 (resp.
proposition 3.2), il suffit de supposer l'existence d'une
application holomorphe surjective $\pi$ de $X_2$ (resp. $X$) dans
une vari\'et\'e k\"ahl\'erienne compacte homog\`ene de m\^eme dimension
(\voir remarques 2.9). 
\end{remarque}
\par\ 
\\
{\bf 3.4}. {\it Familles m\'eromorphes
adapt\'ees de sous-ensembles analytiques.}
\\
\ \par
Soit $F'$ une
transformation m\'eromorphe de codimension $l'$
de $X_2$ dans $X_3$ dont le graphe est irr\'eductible.
Puisque $F'$ est surjective , la r\'eunion de ses fibres
$\H_{x_3}:={F'}^{-1}(x_3)$ est \'egale \`a $X_2$. On dira que
$\H=(\H_{x_3})$ est {\it une famille adapt\'ee d'ensembles 
analytiques de dimension $l'$}. 
Si $x_3$ n'appartient pas au deuxi\`eme
ensemble d'ind\'etermination $I_2(F')$ de $F'$,
${\cal H}_{x_3}$ est de dimension $l'=\codim(F')$. 
D'apr\`es le
th\'eor\`eme de Bertini \cite[p.141]{Shafarevich},
pour $x_3$ g\'en\'erique, les composantes
de $\H_{x_3}$ sont de multiplicit\'e 1. On a donc
$[\H_{x_3}]=(F')^*(\delta_{x_3})$. 
\par
Comme pr\'ec\'edemment,
soit $F$ une transformation m\'eromorphe de codimension $l$
de $X_1$ dans $X_2$.
Supposons que $l+l'<k_1$. 
On dira que ${\cal
H}$ est $F$-r\'eguli\`ere si pour $x_3\in X_3$ g\'en\'erique, la
dimension de $F^{-1}({\cal H}_{x_3}\cap I_2(F))$ est strictement
plus petite que $(l+l')$. Pour un tel $x_3$, le courant
$F^*{F'}^*(\delta_{x_3})$ est bien d\'efini.
C'est un courant
d'int\'egration sur une cha\^{\i}ne holomorphe de dimension
$(l+l')$. La famille $\H$ est dite {\it
r\'eguli\`ere}, si elle est $F$-r\'eguli\`ere pour toute
transformation m\'eromorphe $F$ d'une vari\'et\'e $X_1$ dans $X_2$. 
Les familles adapt\'ees de sous-vari\'et\'es
associ\'ees aux transformations m\'eromorphes
$\Psi_1$, $\Psi_2$, $\Psi_3$ et $F_{l,n}$ que nous allons d\'ecrire au 
paragraphe 3.6 sont r\'eguli\`eres.
\par
Les {\it intersections de familles
adapt\'ees} de sous-ensembles analytiques sont d\'efinies
comme \'etant associ\'ees aux produits de transformations
m\'eromorphes.  
\par\ 
\\
{\bf 3.5}. {\it Correspondances m\'eromorphes.}
\\
\ \par
Supposons que $\dim X_1=\dim
X_2=k$. Une {\it correspondance
(m\'eromorphe)} de $X_1$ dans $X_2$
est une transformation m\'eromorphe $f$ de
codimension $0$ de $X_1$ dans $X_2$.
Notons $\Gamma=\sum \Gamma_i$ le graphe de $f$. 
La correspondance $\overline f$ de $X_2$ dans $X_1$ associ\'ee
\`a $\overline \Gamma$ est appel\'ee
{\it correspondance adjointe} de $f$.
\par
Lorsque la restriction
de $\pi_1$ \`a $\Gamma$
est injective hors d'un sous-ensemble analytique, on
dira que $f$ est une {\it application m\'eromorphe surjective}.
On dit que
$f$ est {\it bim\'eromorphe} si $f$ et son adjoint
$\overline f$ sont des
applications m\'eromorphes surjectives.
\par
Notons $\DD(X_1,\omega_1)$
l'ensemble des fonctions $\psi$ d.s.h. sur $X_1$
pour lesquelles il existe des courants $T^\pm$ positifs ferm\'es de
bidegr\'e $(1,1)$ et de masse 1 v\'erifiant  $\ddc\psi=T^+-T^-$.
Posons
\begin{eqnarray}
\AA(f) & := & \sup_\psi
\left\{\left|\int_{X_2}f_*(\psi)\omega_2^{k}\right|,\ 
\psi\in \DD(X_1,\omega_1), \int \psi\omega_1^{k}=0\right\}
\end{eqnarray}
D'apr\`es les propositions 2.6 et 2.7, cette constante est finie.
Elle mesure combien $f_*$ perturbe la normalisation
$\int\psi\omega_1^{k}=0$. 
\\
\ \par
Consid\'erons le cas o\`u $X_1=X_2=X$.
On notera $f^n$ la correspondance $f\circ\cdots \circ f$ ($n$ fois).
Pour tout $0\leq s\leq k$, on d\'efinit le degr\'e dynamique d'ordre
$s$ de $f$ par la formule suivante:
\begin{eqnarray}
d_s(f) & := & \limsup_{n\rightarrow \infty}
\big[\lambda_s(f^n)\big]^{1/n}
\end{eqnarray}
Le dernier degr\'e dynamique $d_k(f)$ est \'egal au nombre
d'\'el\'ements de la fibre $f^{-1}(z)$ pour un point $z$
g\'en\'erique (ce nombre ne d\'epend pas de $z$). C'est
{\it le degr\'e topologique} de $f$, on le note par $d_t(f)$. On a
aussi 
$d_0(f)=d_t(\overline f)$. 
Observons que si $X$ est une vari\'et\'e projective, d'apr\`es le lemme 2.8,
$\lambda_{s+s'}(f)\leq c\lambda_s(f)\lambda_{s'}(f)$.
Dans ce cas, la
suite $[\lambda_s(f)]^{1/n}$ converge vers sa borne inf\'erieure $\inf_{n\geq 1}
[\lambda_s(f^n)]^{1/n}$ (\voir \cite{DinhSibony3}).
\par
\
\\
{\bf 3.6}. {\it Exemples.}
\\
\ \par
(a)
Notons $\P^k$ l'espace projectif complexe et $\GG(k-l+1,k+1)$ la
grassmannienne qui param\`etre la famille des sous-espaces
projectifs de dimension $k-l$ de $\P^{k}$. Pour
$\scheck\in \GG(k-l+1,k+1)$, soit $\P^{k-l}_\shat$
le sous-espace projectif de
dimension $k-l$ correspondant. Posons
$$\Gamma_1:=\{(z,\shat)\in \P^k\times \GG(k-l+1,k+1),\ z\in
\P^{k-l}_\shat\}.$$
La transformation m\'eromorphe $\Psi_1$ de $\P^k$ dans $\GG(k-l+1,k+1)$
associ\'ee \`a la vari\'et\'e $\Gamma_1$ est de
codimension $k-l$. En effet, si $\shat$ est un point de $\GG(k-l+1,k+1)$,
$\Psi_1^{-1}(\shat)$ est le sous-espace projectif
$\P^{k-l}_{\shat}$ de $\P^k$. 
\\
\ \par
Donnons une autre mani\`ere de voir ces transformations
m\'eromorphes. Soit $\P^{k*}:=\GG(k,k+1)$ le dual de $\P^k$
et soit $\GG^*(l,k+1)$ la
grassmannienne qui param\`etre les
sous-espaces projectifs de dimension $l-1$ de $\P^{k*}$. Elle est
biholomorphe \`a $\GG(k-l+1,k+1)$. Pour
tout $\scheck\in \GG^*(l,k+1)$ notons $\P^{(l-1)*}_\scheck$ le sous-espace
projectif de $\P^{k*}$
associ\'e \`a $\scheck$. On choisit $l$ points $s_1,\ldots, s_l$ de
$\P^{(l-1)*}_\scheck$ qui engendrent $\P^{(l-1)*}_\scheck$. Notons
$\P^{k-1}_{s_i}$ l'hyperplan de $\P^k$ associ\'e \`a $s_i$ et
$\P^{k-l}_\scheck:=\P^{k-1}_{s_1}\cap \ldots \cap \P^{k-1}_{s_l}$. Le
sous-espace projectif $\P^{k-l}_\scheck$ de $\P^k$ est ind\'ependant
du choix des $s_i$.
\par
Posons
$$\Gamma_2:=\{(z,\scheck)\in \P^k\times \GG^*(l,k+1),
\ z\in \P^{k-l}_\scheck\}.$$
La transformation m\'eromorphe $\Psi_2$ de $\P^k$ dans $\GG^*(l,k+1)$
associ\'ee \`a $\Gamma_2$ est de codimension $k-l$.
Si $\scheck$ est un point de $\GG^*(l,k+1)$,
$\Psi_2^{-1}(\scheck)$ est le sous-espace
projectif $\P^{k-l}_{\scheck}$ de
$\P^k$.
\\
\ \par
(b) 
Consid\'erons {\it l'espace multiprojectif}
$\P^{k,l*}:=\P^{k*}\times\cdots\times \P^{k*}$ ($l$ fois). 
Posons
\begin{eqnarray*}
\Gamma_3 & := & \Big\{(\ss,\scheck)\in \P^{k,l*}\times\GG^*(l,k+1),\\
& & \ \ \ 
\ss=(s_1,\ldots,s_l), \P^{k-l}_\scheck\subset \P^{k-1}_{s_i} \mbox{
pour } i=1,\ldots,l\Big\}.
\end{eqnarray*}
Notons $\Pi_l$ la transformation m\'eromorphe de $\P^{k,l*}$ dans
$\GG^*(l,k+1)$ associ\'ee \`a $\Gamma_3$.
C'est une application m\'eromorphe surjective.
Soit $\overline\Pi_l$ son adjoint. La compos\'ee
$\Psi_3:=\overline \Pi_l \circ\Psi_2$ est une transformation
m\'eromorphe de $\P^k$ dans $\P^{k,l*}$. Pour un point
$\ss=(s_1,\ldots,s_l)$
g\'en\'erique de $\P^{k,l*}$, la fibre $\Psi_3^{-1}(\ss)$ est le
sous-espace projectif $\P^{k-l}_\ss:=\P^{k-1}_{s_1}\cap\ldots
\cap \P^{k-1}_{s_l}$ de $\P^k$.
\\
\
\par
(c) Nous allons \'etendre la d\'efinition des
transformations m\'eromorphes de (a) et (b) avec, pour espace
d'arriv\'ee, un espace projectif de sections holomorphes.
Consid\'erons une vari\'et\'e 
projective $X$ et soit $L$ un fibr\'e en droites ample sur $X$. Pour
$n\geq 1$, notons
$\HH^0(X,L^n)$ l'espace des sections holomorphes de
$L^n:=L\otimes \cdots \otimes L$ ($n$ fois),
$\P\HH^0(X,L^n)$ l'espace projectif associ\'e et $k_n$ la
dimension de $\P\HH^0(X,L^n)$.
Pour tout $s^*\in\P\HH^0(X,L^n)^*$ notons
$H_{s^*}$ l'hyperplan projectif de $\P\HH^0(X,L^n)$
associ\'e \`a $s^*$. Rappelons que 
$\P\HH^0(X,L^n)$ est aussi le dual de $\P\HH^0(X,L^n)^*$.
Pour tout $s\in \P\HH^0(X,L^n)$, notons $H^*_s$ l'hyperplan
projectif de $\P\HH^0(X,L^n)^*$ associ\'e \`a $s$.
\par
Pour $x\in X$, notons
$s^*_x\in\P\HH^0(X,L^n)^*$
le point tel que l'hyperplan $H_{s^*_x}$ soit 
l'ensemble des sections s'annulant en $x$.
Consid\'erons
l'application holomorphe $\Phi_n$ de $X$ dans 
$\P\HH^0(X,L^n)^*$ d\'efinie
par $x\mapsto \Phi_n(x):=s^*_x$.
Puisque $L$ est ample, pour $n$ assez grand, 
l'application $\Phi_n$
d\'efinit un plongement de $X$ dans $\P\HH^0(X,L^n)^*$, c'est le
plongement de Kodaira. Observons que
$\Phi_n^{-1}(H^*_s\cap\Phi_n(X))$ est l'ensemble des z\'eros de $s$.
D'apr\`es le th\'eor\`eme de Bertini \cite[p.141]{Shafarevich},
cette intersection est
transverse et d\'efinit une hypersurface
lisse de $X$ pour tout $s$ hors d'un sous-ensemble analytique de
$\P\HH^0(X,L^n)$.
\par
Notons $\GG_{l,n}^X$ la grassmannienne des
sous-espaces projectifs de dimension $l-1$ de $\P\HH^0(X,L^n)$.
On d\'efinit comme dans (a) une transformation m\'eromorphe
$\Psi_{l,n}$ de $\P\HH^0(X,L^n)^*$ dans $\GG_{l,n}^X$. Pour tout
point $\scheck\in \GG_{l,n}^X$, $\Psi_{l,n}^{-1}(\scheck)$ est un sous-espace
projectif de dimension $k_n-l$ de $\P\HH^0(X,L^n)^*$. Posons
$R_{l,n}:=\Psi_{l,n}\circ\Phi_n$. C'est une transformation
m\'eromorphe de codimension $k-l$ de $X$ dans $\GG_{l,n}^X$.
Pr\'ecisons cela.
\par
Notons 
$s_1,\ldots , s_l$ des points qui
engendrent le sous-espace $\P^{l-1}_\scheck$ de dimension $l-1$
de $\P\HH^0(X,L^n)$, 
pour $\scheck\in \GG_{l,n}^X$. Alors
$\Psi_{l,n}^{-1}(\scheck)$ est \'egal au sous-espace projectif
$H_{s_1}^*\cap\ldots\cap H_{s_l}^*$ de dimension $k_n-l$ de
$\P\HH^0(X,L^n)^*$. On en d\'eduit que $R_{l,n}^{-1}(\scheck)$ est
l'ensemble $\ZZ_\scheck$ des z\'eros communs
des sections $s_1$, $\ldots$, $s_l$. Cet ensemble
ne d\'epend pas du choix des $s_i$.
D'apr\`es le th\'eor\`eme de Bertini \cite[p.141]{Shafarevich},
pour $\scheck\in \GG_{l,n}^X$
hors d'un sous-ensemble analytique, l'intersection de
$\Psi_{l,n}^{-1}(\scheck)$ avec $\Phi_n(X)$ est transverse.
Pour un tel $\scheck$, $\ZZ_\scheck$ est lisse.
\\
\ \par
Soit $\Pi_{l,n}$ l'application m\'eromorphe de
$\P^X_{l,n}:=\P\HH^0(X,L^n)\times \cdots\times\P\HH^0(X,L^n)$ ($l$
fois) dans $\GG^X_{l,n}$ d\'efinie comme dans (b).
Soit $\overline \Pi_{l,n}$ sont
adjoint. 
Posons $F_{l,n}:=\overline\Pi_{l,n}\circ R_{l,n}$. C'est une
transformation m\'eromorphe de $X$ dans $\P_{l,n}^X$. Pour tout
$\ss=(s_1,\ldots,s_l)\in \P_{l,n}^X$, la fibre $F_{l,n}^{-1}(\ss)$
est l'ensemble des z\'eros communs des sections holomorphes  
$s_1$, $\ldots$, $s_l$ de $L^n$. Pour un $\ss$
g\'en\'erique, cette fibre est \'egale \`a $\ZZ_\scheck$ avec
$\scheck:=\Pi_{l,n}(\ss)$. En particulier, c'est
un sous-ensemble analytique lisse de dimension $k-l$,
sans multiplicit\'e.
\\
\ \par
Soient $z_1,\ldots,z_m$ des points de $X$.
Pour $n$ suffisamment grand, on peut trouver un sous-espace
projectif de dimension $l$ de $\P\HH^0(X,L^n)^*$ passant par
$\Phi_n(z_1),\ldots, \Phi_n(z_m)$. C'est-\`a-dire qu'on peut
trouver une famille libre de $l$
sections holomorphes de $L^n$ qui s'annulent
simultan\'ement en $z_1,\ldots,z_m$. On peut d\'efinir
les auto-intersections des transformations m\'eromorphes 
$F_{l,n}$ et $R_{l,n}$, en posant
$F_{l,n,m}:=F_{l,n}\cap\ldots\cap F_{l,n}$ ($m$ fois) et
$R_{l,n,m}:=R_{l,n}\cap\ldots \cap R_{l,n}$ ($m$ fois). Ce sont
des transformations m\'eromorphes de $X^m$ dans
$\P_{l,n}^X$ et dans $\GG_{l,n}^X$. Un point $(z_1,\ldots, z_m)\in X^m$
appartient \`a la fibre
$F_{l,n,m}^{-1}(\ss)$ (resp. $R_{l,m,n}^{-1}(\scheck)$) si et
seulement si $\ZZ_\scheck$ passe par
$z_1$, $\ldots$, $z_m$.   
\\
\ \par
(d) Soit $f$ une application m\'eromorphe de $\P^k$ dans
$\P^k$. Pour \'etudier les images r\'eciproques par $f^n$ des
sous-espaces de dimension $k-l$, nous introduisons les
transformations m\'eromorphes $F_n:=\Psi_2\circ f^n$ avec $\Psi_2$
d\'efinie dans l'exemple (a)
On se ram\`ene \`a l'\'etude des images
r\'eciproques des points de $\GG^*(l,k+1)$ par $F_n$.
\sectionbis{Distribution des pr\'eimages de sous-vari\'et\'es}
Dans ce paragraphe, nous donnons des solutions au probl\`eme\break
d'\'equidistribution dans un cadre abstrait. Nous allons, dans les
paragraphes 6 et 7, appliquer ces r\'esultats aux cas
particuliers que nous avons discut\'es dans l'introduction.
Des id\'ees analogues permettent de construire, au
paragraphe 5, les mesures 
d'\'equilibre pour les correspondances m\'eromorphes. 
\\
\ \par
Soit $\sigma_n$ une mesure de probabilit\'e PLB sur $X_n$.
On munit $\XX:=\prod_{n\geq 1} X_n$ de
la mesure de probabilit\'e $\sigma$, \'egale au
produit des $\sigma_n$. 
Consid\'erons
des transformations m\'eromorphes $F_n$ de m\^eme codimension $l$
de $X$ dans $X_n$, $0\leq l<k$.
Soit $\xx=(x_1,x_2,\ldots)\in \XX$.
Si $x_n$ n'appartient pas au deuxi\`eme ensemble d'ind\'etermination
$I_2(F_n)$ de $F_n$,
le courant
$T^\xx_n:=(F_n)^*(\delta_{x_n})$ 
est bien d\'efini. C'est un 
courant d'int\'egration sur une cha\^{\i}ne holomorphe de
dimension $l$ de $X$. Posons $T_n:=(F_n)^*(\sigma_n)$.
\par
Soient $\delta_n$ et $d_n$ les
degr\'es interm\'ediaires d'ordre $k_n-1$ et d'ordre $k_n$ de $F_n$.
Posons 
$\RR_{1,n}:=\RR_1(X_n,\omega_n,\sigma_n)$,
$\RR_{2,n}:=\RR_2(X_n,\omega_n,\sigma_n)$ et
$\Delta_n(t):=\Delta(X_n,\omega_n,\sigma_n,t)$ pour tout $t\in\R$
(\voir paragraphe 2).
On a le th\'eor\`eme suivant.
\begin{theoreme} 
Supposons que la suite $(\RR_{1,n}\delta_n d_n^{-1})_{n\geq 1}$
tend vers $0$
et que l'une des deux propri\'et\'es suivantes soit satisfaite:
\begin{enumerate}
\item[{\rm (1)}] La s\'erie $\sum_{n\geq 1}
\RR_{2,n}\delta_n d_n^{-1}$
converge.
\item[{\rm (2)}]
La s\'erie $\sum_{n\geq 1}
\Delta_n(\delta_n^{-1} d_nt)$
converge pour tout $t>0$.
\end{enumerate}
Alors pour $\sigma$-presque tout $\xx\in \XX$,
la suite $\langle d_n^{-1}(T^\xx_n-T_n),\psi\rangle$ tend
vers $0$ uniform\'ement sur les ensembles
born\'es, en norme ${\cal C}^2$, de
$(l,l)$-formes test $\psi$ sur $X$.
\end{theoreme}
\par
Nous allons montrer les estimations utiles en nous limitant \`a une
transformation m\'eromorphe $F$ de codimension $l$,
$0\leq l\leq k-1$, de
$(X,\omega)$ dans une vari\'et\'e k\"ahl\'erienne compacte
$(X',\omega')$ de dimension $k'$.
Soit $\sigma'$ une mesure de probabilit\'e PLB sur $X'$.
Posons $\widetilde T:=F^*({\omega'}^{k'})$,
$T:=F^*(\sigma')$ et
$T^{x'}:=F^*(\delta_{x'})$
pour tout $x'\in X'\setminus I_2(F)$.
\par
Soient $\delta$ et $d$ les
degr\'es interm\'ediaires d'ordre $k'-1$ et d'ordre $k'$ de $F$.
Posons $\RR_i:=\RR_i(X',\omega',\sigma')$ pour $i=1,2,3$,
$\Delta(t):=\Delta(X',\omega',\sigma',t)$ pour $t\in\R$. Pour
tout $\epsilon>0$ d\'efinissons
$$E(\epsilon):=\bigcup_{\|\psi\|_{{\cal C}^2(X)}\leq 1}
\left\{x'\in X',\ \left|\langle d^{-1}(T^{x'}-T),\psi\rangle
\right| \geq \epsilon \right\}.$$
\par
Posons $S:=F_*(\omega^{l+1})$. C'est un courant de bidegr\'e
$(1,1)$ sur $X'$. 
Par d\'efinition de $\delta$,
on a $\|S\|\leq \delta$. D'apr\`es la proposition 2.2, il existe
une fonction q.p.s.h. $\varphi$ v\'erifiant
\begin{eqnarray}
\int_{X'}\varphi\d\sigma'=0 & \mbox{ et } & 
\ddc\varphi - S \geq -r(X',\omega')\delta \omega'
\end{eqnarray}
Par d\'efinition de $\RR_i$,
on a
$$\varphi\leq \delta\RR_1,\ \ \ \|\varphi\|_{\Lone(\sigma')}
\leq \delta\RR_2\ \ \ \mbox{et}\ \ \
\left|\int \varphi{\omega'}^{k'}\right|\leq \delta\RR_3.$$
\begin{lemme} Soit $\psi$ une $(l,l)$-forme test de classe
${\cal C}^2$ sur $X$. Alors on a 
\begin{enumerate}
\item[{\rm (a)}]
$\displaystyle
\int_{X'}\left|\langle T^{x'}-T,\psi\rangle \right| \d\sigma'(x')
\leq 2 \|\psi\|_{{\cal C}^2} \delta \RR_2$.
\item[{\rm (b)}]
$\displaystyle
\left|\langle T^{x'}-T,\psi\rangle\right| \leq \|\psi\|_{{\cal C}^2}
\big(3\delta \RR_1-\varphi(x')\big)$.
\item[{\rm (c)}]
$\displaystyle
\left|\langle T-\widetilde T,
\psi\rangle\right| \leq 2\|\psi\|_{{\cal C}^2}
\delta \RR_3$.
\item[{\rm (d)}] $\displaystyle
\sigma'(E(\epsilon))\leq
\Delta(\epsilon\delta^{-1}d -3\RR_1)$.
\end{enumerate}
\end{lemme}
\begin{preuve}
Nous devons estimer $\langle T^{x'}-T,\psi\rangle$
et $\langle T-\widetilde T,\psi\rangle$.
Ecrivons $\ddc\psi=\Omega^+-\Omega^-$ avec
$\Omega^\pm$ des $(l+1,l+1)$-formes positives ferm\'ees
telles que $\Omega^\pm\leq \|\psi\|_{{\cal C}^2}\omega^{l+1}$.
Posons $\phi:=F_*(\psi)$ et
$S^\pm:=F_*(\Omega^\pm)$. On a $\ddc\phi=S^+-S^-$ et
$\langle T^{x'},\psi\rangle = \phi(x')$ pour $x'\not\in I_2(F)$.
On a aussi $S^\pm\leq \|\psi\|_{{\cal C}^2}S$
et donc $\|S^\pm\|\leq \|\psi\|_{{\cal C}^2}\delta$.
\par
D'apr\`es la proposition
2.2, on peut choisir
les fonctions q.p.s.h. $\varphi^\pm$
telles que $\int_{X'} \varphi^\pm\d\sigma'=0$
et telles que 
\begin{eqnarray}
-r(X',\omega')\|\psi\|_{{\cal C}^2}\delta \omega' & \leq &
\ddc\varphi^+ - S^+ = \nonumber \\
& = & \ddc\varphi^- - S^-\leq
r(X',\omega')\|\psi\|_{{\cal C}^2}\delta \omega'. 
\end{eqnarray}
Par d\'efinition des $\RR_i$
on a 
\begin{eqnarray}
\varphi^\pm \leq \|\psi\|_{{\cal C}^2}\delta\RR_1,\ 
\|\varphi^\pm\|_{\Lone(\sigma')}\leq \|\psi\|_{{\cal
C}^2}\delta\RR_2,\ 
\left|\int\varphi^\pm{\omega'}^{k'}\right|\leq \|\psi\|_{{\cal
C}^2}\delta\RR_3.
\end{eqnarray}
La fonction $\phi-(\varphi^+-\varphi^-)$ est
constante car elle est pluriharmonique.
Par cons\'equent, on a pour $x'\not\in I_2(F)$
\begin{eqnarray}
\langle T^{x'}-T,\psi\rangle & = &
\langle \delta_{x'}-\sigma',\phi\rangle =
\langle \delta_{x'}-\sigma',\varphi^+-\varphi^-\rangle
\end{eqnarray}
Puisque $\int\varphi^\pm\d\sigma'=0$, on  d\'eduit de
la relation (4.4) que 
$$\langle T^{x'}-T,\psi\rangle = \varphi^+(x')-\varphi^-(x').$$
\par
(a) On a
\begin{eqnarray*}
\int_{X'}\left|\langle T^{x'}-T,\psi\rangle \right| \d\sigma'(x')
& = & \int \left|\varphi^+(x')-\varphi^-(x')\right| \d\sigma'(x') \\
&\leq & \|\varphi^+\|_{\Lone(\sigma')} +
\|\varphi^-\|_{\Lone(\sigma')}\\ 
& \leq & 2 \|\psi\|_{{\cal C}^2} \delta \RR_2.
\end{eqnarray*}
\par
(b)
Posons $h:=\|\psi\|_{{\cal C}^2}\varphi-\varphi^+$. On a par
d\'efinition de $\varphi$ et $\varphi^+$:
$$\ddc h -(\|\psi\|_{{\cal C}^2}S-S^+) \geq
-2r(X',\omega') \|\psi\|_{{\cal C}^2}\delta\omega'.$$
D'une part, la fonction $h$ v\'erifie
$\int_{X'}h\d\sigma'=0$. D'autre part,
$\|\psi\|_{{\cal C}^2}S-S^+$
est un courant positif ferm\'e dont la masse est major\'ee par
$\|\psi\|_{{\cal C}^2}\delta$.
D'o\`u on d\'eduit, par d\'efinition de $\RR_1$,
que $\sup_{X'} h\leq 2\|\psi\|_{{\cal C}^2}
\delta \RR_1$. On d\'eduit de cette in\'egalit\'e et de (4.3) les in\'egalit\'es
suivantes 
$$\|\psi\|_{{\cal C}^2}\varphi(x')-2\|\psi\|_{{\cal C}^2}
\delta\RR_1 \leq
\varphi^+(x')\leq \|\psi\|_{{\cal C}^2}\delta
\RR_1.$$
Avec une estimation analogue pour $\varphi^-$, on obtient finalement
$$\left|\langle T^{x'}-T,\psi\rangle\right| =
\left|\varphi^+(x')-\varphi^-(x')\right| \leq \|\psi\|_{{\cal
C}^2} (3\delta \RR_1 -\varphi(x')).$$
\par
(c) D'apr\`es (4.3), on a 
\begin{eqnarray*}
\left|\langle T-\widetilde T,\psi\rangle\right| & = &
\left|\langle
\sigma-{\omega'}^{k'},\varphi^+-\varphi^-\rangle\right| 
 =  \left|\langle
{\omega'}^{k'},\varphi^+-\varphi^-\rangle\right| \\
& \leq & \left|\langle {\omega'}^{k'},\varphi^+\rangle\right| +
\left|\langle {\omega'}^{k'},\varphi^-\rangle\right|
 \leq  2\delta\RR_3.
\end{eqnarray*}
\par
(d) 
D'apr\`es (b) appliqu\'e au cas o\`u $\|\psi\|_{{\cal C}^2}\leq1$,
l'ensemble $E(\epsilon)$ est contenu dans
$$E'(\epsilon)
:=\{x'\in X',\ \varphi(x')\leq -\epsilon d +3\delta \RR_1\}.$$
Par d\'efinition de $\Delta(t)$, les relations (4.1)
entra\^{\i}nent que
$$\sigma'(E(\epsilon))\leq
\sigma'(E'(\epsilon))\leq \Delta(\epsilon \delta^{-1} d -3\RR_1).$$
\end{preuve}
{\bf Fin de la d\'emonstration du th\'eor\`eme 4.1.}
Posons $S_n:=(F_n)_*(\omega^{l+1})$.
On a $\|S_n\|\leq \delta_n$. D'apr\`es la proposition 2.2, il existe
une fonction q.p.s.h. $\varphi_n$ v\'erifiant
$\int_{X_n}\varphi_n\d\sigma_n=0$ telle que
$$\ddc\varphi_n - S_n \geq -r(X_n,\omega_n)\delta_n \omega_n.$$
Par d\'efinition de $\RR_{1,n}$ et $\RR_{2,n}$,
on a $\varphi_n\leq \delta_n\RR_{1,n}$ et
$\|\varphi_n\|_{\Lone(\sigma_n)}
\leq \delta_n\RR_{2,n}$. 
\\
\ \par
(1) Consid\'erons la
fonction r\'eelle positive $\Phi$ sur $\XX$ 
$$\Phi(\xx):=\sum_{n\geq 1}d_n^{-1}|\varphi_n(x_n)|.$$
On a
$$\int_\XX \Phi\d\sigma = \sum_{n\geq 1} d_n^{-1}
\|\varphi_n\|_{\Lone(\sigma_n)}
\leq \sum_{n\geq 1}\RR_{2,n}\delta_nd_n^{-1}.$$
Par hypoth\`ese, la derni\`ere s\'erie converge, donc
$\Phi(\xx)$ est finie $\sigma$-presque partout, et 
$d_n^{-1}\varphi_n(x_n)$ tend vers $0$ pour $\sigma$-presque tout
$\xx\in \XX$.
\par
Fixons $\xx=(x_n)\in \XX$ tel que $x_n\not\in I_2(F_n)$ et tel que
$d_n^{-1}\varphi_n(x_n)$ tende vers $0$.
Soit $\psi$ une $(l,l)$-forme de classe ${\cal C}^2$ sur $X$.
D'apr\`es le lemme 4.2, on a
$$\left|\langle T^\xx_n-T_n, \psi\rangle \right|  \leq  
\|\psi\|_{{\cal C}^2}(3\delta_n\RR_{1,n}-\varphi_n(x_n)).$$
Comme $\RR_{1,n}\delta_nd_n^{-1}$ et $d_n^{-1}\varphi_n(x_n)$
tendent vers $0$, la suite 
$d_n^{-1}\langle T^\xx_n-T_n, \psi\rangle$ tend aussi vers $0$
uniform\'ement sur les ensembles born\'es en norme ${\cal C}^2$ de
$(l,l)$-formes test $\psi$ sur $X$.
\\
\ \par
(2) Posons pour tout $\epsilon>0$
$$E_n(\epsilon):=\bigcup_{\|\psi\|_{{\cal C}^2(X)}\leq 1}
\left\{x_n\in X_n,\ \left|\big\langle d_n^{-1}(T^{x_n}-T),\psi
\big\rangle
\right| \geq \epsilon \right\}.$$
Par hypoth\`ese, on a
$\RR_{1,n}=\o(\delta_n^{-1}d_n)$. D'apr\`es le lemme 4.2(d), pour
$n$ assez grand, on a
\begin{eqnarray}
\sigma_n(E_n(\epsilon))
& \leq &  \Delta_n(\epsilon \delta_n^{-1} d_n-3\RR_{1,n})
\leq \Delta_n(\epsilon \delta_n^{-1}
d_n/2)
\end{eqnarray}
Par hypoth\`ese, la s\'erie
$\sum\Delta(\epsilon \delta_n^{-1}d_n/2)$ 
converge. On en d\'eduit que la s\'erie $\sum
\sigma_n(E_n(\epsilon))$ converge pour tout $\epsilon>0$. Ceci
implique la convergence annonc\'ee.
\par
\hfill $\square$
\\
\par
La proposition suivante permet de comparer les courants obtenus
en prenant les images r\'eciproques
de $\sigma_n$ et de la forme volume $\omega_n^{k_n}$.
\begin{proposition}  Supposons que la suite
$\RR_3(X_n,\omega_n,1)\delta_n d_n^{-1}$ tend
vers $0$. Alors 
$\langle d_n^{-1}(T_n -F_n^*(\omega_n^{k_n})),\psi\rangle$
tend vers $0$ uniform\'ement sur les ensembles
born\'es, en norme ${\cal
C}^2$, de $(l,l)$-formes $\psi$ sur $X$.
\end{proposition}
\begin{preuve} Il suffit d'appliquer le lemme 4.2(c) avec
les estimations comme dans la d\'emonstration du th\'eor\`eme 4.1.
\end{preuve}
\par
Posons
$\RR^*_n:=\RR^*_2(X_n,\omega_n,2)$.
Soient
$\nu_n=h_n\omega_n^{k_n}$ et $\nu_n'=h_n'\omega_n^{k_n}$ des mesures
de probabilit\'e sur $X_n$ o\`u $h_n$ et $h_n'$ sont des
fonctions dans ${\rm L}^2(X_n)$. On a le r\'esultat suivant.
\begin{theoreme} 
Supposons que 
$\|h_n-h_n'\|_{\LL^2(X_n)}=\o(\delta_n^{-1}d_n(\RR^*_n)^{-1})$.
Alors la suite 
$$d_n^{-1}\big\langle (F_n)^*(\nu_n)-(F_n)^*(\nu_n'),\psi
\big\rangle$$
tend vers $0$ quand $n$ tend
vers l'infini, uniform\'ement
sur les ensembles born\'es, en norme ${\cal
C}^2$, de $(l,l)$-formes test $\psi$ sur $X$.
\end{theoreme}
\begin{preuve} Utilisons les notations du th\'eor\`eme 4.1.
Posons $\phi_n:=(F_n)_*(\psi)$.  Il
existe des constantes $a_n$ et des fonctions q.p.s.h.
$\varphi_n^\pm$ telles que $\phi_n=\varphi^+_n-\varphi^-_n+a_n$ et
$\|\varphi_n^\pm\|_{\LL^2(X_n)} \leq \|\psi\|_{{\cal C}^2}
\RR^*_n\delta_n$. 
Puisque $\nu_n$ et $\nu_n'$ ont la
m\^eme masse, l'in\'egalit\'e de Cauchy-Schwarz entra\^{\i}ne que
\begin{eqnarray*}
d_n^{-1}\big|\big\langle (F_n)^*(\nu_n)-(F_n)^*(\nu_n'), \psi
\big\rangle\big|
& = &  d_n^{-1}\big|\big\langle (h_n-h_n')\omega_n^{k_n}, \phi_n-a_n
\big\rangle\big| \\
& \leq & 
d_n^{-1}\|h_n-h_n'\|_{\LL^2(X_n)} \|\phi_n-a_n\|_{\LL^2(X_n)}\\
& \leq & 2\|\psi\|_{{\cal C}^2}
\delta_n d_n^{-1}\RR^*_n \|h_n-h_n'\|_{\LL^2(X_n)}.
\end{eqnarray*}
Par hypoth\`ese, la derni\`ere expression tend vers $0$. Ceci
implique le th\'eor\`eme.
\end{preuve}
\begin{remarque} \rm
Dans le cas o\`u  les $(X_n,\omega_n,\sigma_n)$ appartiennent \`a
une famille compacte lisse,
les constantes $\RR_{1,n}$ et $\RR^*_n$ sont uniform\'ement
born\'ees (\voir remarque 2.3).
Les hypoth\`eses des th\'eor\`emes 4.1 et
4.4 ne font alors intervenir que les degr\'es interm\'ediaires de
$F_n$. Ces degr\'es sont calcul\'es cohomologiquement.
\end{remarque}
\par
Dans la suite, nous
consid\'erons des transformations m\'eromorphes $F_n$ de
$X$ dans une m\^eme vari\'et\'e $X'$, \cad qu'on suppose que
$(X_n,\omega_n)=(X',\omega')$ pour tout $n\geq 1$.
Nous avons alors le th\'eor\`eme suivant.
\begin{theoreme} Supposons
que la s\'erie $\sum_{n\geq 1} d_n^{-1}\delta_n$
converge. Alors il existe un sous-ensemble pluripolaire $\E$ de
$X'$ tel que 
pour tout $x'\in X'\setminus \E$ la suite 
$$\big\langle d_n^{-1}((F_n)^*(\delta_{x'})-
(F_n)^*({\omega'}^{k'})),
\psi \big\rangle$$
tend vers $0$ uniform\'ement sur les ensembles
born\'es, en
norme ${\cal C}^2$, de $(l,l)$-formes test $\psi$ sur $X$.
\end{theoreme}
\begin{preuve}
Posons $r:=r(X',\omega')$, $\RR^*_1:=\RR^*_1(X',\omega')$,
$\RR^*_2:=\RR^*_2(X',\omega',1)$ et
$S_n:=(F_n)_*(\omega^{l+1})$.
On a $\|S_n\|\leq \delta_n$. D'apr\`es la proposition 2.2, il existe
une fonction q.p.s.h. $\varphi_n$ v\'erifiant
$$\int_{X'}\varphi_n{\omega'}^{k'}=0 \ \ \mbox{ et } \ \ 
\ddc\varphi_n - S_n \geq -r\delta_n \omega'.$$
Par d\'efinition de $\RR^*_1$ et $\RR^*_2$,
on a $\varphi_n\leq \RR^*_1\delta_n$ et
$\|\varphi_n\|_{\Lone(X')}
\leq \RR^*_2\delta_n$. 
On en d\'eduit que 
la s\'erie
$$\Phi(x'):=\sum_{n\geq 1} d_n^{-1}\varphi_n(x')$$
converge vers une fonction q.p.s.h. Posons
$$\E:=\bigcup_{n\geq 1}I_2(F_n)\cup
(\Phi=-\infty).$$
D'apr\`es la proposition A.1, c'est un ensemble pluripolaire de $X'$.
Pour $x'\in X'\setminus \E$,
on a $\lim d_n^{-1}\varphi_n(x')=0$.
D'apr\`es le lemme 4.2(b) appliqu\'e \`a la mesure
$\sigma':={\omega'}^{k'}$, la suite
$\langle d_n^{-1}((F_n)^*(\delta_{x'})-(F_n)^*({\omega'}^{k'})),
\psi\rangle$ tend vers $0$ uniform\'ement sur les ensembles
born\'es, en
norme ${\cal C}^2$, de $(l,l)$-formes test $\psi$ sur $X$.
\end{preuve}
\begin{remarques} \rm
Supposons que $\sum_{n\geq 1}\delta_n^pd_n^{-p}<+\infty$ pour un
$p>1$. Soit $\nu$ une mesure de probabilit\'e telle que $\int_X
|\varphi|^p\d\nu<+\infty$ pour toute fonction q.p.s.h. $\varphi$ sur
$X$. En consid\'erant la s\'erie
$\Phi_p(x'):=\sum_{n\geq 1}d_n^{-p}|\varphi_n(x')|^p$,
on montre comme au th\'eor\`eme 4.6 que pour
$\nu$-presque tout $x'\in X'$ la suite de courants
$d_n^{-1}(F_n^*(\delta_{x'})
-F_n^*({\omega'}^{k'}))$ tend faiblement vers $0$. On a \'egalement
la convergence uniforme sur les formes born\'ees en norme ${\cal
C}^2$.
\par
Soit $(c_n)$ une suite de nombres r\'eels positifs telle
que la s\'erie $\sum_{n\geq 1} c_n\delta_n$
converge. Alors il existe un sous-ensemble pluripolaire $\E$ de
$X'$ tel que 
pour $x'\in X'\setminus \E$ la suite 
$c_n\big\langle (F_n)^*(\delta_{x'})-
(F_n)^*({\omega'}^{k'}),\psi \big\rangle$
tend vers $0$ uniform\'ement sur les ensembles
born\'es, en
norme ${\cal C}^2$, de $(l,l)$-formes test sur $X$. Dans ce cas, l'hypoth\`ese
$\sum_{n\geq 1} \delta_n d_n^{-1}<+\infty$ n'est pas n\'ecessaire.
La m\^eme remarque est valable pour le th\'eor\`eme 4.1.
\end{remarques}
\begin{theoreme} Soit $\sigma'$ une mesure de probabilit\'e PLB
sur $X'$. Supposons que la s\'erie
$\sum_{n\geq 1} \Delta(X',\omega',\sigma',t\delta_n^{-1}d_n)$
converge pour tout $t>0$.
Alors pour $\sigma'$-presque tout $x'\in X'$,
la suite 
$\langle d_n^{-1}((F_n)^*(\delta_{x'})-(F_n)^*({\omega'}^{k'})),
\psi \rangle$ tend vers $0$ uniform\'ement sur les ensembles
born\'es, en
norme ${\cal C}^2$, de $(l,l)$-formes test $\psi$ sur $X$.
\end{theoreme}
\begin{preuve} Soient $r:=r(X',\omega')$,
$\RR^*_1:=\RR^*_1(X',\omega')$ et
$\Delta(t):=\Delta(X',\omega',\sigma',t)$.
Posons pour tout $\epsilon>0$
$$E_n(\epsilon):=\bigcup_{\|\psi\|_{{\cal C}^2(X)}\leq 1}
\left\{x'\in X',\ \left|\big\langle d_n^{-1}(F_n^*(\delta_{x'})
-F_n^*({\omega'}^{k'}))
,\psi\big\rangle\right| \geq \epsilon \right\}.$$
Le lemme 4.2(d)  entra\^{\i}ne que
$$\sigma'(E_n(\epsilon))\leq \Delta(\epsilon \delta_n^{-1}d_n
-3\RR^*_1).$$
L'hypoth\`ese entra\^{\i}ne que $\lim\delta_n^{-1}d_n =+\infty$ et la
s\'erie $\sum_{n\geq 1} \Delta(\epsilon \delta_n^{-1}d_n
-3\RR^*_1)$ converge. On en d\'eduit que la s\'erie 
$\sum_{n\geq 1} \sigma'(E_n(\epsilon))$ converge pour tout
$\epsilon>0$.
Le th\'eor\`eme en d\'ecoule.
\end{preuve}
\begin{corollaire} Soit $\sigma'$ une mesure de probabilit\'e
$(c,\alpha)$-mod\'er\'ee
sur $X'$ avec $c>0$ et $\alpha>0$. Supposons que
$\sum_{n\geq 1} \exp(-\delta_n^{-1}d_nt)$
converge pour tout $t>0$ (par exemple si
$\delta_nd_n^{-1}=\o(1/\log n)$).
Alors pour $\sigma'$-presque tout $x'\in X'$,
la suite 
$\langle d_n^{-1}((F_n)^*(\delta_{x'})-(F_n)^*({\omega'}^{k'})),
\psi \rangle$ tend vers $0$ uniform\'ement sur les ensembles
born\'es, en
norme ${\cal C}^2$, de $(l,l)$-formes test $\psi$ sur $X$.
\end{corollaire}
\begin{preuve} Soit $r:=r(X',\omega')$.
La mesure $\sigma'$ \'etant $(c,\alpha)$-mod\'er\'ee,
la proposition 2.5 entra\^{\i}ne que
$$\Delta(t \delta_n^{-1}d_n)\leq 
c\exp(-\alpha r^{-1} t\delta_n^{-1}d_n).$$
Par cons\'equent, la s\'erie $\sum\Delta(t\delta_n^{-1}d_n)$ converge pour tout $t>0$.
On peut appliquer le th\'eor\`eme 4.8.
\end{preuve}
\par
Soit ${\cal H}=(\H_y)_{y\in Y}$
une famille m\'eromorphe adapt\'ee de sous-ensembles
analytiques de dimension $m$ d'une vari\'et\'e projective 
$X'$ associ\'ee \`a une transformation
m\'eromorphe $G:X'\longrightarrow Y$. Supposons que
$l+m<k$ et que pour tout $n$, $\H$ soit $F_n$-r\'eguli\`ere. 
Pour $y\in Y$ g\'en\'erique, d'apr\`es le lemme 2.8, les courants
$[\lambda_{k'-m}(F_n)]^{-1} (F_n)^*[\H_y]$ sont bien d\'efinis et de masse born\'ee
ind\'ependemment de $n$. On a le corollaire suivant.
\begin{corollaire}
Supposons que
la s\'erie $\sum_{n\geq 1} \lambda_{k'-m-1}(F_n)[\lambda_{k'-m}(F_n)]^{-1}$
converge et que ${\cal H}$ soit $F_n$-r\'eguli\`ere pour
tout $n\geq 1$. Supposons aussi que $X'$ est une vari\'et\'e
projective. Alors la suite de courants 
$$\frac{1}{\lambda_{k'-m}(F_n)}
\Big((F_n)^*[\H_{y}]-(F_n)^*[\H_{y'}]\Big)$$
tend faiblement
vers $0$ pour $y$ et $y'$ hors d'un
sous-ensemble pluripolaire $\E$ de $Y$.
\end{corollaire}
\begin{preuve} Il suffit d'appliquer le th\'eor\`eme 4.6
et la
proposition 3.1 pour les transformations m\'eromorphes
$G\circ F_n$ (\voir aussi remarques 4.7). 
\end{preuve}
\begin{remarque} \rm
Lorsque $X$, $X'$ sont des espaces projectifs, $\H$ la
famille des sous-espaces projectifs de $X$ et 
$F_n$ des 
applications rationnelles de $X$ dans $X'$,
des versions du corollaire 4.10 sont prouv\'ees par
Russakovkii-Sodin \cite{RussakovskiiSodin}
et Russakovskii-Shiffman \cite{RussakovskiiShiffman}.
\end{remarque} 
\sectionbis{Mesures d'\'equilibre de correspondances}
Dans ce paragraphe, on suppose que les vari\'et\'es $X_n$ sont de
m\^eme dimension $k$. Nous \'etudions l'it\'eration
al\'eatoire d'une suite de correspondances
$f_n:X_{n-1}\longrightarrow X_n$. 
Notons $d_n$ le degr\'e 
topologique de $f_n$, pour $n\geq 1$.
Posons $\RR^*_n:=\RR^*_2(X_n,\omega_n,2)$,
$\AA_n:=\AA(f_n)$ et $\delta_n$
le degr\'e interm\'ediaire d'ordre $k-1$ de
$f_n\circ\cdots\circ f_1$ (\voir
d\'efinitions aux paragraphes 2.1 et 3.5). 
Soient $h_n$ des fonctions positives
dans
$\LL^2(X_n)$ telles que $\int_{X_n}h_n\omega_n^k=1$.
Posons $\nu_n:=h_n\omega_n^k$ et $\mu_n:=d_1^{-1}\ldots d_n^{-1}
(f_n\circ\cdots \circ f_1)^*(\nu_n)$.
Les mesures de probabilit\'e
$\nu_n$ \'etant absolument continues par
rapport aux mesures de Lebesgue,
les mesures de probabilit\'e $\mu_n$
le sont aussi.
On a le th\'eor\`eme suivant.
\begin{theoreme} Supposons que
$\delta_n\RR^*_n\|h_n\|_{\LL^2(X_n)}
=\o(d_1\ldots d_n)$
et que la s\'erie $\sum_{n\geq
2} d_1^{-1}\ldots d_n^{-1}\delta_{n-1}\AA_n$ converge.
Alors la suite
de mesures $\mu_n$ tend faiblement vers une mesure de
probabilit\'e PLB $\mu$ sur $X_0$. De plus, la mesure $\mu$ est
ind\'ependante de la suite $(h_n)$.
%
%
\end{theoreme}
\begin{remarques} \rm
Si les $(X_n,\omega_n,f_n)_{n\geq 0}$
appartiennent \`a une famille compacte lisse
(par exemple une famille
finie), 
les constantes $\RR^*_n$ et $\AA_n$
sont uniform\'ement born\'ees en $n$ (\voir remarque
2.3). Dans
ce cas, il suffit de supposer
que $\delta_n\|h_n\|_{\LL^2(X)}=\o(d_1\ldots d_n)$ et 
que la s\'erie $\sum d_1^{-1}\ldots
d_n^{-1}\delta_{n-1}$ converge. Si les $f_n$ sont des applications
rationnelles dominantes de degr\'e alg\'ebrique
$s_n$ de $\P^k$ dans $\P^k$, on peut
majorer $\delta_n$ par $(s_1\ldots s_n)^{k-1}$.
\par
Le th\'eor\`eme 5.1 est aussi valable pour les transformations
m\'eromorphes $f_n: X_{n-1}\longrightarrow X_n$ de codimension
$0$ (dans ce cas, on ne suppose pas que les $X_n$ ont la m\^eme
dimension). 
\end{remarques}
\begin{preuve}
Soit $\varphi$ une fonction q.p.s.h. sur $X_0$, $\ddc\varphi\geq
-\omega_0$.
Il nous suffit de montrer que la suite $\langle \mu_n,\varphi\rangle$
converge vers une constante $c_\varphi$ ind\'ependante de
$(h_n)$ (on posera alors $\langle\mu,\varphi\rangle := c_\varphi$
pour $\varphi$ continue). 
\par
Posons $F_n:=f_n\circ\cdots\circ f_1$, $T_0^-:=\omega_0$ et
$T_0^+:=\ddc\varphi +\omega_0$. Le courant
$T^+_0$ est positif ferm\'e et
cohomologue \`a $\omega_0$. 
On d\'efinit par r\'ecurrence les nombres $b_n$ et les
fonctions $\varphi_n$. Posons
$$b_0:=\int_{X_0}\varphi\omega_0^k \ \ \ \mbox{ et } \ \ \
\varphi_0:=\varphi-b_0.$$
D'apr\`es les propositions 2.6 et 2.7, on peut
poser, pour tout $n\geq 1$,
$$b_n:=\int_{X_n}(f_n)_*(\varphi_{n-1})\omega_n^k \ \ \ \mbox{ et }
\ \ \ 
\varphi_n:=(f_n)_*(\varphi_{n-1})-b_n.$$
On a
$$\ddc\varphi_n=T_n^+-T_n^-\  \mbox{ avec } \
T_n^\pm=(F_n)_*(T_0^\pm).$$
De plus
$$\cl(T_n^\pm)=\cl((F_n)_*(\omega_0)) \ \mbox{ et }\
\int_{X_n} \varphi_n\omega_n^k=0.$$
On en d\'eduit que $\|T_n^\pm\|= \|(F_n)_*(\omega_0)\|
=\delta_n$.
\\
\ \par
D'apr\`es la proposition 2.2, il existe
des fonctions q.p.s.h. $\varphi_n^\pm$ v\'erifiant
$\int_{X_n}\varphi_n^\pm\omega_n^k=0$ et 
$$\ddc\varphi_n^+- T^+_n=\ddc\varphi^-_n-T^-_n\geq
-r(X_n,\omega_n)\delta_n\omega_n.$$
Par d\'efinition de $\AA_n$ et $\RR^*_n$, on a
\begin{eqnarray}
\left|\int_{X_n}(F_n)_*(\varphi^\pm_{n-1})\omega_n^k\right|
\leq \delta_{n-1}\AA_n  & \ \ \ \mbox{ et }\ \ \ &
\|\varphi^\pm_n\|_{\Ltwo(X_n)}\leq \delta_n\RR^*_n
\end{eqnarray}
On en d\'eduit que $b_n\leq 2\delta_{n-1}\AA_n$. 
L'hypoth\`ese du th\'eor\`eme implique que la s\'erie
$\sum d_1^{-1}\ldots d_n^{-1} b_n$ converge. Notons $c_\varphi$
la somme de cette s\'erie. 
\\
\ \par
Dans la suite, on int\'egre seulement sur un ouvert de volume
total car les mesures sont absolument continues par rapport
aux mesures de Lebesgue. On a
\begin{eqnarray*}
\langle \mu_n,\varphi \rangle & = & \langle
d_1^{-1}\ldots d_n^{-1}f_1^*\ldots f_n^*(\nu_n),
b_0+\varphi_0 \rangle \\
& = & b_0 + \langle d_1^{-1}\ldots d_n^{-1}f_2^*\ldots
f_n^*(\nu_n),(f_1)_*(\varphi_0)\rangle\\
& = & b_0 + \langle d_1^{-1}\ldots d_n^{-1}f_2^*\ldots
f_n^*(\nu_n),b_1+\varphi_1\rangle\\
& = & b_0 +d_1^{-1}b_1 + 
\langle d_1^{-1}\ldots d_n^{-1}f_2^*\ldots
f_n^*(\nu_n),\varphi_1\rangle.
\end{eqnarray*}
Par r\'ecurrence, on obtient
\begin{eqnarray}
\langle \mu_n,\varphi \rangle & = & b_0+d_1^{-1}b_1+\cdots +
d_1^{-1}\ldots d_n^{-1}b_n +d_1^{-1}\ldots d_n^{-1}
\langle \nu_n, \varphi_n\rangle.\ \ \ \
\end{eqnarray}
V\'erifions que $\langle \mu_n,\varphi\rangle$ tend vers
$c_\varphi$. 
L'in\'egalit\'e de Cauchy-Schwarz et les relations (5.1) impliquent que
\begin{eqnarray*}
|\langle \nu_n, \varphi_n \rangle|  & = &  |\langle
h_n\omega_n^k,\varphi_n \rangle|
\leq  \|h_n\|_{\LL^2(X_n)}\|\varphi_n\|_{\LL^2(X_n)}\\
& \leq &  \|h_n\|_{\LL^2(X_n)}(\|\varphi^+_n\|_{\LL^2(X_n)}+
\|\varphi^+_n\|_{\LL^2(X_n)})\\
& \leq &  2\|h_n\|_{\LL^2(X_n)}\delta_n \RR^*_n.
\end{eqnarray*}
Par hypoth\`ese, la derni\`ere expression est
d'ordre $\o(d_1\ldots d_n)$.
Donc $\lim \langle \mu_n,\varphi\rangle =c_\varphi$.
\\
\ \par
D\'efinissons la mesure $\mu$ par
$$\langle\mu,\varphi\rangle:=c_\varphi \ \mbox{ pour } \varphi
\mbox{ lisse}.$$
On a montr\'e que $\mu_n$ tend faiblement vers $\mu$.
\\
\ \par
Si $\varphi$ est une fonction q.p.s.h. quelconque, par
semi-continuit\'e sup\'erieure, on a
$$\langle\mu,\varphi\rangle \geq \limsup \langle
\mu_n,\varphi\rangle =c_\varphi.$$
Puisque $\varphi$ est born\'ee sup\'erieurement, elle est
$\mu$-int\'egrable. Donc $\mu$ est PLB.
\end{preuve}
\par
Nous allons pr\'eciser notre r\'esultat pour l'it\'eration
d'une correspondance $f$ de degr\'e topologique $d_t$ 
de $X$ dans elle-m\^eme.
D'apr\`es la proposition A.1, les mesures PLB sur $X$
ne chargent pas les sous-ensembles analytiques. On peut donc d\'efinir
l'image de ces mesures par $f^*$. On a le corollaire suivant. 
\begin{corollaire} Soit $f$ une correspondance m\'eromorphe de
degr\'e topologique $d_t$ sur une vari\'et\'e k\"ahl\'erienne
compacte $(X,\omega)$. 
Supposons que le degr\'e dynamique d'ordre $k-1$
de $f$ v\'erifie $d_{k-1}<d_t$.
Soient $h_n$ des fonctions positives v\'erifiant $\int_X
h_n\omega^k=1$ et
$\|h_n\|_{\Ltwo(X)}^{1/n}=\o(d_{k-1}^{-1}d_t)$.
Alors la suite de
mesures $\mu_n:=d_t^{-n} (f^n)^*(h_n\omega^k)$
converge vers une mesure
PLB $\mu$ ind\'ependante de $(h_n)$. De
plus, $\mu$ v\'erifie la relation d'invariance
$f^*(\mu)=d_t\mu$ et on a $\langle \mu_n,
\varphi\rangle\rightarrow \langle \mu,\varphi\rangle$ pour
toute fonction q.p.s.h. $\varphi$ sur $X$. 
\end{corollaire}
\begin{preuve} La convergence de $(\mu_n)$
se d\'eduit du th\'eor\`eme 5.1.
Soit $\Omega$ une forme volume lisse telle que $\int_X\Omega=1$. 
Soit $\varphi$ une fonction lisse, on a
\begin{eqnarray*}
\langle f^*(\mu),\varphi \rangle & = &
\langle \mu, f_*(\varphi)\rangle = \lim_{n\rightarrow\infty}
\langle d_t^{-n}\Omega, (f^{n+1})_*(\varphi) \rangle \\
& = & d_t \langle d_t^{-n-1}\Omega, (f^{n+1})_*(\varphi) \rangle
=d_t \langle \mu,\varphi\rangle.
\end{eqnarray*}
D'o\`u la relation d'invariance.
\\
\ \par
Soit $\varphi$ une fonction q.p.s.h. sur $X$ avec
$\ddc\varphi\geq -\omega$. 
Montrons que $\langle\mu,\varphi\rangle= c_\varphi$.
Utilisons les notations du th\'eor\`eme 5.1.
Puisque $\mu$ est $f^*$-invariante, on montre comme pour la
relation (5.2) que
\begin{eqnarray}
\langle \mu,\varphi \rangle & = & b_0+d_t^{-1}b_1+\cdots +
d_t^{-n}b_n +d_t^{-n}
\langle \mu, \varphi_n^+-\varphi_n^-\rangle.
\end{eqnarray}
D'autre part, puisque $\int\varphi_n^\pm\omega^k=0$, d'apr\`es la relation (2.6), 
on a
$$|\langle\mu,\varphi^\pm_n\rangle|\leq
\lambda_{k-1}(f^n) \RR_3(X,\omega,\mu)$$
o\`u $\lambda_{k-1}(f^n)$ est le degr\'e interm\'ediaire d'ordre $k-1$ de $f^n$.
On en d\'eduit que le membre \`a droite de (5.3) tend vers
$c_\varphi$. D'o\`u $\langle \mu,\varphi\rangle =c_\varphi$.
\end{preuve}
\begin{remarques} \rm
Soit $f$ comme au corollaire 5.3. Notons $\E$ l'ensemble des
points $x\in X$ tels que la suite des mesures $\mu^x_n:=d_t^{-n}
(f^n)^*(\delta_x)$ ne converge pas vers
{\it la mesure d'\'equilibre} $\mu$ de $f$. On montre comme dans
\cite{Dinh} (\voir aussi \cite{BriendDuval, DinhSibony2}) que
$\E$ est une union finie ou d\'enombrable d'ensembles
analytiques et que les points p\'eriodiques r\'epulsifs de $f$ sont
denses dans le support de $\mu$. Si le nombre de points
p\'eriodiques r\'epulsifs est d'ordre $d_t^n+\o(d_t^n)$
alors ils sont \'equidistribu\'es sur
le support de $\mu$. Les arguments utilis\'es dans
\cite{DinhSibony2} permettent de traiter le cas des vari\'et\'es
k\"ahl\'eriennes non projectives.
\end{remarques}
\begin{remarques} \rm
Soit $f$ une application
m\'eromorphe dominante de degr\'e topologique $d_t$ 
de $X$ dans $X$ telle que son degr\'e dynamique d'ordre
$k-1$ v\'erifie $d_{k-1}<d_t$. Lorsque $X$ est projective, 
Vincent Guedj \cite{Guedj} a r\'ecemment construit
pour $f$ la mesure d'\'equilibre $\mu$. Il a montr\'e que $\mu$ est
PLB et m\'elangeante.
Dans \cite{DinhSibony3}, nous avons montr\'e que cette mesure est
d'entropie maximale $\log d_t$. 
V. Guedj utilise une m\'ethode de th\'eorie du
potentiel pour construire $\mu$.
Notre construction ci-dessus, par dualit\'e, 
donne dans ce cas deux informations
suppl\'ementaires. D'une part, la suite de mesures 
$d_t^{-n}(f^n)^*(h_n\omega^k)$ converge
vers $\mu$ uniform\'ement en $(h_n)$; on a aussi la convergence pour
toute fonction test q.p.s.h.
D'autre part, on
obtient une estimation de la vitesse de m\'elange.
\end{remarques}
\begin{theoreme} Soit $(X,\omega)$ une vari\'et\'e k\"ahl\'erienne
  compacte de dimension $k$. Soit $f$ une application m\'eromorphe
  dominante de
  degr\'e topologique $d_t$ de $X$ dans elle-m\^eme.
Supposons que son degr\'e dynamique
  d'ordre $k-1$ v\'erifie $d_{k-1}<d_t$. Alors la mesure d'\'equilibre
  $\mu$ de $f$ est m\'elangeante avec une vitesse de m\'elange
d'ordre $d_t^{-n}(d_{k-1}+\epsilon)^n$
  pour tout $\epsilon>0$. Plus pr\'ecis\'ement,
si $\varphi$ est une fonction de classe ${\cal C}^2$ et $\psi$
est une fonction born\'ee, posons 
$$I_n(\varphi,\psi):=\int_X \varphi(\psi\circ f^n) \d\mu - 
\left(\int_X\varphi\d\mu\right)\left(\int_X\psi\d\mu\right).$$
Il existe $c>0$ ind\'ependante de $\varphi$ et de $\psi$ telle que
$$|I_n(\varphi,\psi)|\leq cd_t^{-n}(d_{k-1}+\epsilon)^n
\|\varphi\|_{{\cal C}^2}\|\psi\|_{\Linfty(\mu)}.$$ 
\end{theoreme}
\begin{preuve} 
Puisque $\varphi$ s'\'ecrit comme diff\'erence de deux fonctions
q.p.s.h., on peut supposer que $\varphi$ est q.p.s.h. avec
$\ddc\varphi\geq -\omega$.
Du fait que $I_n(\varphi,\psi)=-I_n(\varphi,-\psi)$, il
  suffit de majorer $I_n(\varphi,\psi)$. 
Comme $I_n(\varphi,\psi+A)=I_n(\varphi,\psi)$ pour toute constante $A$, 
on peut supposer que $\psi$ est positive.
On peut \'egalement supposer que $\|\psi\|_{\Linfty(\mu)}=1$. 
Posons  $c_\varphi:=\langle\mu,\varphi\rangle$.
Comme $\mu$ est invariante,
  on a 
$$I_n(\varphi,\psi)=\int_X
\Big(d_t^{-n}(f^n)_*(\varphi)-c_\varphi\Big)
\psi\d\mu\leq \|d_t^{-n}(f^n)_*(\varphi) -c_\varphi\|_{\Lone(\mu)}.$$
On reprend les calculs d\'ej\`a faits au th\'eor\`eme 5.1. On a
$$d_t^{-n}(f^n)_*(\varphi) -c_\varphi=d_t^{-n}\varphi_n^+-
d_t^{-n}\varphi_n^-
-\sum_{i\geq n+1} d_t^{-i}b_i.$$
Comme au th\'eor\`eme 5.1, on a
$$|b_i|\leq
c_1(d_{k-1}+\epsilon)^{i-1}\ \mbox{ et }\
\|\varphi_n^\pm\|_{\Ltwo(X)}\leq c_1(d_{k-1}+\epsilon)^n$$
pour une constante $c_1>0$.
La mesure $\mu$ \'etant PLB, d'apr\`es la proposition 2.4,
il existe une constante $c_2>0$ telle que
$\|\varphi^\pm_n\|_{\Lone(\mu)}\leq c_2\|\varphi^\pm_n\|_{\Ltwo(X)}$.
On en d\'eduit que
$\|d_t^{-n}(f^n)_*(\varphi) -c_\varphi\|_{\Lone(\mu)}\leq
c(d_{k-1}+\epsilon)^nd_t^{-n}$ pour une constante $c>0$.
Ceci termine
la preuve du th\'eor\`eme.
\end{preuve}
\begin{remarque}\rm
Si $f$ est une correspondance, la mesure $\mu$ n'est pas $f_*$ invariante 
en g\'en\'eral. On a cependant montr\'e dans ce cas que
$$\|d_t^{-n}(f^n)_*\varphi -c_\varphi\|_{\Lone(\mu)} \leq cd_t^{-n}(d_{k-1}+\epsilon)^n
\|\varphi\|_{{\cal C}^2}.$$ 
\end{remarque}
\begin{exemple}\rm
Nous avons montr\'e dans \cite{DinhSibony3} qu'\'etant donn\'e
une correspondance 
$f$ d'une vari\'et\'e projective $X$, la suite
$[\lambda_l(f^n)]^{1/n}$ converge 
vers sa borne inf\'erieure $d_l(f)=\inf_n[\lambda_l(f^n)]^{1/n}$.
Il en r\'esulte que pour v\'erifier 
l'hypoth\`ese du corollaire 5.6, il suffit de
montrer que $\lambda_{k-1}(f)<d_t(f)$. 
On peut donc exhiber pour toute vari\'et\'e
projective $X$ des correspondances
v\'erifiant cette derni\`ere in\'egalit\'e.
Soient $h$ et $g$ deux projections
holomorphes surjectives de $X$ sur $\P^k$. Soit $u$ un
endomorphisme holomorphe de degr\'e 
\'elev\'e de $\P^k$. Posons
$f:=\overline g\circ u\circ h$ o\`u $\overline g$ est
l'adjoint de $g$. C'est une correspondance sur $X$. 
Le lemme 2.8 montre
(\voir aussi \cite[remarques 8]{DinhSibony3}) que
$$\lambda_{k-1}(f)\leq c\lambda_{k-1}(\overline g)\lambda_{k-1}(u)\lambda_{k-1}(h)$$
pour une certaine constante $c>0$.
Si le degr\'e de $u$ est suffisamment
\'elev\'e, on a $\lambda_{k-1}(f)<d_t(f)$. 
En effet, $d_t(f)=d_t(\overline g)d_t(u)d_t(h)$ et 
$d_t(u)>>\lambda_{k-1}(u)$ si le degr\'e
de $u$ est suffisamment \'elev\'e. On trouve dans \cite{DinhSibony3} quelques 
exemples d'applications rationnelles v\'erifiant l'hypoth\`ese du corollaire 5.3. 
\end{exemple}
\sectionbis{Distribution des intersections de vari\'et\'es}
Consid\'erons des transformations
m\'eromorphes $F_{1,n}:X\longrightarrow X_1$ et
$F_{2,n}:X\longrightarrow X_2$
de codimensions respectives $l_1$, $l_2$.
On suppose que $l_1+l_2\geq k$ et que
les produits de $F_{1,n}$ et $F_{2,n}$ sont bien d\'efinis
(\voir paragraphe 3.3). Posons
$\Phi_n:=F_{1,n}\bullet F_{2,n}$,
$\delta_{i,n}:=\lambda_{k_i-1}(F_{i,n})$ et 
$d_{i,n}:=\lambda_{k_i}(F_{i,n})$. Rappelons que $\Phi_n$ est une
transformation m\'eromorphe de codimension $l_1+l_2-k$
de $X$ dans $X_1\times X_2$. 
\begin{theoreme} 
Supposons que la vari\'et\'e $X$ soit projective et que
les s\'eries $\sum
\delta_{i,n}d_{i,n}^{-1}$ soient convergentes pour $i=1,2$.
Alors il
existe un sous-ensemble pluripolaire $\E$ de $X_1\times X_2$ tel
que pour tout $(x_1,x_2)\in (X_1\times X_2)\setminus\E$, la suite
de courants
$$\frac{1}{d_{1,n}d_{2,n}}
\Big(F_{1,n}^*(\delta_{x_1})\wedge F_{2,n}^*(\delta_{x_2}) -
F_{1,n}^*(\omega_1^{k_1})\wedge F_{2,n}^*(\omega_2^{k_2})\Big)$$
tend faiblement vers $0$.
\end{theoreme}
\begin{preuve} 
Soit $\delta_n$ le degr\'e interm\'ediaire d'ordre $k_1+k_2-1$ de $\Phi_n$. 
D'apr\`es, la proposition 3.2, il existe une constante $c>0$ telle que
$$\delta_n\leq c(\delta_{1,n}d_{2,n} + \delta_{2,n}d_{1,n}).$$ 
Le degr\'e interm\'ediaire d'ordre $k_1+k_2$ de $\Phi_n$ est
\'egal \`a $d_{1,n}d_{2,n}$. Il suffit d'appliquer
le th\'eor\`eme 4.6 pour $\Phi_n$.
\end{preuve}
\par
Soit $f$ une application bim\'eromorphe de $X$ dans
$X$ et soit $f^{-1}$ son inverse. Soient $\H^+$, $\H^-$ deux
familles m\'eromorphes
adapt\'ees r\'eguli\`eres de sous-ensembles analytiques
de dimensions respectives $k-l^+$ et $k-l^-$
avec $l^++l^-\leq k$. Notons $P^\pm:X\longrightarrow Y^\pm$ les
transformations m\'eromorphes associ\'ees,
$d_n^\pm$ les degr\'es dynamiques d'ordre $l^\pm$ de $f^{\pm
n}$ et $\delta_n^\pm$ les degr\'es dynamiques d'ordre $l^\pm-1$
de $f^{\pm n}$.
\begin{corollaire}
Supposons que la vari\'et\'e $X$ soit projective et
que les s\'eries $\sum \delta_n^\pm[d_n^\pm]^{-1}$
soient convergentes. Alors il existe un
sous-ensemble pluripolaire $\E$ de $Y^+\times Y^-$ tel que
pour $(a_1,b_1)$ et $(a_2,b_2)$ dans
$(Y^+\times Y^-)\setminus \E$, la suite
de courants
$$\frac{1}{d_n^+d_n^-}\Big([f^{-n}(\H^+_{a_1})\cap f^n(\H^-_{b_1})]
-[f^{-n}(\H^+_{a_2})\cap f^n(\H^-_{b_2})
]\Big)$$
tend faiblement vers $0$.
\end{corollaire}
\begin{preuve}
Posons $F^\pm_n:=P^\pm\circ f^{\pm n}$. Ce sont des
transformations m\'eromorphes de codimension $k-l^\pm$ de
$X$ dans $Y^\pm$.
Notons $\delta_n$ et $d_n$
leurs degr\'es interm\'ediaire
d'ordre $\dim Y^++\dim Y^--1$ et d'ordre $\dim Y^++\dim Y^-$. 
D'apr\`es les propositions 3.1 et 3.2, il existe $c>0$ telle que
$$\delta_n\leq c(\delta_n^+d_n^-+\delta_n^- d^+) \ \mbox{ et } \
d_n\leq cd_n^+d_n^-.$$
Il suffit d'appliquer le th\'eor\`eme 6.1 (\voir aussi
remarques 4.7). 
\end{preuve}
\par
Nous allons expliciter ce r\'esultat dans le cadre des
automorphismes r\'eguliers de $\C^k$ introduits par le second
auteur \cite{Sibony2}. Soit $f$ un
automorphisme polynomial de $\C^k$. On note aussi $f$ son prolongement
en application birationnelle de $\P^k$ dans $\P^k$. Notons $I^+$
(resp. $I^-$) l'ensemble d'ind\'etermination de $f$ (resp. de
$f^{-1}$). Ce sont des sous-ensembles analytiques de l'hyperplan
\`a l'infini. L'automorphisme $f$ est dit {\it r\'egulier} si
$I^+\cap I^-=\emptyset$ (en dimension 2, les automorphismes
r\'eguliers sont ceux du type H\'enon). On a 
alors $\dim I^++\dim I^-=k-2$. Posons
$s:=\dim I^++1$. Notons $d_+$ et $d_-$ les degr\'es alg\'ebriques
de $f$ et $f^{-1}$. Ces degr\'es sont li\'es par la relation
$(d_+)^s=(d_-)^{k-s}$. 
\par
On peut construire deux courants $T^+$, $T^-$
positifs ferm\'es de bidegr\'e $(1,1)$ de masse 1 de $\P^k$, \`a
potentiel continu dans $\C^k$ et tels que $f^*(T^+)=d_+T^+$,
$f_*(T_-)=d_-T^-$. Pour $0\leq l\leq s$ et $0\leq l'\leq
k-s$, le courant $T_{l,l'}:=
(T^+)^l\wedge (T^-)^{l'}$ est bien d\'efini. 
Quand $l=s$, $l'=k-s$, on obtient une mesure de
probabilit\'e invariante \`a support compact dans $\C^k$. On a
\'egalement le th\'eor\`eme de convergence suivant:
\begin{eqnarray}
\lim_{n,m\rightarrow\infty} (d_+)^{-nl}(d_-)^{-ml'}
(f^n)^*(\omega_\FS^l)\wedge
(f^m)_*(\omega_\FS^{l'}) & = & T_{l,l'}.
\end{eqnarray}
\par
Notons $\GG_l$ et $\GG_{l'}$ les grassmanniennes qui
param\`etrent les sous-espaces projectifs de dimension $k-l$ et $k-l'$
de $\P^k$. Notons $\P^{k-l}_x$ et $\P^{k-l'}_{x'}$ les sous-espaces
projectifs associ\'es aux points
$x\in\GG_l$ et $x'\in\GG_{l'}$.
Nous avons le th\'eor\`eme suivant.
\begin{theoreme} Soit $f$ un automorphisme r\'egulier de $\C^k$
comme ci-dessus. Il existe un sous-ensemble pluripolaire
$\E$ de $\GG_l\times \GG_{l'}$ tel que pour tout
$(x,x')\in(\GG_l\times \GG_{l'})\setminus\E$ la suite de courants
$$(d_+)^{-nl}(d_-)^{-ml'} [f^{-n}(\P^{k-l}_x)\cap f^m(\P^{k-l'}_{x'})]$$ 
tend faiblement vers le courant invariant $T_{l,l'}$ quand $n$ et $m$
tendent vers l'infini.
\end{theoreme}
\begin{preuve} 
Avec les notations du corollaire 6.2, on a $d_n^+=d_+^{ln}$, 
$\delta_n^+=d_+^{(l-1)n}$ et
des relations semblables pour les inverses.
Il suffit d'appliquer le corollaire 6.2 apr\`es avoir int\'egr\'e 
par rapport aux variables $a_2$, $b_2$. On utilise ensuite la relation (6.1). 
\end{preuve}
\begin{remarque} \rm
On peut obtenir des r\'esultats analogues pour les
familles de sous-ensembles analytiques de degr\'e quelconque ou pour des applications
rationnelles ou birationnelles
$f$ d'une vari\'et\'e projective $X$ dans elle-m\^eme.
\end{remarque}
\sectionbis{Z\'eros des sections de fibr\'es
en droites}
Soit $X$ une vari\'et\'e projective de dimension $k$
et soit $L$ un
fibr\'e en droites ample sur $X$. On munit $L$ d'une
m\'etrique hermitienne $h$. Pour toute section holomorphe
locale $e_L$ de $L$, on d\'efinit la norme de $e_L$ en chaque
point par
$\|e_L\|_h:=h(e_L,e_L)^{1/2}$. Soit 
$$c_1(h):=-\ddc \log\|e_L\|_h$$
la forme de courbure de $(L,h)$.
Elle repr\'esente dans la cohomologie de Rham 
la classe de Chern $c_1(L)\in \H^2(X,\Z)$ de $L$.
Puisque $L$ est ample, on peut choisir $h$ de sorte que $c_1(h)$
soit une $(1,1)$-forme strictement positive.
La vari\'et\'e $X$ est donc munie
de la forme de K\"ahler $\omega:=c_1(h)$ et
$\int_X\omega^k=c_1(L)^k\in\Z^+$.
\par
Le fibr\'e $L^n$ est \'egalement muni d'une m\'etrique
hermitienne $h_n$, induite par la m\'etrique $h$ sur $L$. 
Plus pr\'ecis\'ement, $h_n$ est d\'efinie localement par
$\|s^n\|_{h_n}=\|s\|_h^n$.
L'espace
$\HH^0(X,L^n)$ des sections holomorphes de $L^n$
est muni du produit hermitien naturel
$$\langle s_1,s_2 \rangle :=\frac{1}{c_1(L)^k}
\int_Xh_n(s_1,s_2)\omega^k\ \ \ \
(s_1,s_2\in \HH^0(M,L^n)).$$
Notons $\omega_\FS$ la m\'etrique de Fubini-Study de
$\P\HH^0(X,L^n)$.
\par
Le lecteur trouvera les autres notations dans l'exemple 3.6(c).
La dimension $k_n$ de $\P\HH^0(X,L^n)$ est donn\'ee par le
polyn\^ome de Hilbert dont le terme dominant est \'egal \`a 
$c_1(L)^kn^k/k!$ \cite[p.386]{Kobayashi}. 
Rappelons que pour $n$ assez grand, $\Phi_n$ est le plongement de Kodaira
de $X$ dans $\P\HH^0(X,L^n)^*$, $\Psi_{l,n}$ la transformation
m\'eromorphe naturelle de $\P\HH(X,L^n)^*$ dans la grassmannienne
$\GG_{l,n}^X$ des sous-espaces de dimension $l$ de 
$\P\HH^0(X,L^n)$,
$\overline\Pi_{l,n}$ est la transformation m\'eromorphe naturelle
de $\GG_{l,n}^X$ dans
$\P^X_{l,n}:=\P\HH^0(X,L^n)\times\cdots\times \P\HH^0(X,L^n)$ 
($l$ fois) et $F_{l,n}=\overline \Pi_{l,n}\circ 
\Psi_{l,n}\circ \Phi_n$.
\begin{lemme}
Soient $\delta_{l,n}$ et $d_{l,n}$
les degr\'es interm\'ediaires d'ordre $lk_{n}-1$ et
d'ordre $lk_n$ de $F_{l,n}$. On a
$d_{l,n}=n^lc_1(L)^k$ et $\delta_{l,n}=\O(n^{l-1})$.
\end{lemme}
\begin{preuve}
L'invariance des m\'etriques par l'action du groupe unitaire
implique que
\begin{eqnarray}
\Psi_{l,n}^*\overline\Pi_{l,n}^*(\omega_\MP^{lk_n})=\alpha_{l,n}
\omega_\FS^l
&\  \mbox{ et }\   &
\Psi_{l,n}^*\overline\Pi_{l,n}^*(\omega_\MP^{lk_n-1})=
\beta_{l,n}\omega_\FS^{l-1}\ 
\end{eqnarray}
o\`u $\omega_\MP$ est la forme k\"ahl\'erienne naturelle
associ\'ee \`a $\P^X_{l,n}$ (\voir appendice A.3) et 
$\alpha_{l,n}$, $\beta_{l,n}$ sont des constantes positives.
On calcule ces constantes cohomologiquement.
\par
Pour calculer $\alpha_{l,n}$, on
remplace $\omega_\MP^{lk_n}$ dans (7.1) par une masse de
Dirac $\delta_\ss$. Son image
$\Psi_{l,n}^*\overline\Pi_{l,n}^*(\delta_\ss)$ sera le courant
d'int\'egration sur un sous-espace projectif de codimension $l$
de $\P\HH(X,L^n)^*$ qui est cohomologue \`a $\omega_\FS^l$. Ceci
implique que $\alpha_{l,n}=1$.
\par
Soit $T$ le courant
d'int\'egration sur une droite ${\cal D}\times \{s_2\}
\times\cdots\times
 \{s_l\}$ de $\P^X_{l,n}$. Il est de masse $c_{k_n,l}$
(\voir (A.3) pour la d\'efinition de $c_{k_n,l}$).
Son image $\Psi_{l,n}^*\overline\Pi_{l,n}^*(T)$ est
le courant d'int\'egration sur un sous-espace
projectif de codimension $l-1$
de $\P\HH(X,L^n)^*$. La masse de ce dernier courant
est \'egale \`a 1.
On en d\'eduit que 
$\beta_{l,n}=c_{k_n,l}^{-1}$. En particulier, il est major\'e
par une constante qui ne d\'epend que de $l$ (\voir (A.3)).
\\
\ \par
Puisque la classe de $\Phi_n^*(\omega_\FS)$ est \'egale \`a
$nc_1(L)$, on a 
$$d_{l,n}=\int_X
\Phi_n^*\Psi_n^*\overline\Pi_{l,n}^*(\omega_\MP^{lk_n})
\wedge \omega^{k-l}
=\int_X \Phi_n^*(\omega_\FS^l)\wedge \omega^{k-l}
=n^lc_1(L)^k$$
et
\begin{eqnarray*}
\delta_{l,n} & = & \int_X
\Phi_n^*\Psi_n^*\overline\Pi_{l,n}^*(\omega_\MP^{lk_n-1})\wedge \omega^{k-l+1}
=\int_X \beta_{l,n}\Phi_n^*(\omega_\FS^{l-1})\wedge \omega^{k-l+1}\\
 & = & \beta_{l,n}n^{l-1}c_1(L)^k.
\end{eqnarray*}
\end{preuve}
On voit que
la s\'erie $\sum \delta_{l,n}d_{l,n}^{-1}$ ne converge pas.
C'est donc le th\'eor\`eme 4.1(2) que nous appliquerons.
\\
\ \par
Le th\'eor\`eme suivant, d\^u \`a Zelditch \cite{Zelditch},
est une
am\'elioration d'un th\'eor\`eme de Tian \cite{Tian}, il donne la convergence 
en moyenne des courants
$n^{-l}[\ZZ_{\ss_n}]$ avec $\ss_n\in\P^X_{l,n}$.
\begin{theoreme}{\bf \cite{Zelditch}}
Pour tout $r\geq 0$, on a
$$\|n^{-l}\Phi_n^*(\omega_\FS^l)-\omega^l\|_{{\cal
C}^r}=\O(n^{-1}).$$
\end{theoreme}
\par
Soient $\sigma_n$ des mesures de probabilit\'e
PLB sur $\P_{l,n}^X$. Posons
$\RR_{1,n}:=\RR_1(\P_{l,n}^X,\omega_\MP,\sigma_n)$,
$\RR_{3,n}:=\RR_3(\P_{l,n}^X,\omega_\MP,\sigma_n)$ et 
$\Delta_n(t):=\Delta(\P_{l,n}^X,\omega_\MP,\sigma_n,t)$.
On munit $\P^X_l:=\prod_{n\geq 1} \P_{l,n}^X$ de la mesure $\sigma$, 
produit des
$\sigma_n$. Faisons des hypoth\`eses sur les mesures $\sigma_n$. 
\begin{theoreme} Supposons que la s\'erie
$\sum_{n\geq 1}\Delta_n(nt)$ converge pour tout $t>0$. Supposons
aussi que
$\RR_{1,n}=\o(n)$ et $\RR_{3,n}=\o(n)$.
Alors pour $\sigma$-presque tout $\ss=(\ss_n)\in \P^X_l$,
la suite de courants $n^{-l}[\ZZ_{\ss_n}]$ tend faiblement
vers $\omega^l$. 
\end{theoreme}
\begin{preuve} Les relations (7.1) et le th\'eor\`eme 7.2
entra\^{\i}nent que 
la suite des formes
$n^{-l}F_{l,n}^*(\omega_\MP^{lk_n})$ tend vers $\omega^l$
dans ${\cal C}^r$ pour tout $r\geq 0$. D'apr\`es la proposition
4.3, la suite de courants
$n^{-l}F_{l,n}^*(\sigma_n)$ tend faiblement vers $\omega^l$ car
$n^{-l}\delta_{l,n}\RR_{3,n}$ tend vers $0$ par hypoth\`ese.
D'apr\`es le th\'eor\`eme 4.1(2), pour $\sigma$-presque
tout $\ss\in \P^X_l$, la suite de courants
$n^{-l}([\ZZ_{\ss_n}]-F_{l,n}^*(\sigma_n))$
tend faiblement vers $0$.
\end{preuve}
\begin{remarque}\rm
Si au th\'eor\`eme 7.3, on suppose que $\RR_{3,n}=\O(\log n)$
alors il existe des constantes $c>0$, $\alpha>0$ et $m\geq 0$
telles que pour tout $\epsilon>0$
$$\sigma_n\big(\ss\in\P^X_l,\ |\langle
n^{-l}[\ZZ_{\ss_n}]-\omega^l,\psi\rangle|\geq \epsilon\big) \leq
c \|\psi\|_{{\cal C}^2} n^{mk}\exp(-\alpha \epsilon n).$$
Cette estimation r\'esulte de l'in\'egalit\'e (4.5), du
th\'eor\`eme 7.2 et de la proposition 4.2(c).
\end{remarque}
\par
Posons $\GG_l^X:=\prod_{n\geq 1}\GG^X_{l,n}$. Soit $\Omega_{l,n}$ la mesure
de probabilit\'e invariante sur $\GG^X_{l,n}$. Notons $\Omega_l$ la mesure
de probabilit\'e sur $\GG^X_l$, produit des $\Omega_{l,n}$. 
Fixons des sous-espaces r\'eels $\R\P\HH^0(X,L^n)$
de $\P\HH^0(X,L^n)$ 
invariants par l'action du groupe orthogonal associ\'e. 
Soit $\R\GG_{l,n}^X$ la sous-grassmannienne totalement r\'eelle
de $\GG_{l,n}^X$ correspondante.
Soit $m_{l,n}$ la mesure invariante de masse 1 sur $\R\GG_{l,n}^X$. 
Notons $m_l$ le
produit des $m_{l,n}$ qui est une mesure de probabilit\'e sur
$\GG_l$. 
\begin{corollaire} Soient $\mu_{l,n}:=\Omega_{l,n}$ ou $m_{l,n}$ et 
$\mu_l:=\Omega_l$ ou $m_l$. Alors
pour $\mu_l$-presque tout $\scheck=(\scheck_n)\in \GG_l$,
la suite de courants $n^{-l}[\ZZ_{\scheck_n}]$
tend faiblement vers $\omega^l$.
De plus, on a
$$\mu_{l,n}\big(\scheck\in\GG^X_l,\ |\langle
n^{-l}[\ZZ_{\scheck_n}]-\omega^l,\psi\rangle|\geq \epsilon\big) \leq
c \|\psi\|_{{\cal C}^2} n^{mk}\exp(-\alpha \epsilon n)$$
o\`u $m\geq 0$, $\alpha>0$ et $c>0$ sont des constantes ind\'ependantes
de $\epsilon$.  
\end{corollaire} 
\begin{preuve}
Notons $\R\P^X_{l,n}:= \R\P\HH^0(X,L^n)\times\cdots \times
\R\P\HH^0(X,L^n)$ ($l$ fois).
Notons $\widetilde\Omega_{l,n}$ (resp. $\widetilde m_{l,n}$) la mesure de
probabilit\'e invariante naturelle sur $\P^X_{k,l}$ (resp.
sur $\R\P^X_{l,n}$) (\voir appendice pour les d\'etails)
et $\widetilde\Omega_l$ (resp. $\widetilde m_l$) la mesure de
probabilit\'e sur $\P^X_l$, produit des $\widetilde\Omega_{l,n}$
(resp. des $\widetilde m_{l,n}$). L'invariance des mesures consid\'er\'ees
implique que $\overline \Pi_{l,n}^*(\widetilde\Omega_{l,n})=\Omega_{l,n}$
et $\overline \Pi_{l,n}^*(\widetilde m_{l,n})= m_{l,n}$.
\par
Soient $\widetilde\mu_{l,n}=\widetilde\Omega_{l,n}$ ou $\widetilde m_{l,n}$ et
$\widetilde\mu_l=\widetilde\Omega_l$ ou $\widetilde m_l$.
Il suffit de montrer que pour $\widetilde\mu$-presque tout
$\ss=(\ss_1,\ss_2,\ldots) \in\P^X_l$ on a
$n^{-l}[\ZZ_{\ss_n}]\rightarrow \omega^l$ et 
$$\widetilde \mu_{l,n}\big(\ss\in\P^X_l,\ |\langle
n^{-l}[\ZZ_{\ss_n}]-\omega^l,\psi\rangle|\geq \epsilon\big) \leq
c \|\psi\|_{{\cal C}^2} n^{mk}\exp(-\alpha \epsilon n).$$
Cela r\'esulte du th\'eor\`eme 7.3 et la remarque 7.4.
Les estimations sur
$\RR_{1,n}$, $\RR_{3,n}$ et $\Delta_n$ sont fournies par la
proposition A.11. La dimension de $\P^X_{l,n}$ est de l'ordre
$n^k$. 
\end{preuve}
\begin{remarques}\rm
Soit $(c_n)$ une suite de nombre r\'eels positifs v\'erifiant 
$c_n=\o(n/\log n)$. On peut montrer que
$c_n(n^{-l}[\ZZ_{\scheck}]-\omega^l)$ 
tend vers $0$ pour $\mu_l$-presque tout $\scheck$
(\voir aussi remarques 4.7). Ceci montre que
$n^{-l}[\ZZ_{\scheck}]-\omega^l$ tend vers $0$ \`a vitesse
$\simeq \log n/n$. Obsevons que la multiplication par $c_n$ revient \`a diviser 
$\epsilon$ par $c_n$. 
\par
Dans le cas o\`u $l=1$, 
Shiffman et Zelditch \cite{ShiffmanZelditch1}
ont d\'emontr\'e le r\'esultat 
de convergence $n^{-1}[\ZZ_{s_n}]\rightarrow \omega$ pour 
$\Omega_1$-presque toute suite $s=(s_n)$, $s_n\in\P\HH^0(X,L^n)$. Il ont
prouv\'e une vitesse de convergence major\'ee par $1/\sqrt{n}$.
\end{remarques}
\begin{corollaire} Soit $U$ un ouvert de $X$ dont le bord $\partial U$
est de mesure nulle. 
Alors pour
$\Omega_l$-presque tout (resp. $m_l$-presque tout) $\scheck=(\scheck_n)
\in\GG^X_l$ on a
$$\lim_{n\rightarrow\infty} n^{-l}\vol_{2k-2l}(\ZZ_{\scheck_n}\cap U)
=\frac{k!}{(k-l)!}\vol_{2k}(U).$$
\end{corollaire}
La d\'emonstration est laiss\'ee au lecteur.
\begin{remarque}\rm
Soit $(c_n)$ une suite de nombre r\'eels positifs v\'erifiant 
$c_n=\o(n/\log n)$. On peut montrer (\voir remarques 7.6 et 4.7) que
$$\lim_{n\rightarrow\infty} c_n\Big[n^{-l}\vol_{2k-2l}(\ZZ_{\scheck_n}\cap U)
-\frac{k!}{(k-l)!}\vol_{2k}(U)\Big]=0.$$
\end{remarque}
\appendix {Appendice: estimations des constantes}
{\bf A.1}. {\it Ensembles pluripolaires, capacit\'es, mesures 
mod\'er\'ees}.
\\
\ \par
Un sous-ensemble $E$ de $X$ est {\it pluripolaire} si
$E\subset(\varphi=-\infty)$ o\`u $\varphi$ est une fonction
q.p.s.h. On peut supposer $\ddc\varphi\geq -\omega$. 
Observons qu'une r\'eunion d\'enombrable d'ensembles pluripolaires est pluripolaire.
On dit que
$E$ est {\it localement pluripolaire} si pour tout $a\in X$ il
existe un voisinage $U$ de $a$ tel que $E\cap U$ soit
pluripolaire dans $U$.
Nous renvoyons \`a Demailly \cite{Demailly2}
pour les propri\'et\'es des fonctions
q.p.s.h. 
\par
Josefson \cite{Josefson} a montr\'e qu'un sous-ensemble
localement pluripolaire dans $\C^k$ est pluripolaire. Alexander a
\'etendu ce r\'esultat \`a $\P^k$ \cite{Alexander}.
Il en r\'esulte que pour toute vari\'et\'e projective $X$
de dimension $k$, un
ensemble $E$ localement pluripolaire l'est globalement: il suffit
d'utiliser une application holomorphe surjective
$\pi:X\longrightarrow \P^k$. Si $\varphi$ est associ\'ee \`a
$\pi(E)$, $\pi^*\varphi$ est associ\'ee \`a $E$. Il serait utile
de montrer ce r\'esultat pour toute vari\'et\'e complexe
compacte. On utilise \`a plusieurs reprises le r\'esultat suivant.
\begin{proposition} Soit $(X,\omega)$ une vari\'et\'e
k\"ahl\'erienne compacte de dimension $k$. Tout sous-ensemble
analytique propre $Y\subset X$ est pluripolaire.
\end{proposition}
\begin{preuve}
On peut supposer que $Y$ est irr\'eductible.
Si $Y$ est une hypersurface de $X$, il existe une fonction
q.p.s.h. $\varphi$ telle que $\ddc\varphi=[Y]-\alpha$ o\`u
$\alpha$ est une $(1,1)$-forme lisse cohomologue \`a $[Y]$. Il
est clair que $Y=(\varphi=-\infty)$. La fonction $\varphi$ est lisse sur 
$X\setminus Y$.
\par
Si $\dim Y<k-1$, on construit une vari\'et\'e lisse $\widehat X$
par des \'eclatements successifs le long $Y$ et le long des singularit\'es de $Y$.
Notons $\pi$ la projection de $\widehat X$ sur $X$. L'ensemble
$\pi^{-1}(Y)$ est alors une hypersurface de $\widehat X$.  
D'apr\`es Blanchard \cite{Blanchard} $\widehat X$ est k\"ahl\'erienne. 
\par
Soit $\psi$ une fonction q.p.s.h. sur $\widehat X$, lisse sur $\widehat X\setminus 
\pi^{-1}(Y)$ telle que
$\pi^{-1}(Y)=(\psi=-\infty)$. Puisque $\ddc\psi$ s'\'ecrit comme diff\'erence 
de deux courants positifs ferm\'es, lisses sur $\widehat{X}\setminus \pi^{-1}(Y)$, 
on a $\pi_*\psi=u_1-u_2$ dans $\Lone(X)$ avec $u_i$
q.p.s.h. lisses sur $X\setminus Y$.  
On en d\'eduit que $u_1-u_2=\pi_*(\psi)$ sur $X\setminus Y$. Puisque $\pi_*\psi(x)$
tend vers $-\infty$ quand $x$ tend vers $Y$, on a 
$Y\subset (u_1=-\infty)$ car $u_2$ est born\'ee sup\'erieurement.   
\end{preuve}
\par
Bien que ce ne soit pas indispensable, introduisons une capacit\'e dont 
les ensembles de capacit\'e nulle sont les ensembles pluripolaires. 
Cela a \'et\'e fait par
Alexander pour
$\P^k$. Notons $\pi:\C^{k+1}\setminus\{0\}\longrightarrow \P^k$ la projection
canonique. Soient $\S^{2k+1}$ la sph\`ere unit\'e de $\C^{k+1}$ et $\sigma_{2k+1}$ 
la mesure de probabilit\'e invariante sur $\S^{2k+1}$.
Alexander a pos\'e pour un compact $K$
de $\P^k$:
\begin{eqnarray}
\capacity'(K) & := & \inf_f\Big\{\sup_{\pi^{-1}(K)\cap\S^{2k+1}} 
|f|^{1/n},\
 f \mbox{ polyn\^ome
homog\`ene de degr\'e } n\nonumber \\
& & 
\ \ \ \ \mbox{ de } \C^{k+1},\ 
\int_{\S^{2k+1}} \big(\log|f|^{1/n}-\log|z_1|\big) \d\sigma_{2k+1}=0
\Big\}
\end{eqnarray}
Etant donn\'e un compact $K$ dans une vari\'et\'e k\"ahl\'erienne
compacte $(X,\omega)$
nous d\'efinissons {\it la capacit\'e de
$K$} par 
$$\capacity(K):=\inf_\varphi\left\{\exp\Big(\sup_K\varphi\Big),
\ \varphi \mbox{ q.p.s.h., } \ddc\varphi\geq -\omega,
\max_X\varphi=0\right\}.$$
Dans $\P^k$
toute fonction q.p.s.h. $\varphi$, v\'erifiant $\ddc\varphi\geq
-\omega_\FS$, est limite de fonctions sur $\P^k$ de la forme
$\log|f|^{1/n}-\log\|z\|$ o\`u $f$ est homog\`ene de degr\'e
$n$. En utilisant la proposition A.5 ci-dessous
on peut montrer que
$$\capacity'(K)\leq \capacity(K)\leq \sqrt{k\mbox{e}}\ \capacity'(K).$$
\par
Dans la suite nous utilisons plut\^ot la capacit\'e $\capacity$
qui a un sens pour toute vari\'et\'e k\"ahl\'erienne compacte.
Elle peut \^etre d\'efinie pour toute vari\'et\'e complexe compacte
hermitienne, en rempla\c cant $\omega$ par une forme hermitienne
positive.
Avec notre normalisation, on a toujours $\capacity(X)=1$. 
\begin{proposition} On a $\capacity(K)=0$ si et seulement si $K$
est pluripolaire.
\end{proposition}
\begin{preuve} Il est clair que si $K$ est pluripolaire on a
$\capacity(K)=0$. 
Si $\capacity(K)=0$, il existe $\varphi_n$ q.p.s.h.,
$\ddc\varphi_n\geq -\omega$, $\max_K\varphi_n\leq -n^3$,
$\max_X\varphi_n=0$. 
D'apr\`es la proposition 2.1, la s\'erie $\sum n^{-2}\varphi_n$ converge 
ponctuellement vers une fonction 
$\varphi$ q.p.s.h. On a $K\subset (\varphi=-\infty)$. 
\end{preuve}
\begin{proposition} Soit $(X,\omega)$ une vari\'et\'e
k\"ahl\'erienne compacte de dimension $k$.
Soit $\sigma$ la mesure associ\'ee \`a la forme volume $\omega^k$.
Alors, il existe des constantes $c>0$ et $\alpha>0$
telles que la mesure $\sigma$ soit $(c,\alpha)$-mod\'er\'ee. Plus 
pr\'ecis\'ement,
$\int_X\exp(-\alpha\varphi)\omega^k\leq c$ pour 
toute fonction $\varphi$ q.p.s.h. v\'erifiant $\ddc\varphi\geq
-\omega$ et $\max_X\varphi=0$.
En particulier, on a
$\Delta(X,\omega,\sigma,t)\leq c\exp(-\alpha r^{-1} t)$ pour
tout $t\in\R$, o\`u $r:=r(X,\omega)$.
\end{proposition}
\begin{preuve} 
Notons $\Ball(a,r)$ (resp. $\Ball_r$)
la boule de $\C^k$ de rayon
$r$ centr\'ee en $a$ (resp. centr\'ee en $0$).
Posons $\omega_0:=\ddc\|z\|^2$ la forme
euclidienne sur $\C^k$ et $\sigma_0:=\omega_0^k$.
Soit $\Psi_n:\Ball_4\longrightarrow X$ une famille finie
d'applications holomorphes injectives telles que les
ouverts $\Psi_n(\Ball_1)$ recouvrent $X$. Soit $A_1>0$,
une constante, telle que pour tout $n$, $\Psi_n^*(\omega)\leq A_1\omega_0$ et
$\Psi_n^*(\sigma)\leq A_1\sigma_0$ sur $\Ball_3$.
\par
Pour $\varphi$ comme dans la proposition,
posons $\varphi_n:=\varphi\circ\Psi_n$. Il suffit de montrer que
$\int_{\Ball_1}\exp(-\alpha'\varphi_n)\omega_0^k\leq c'$ pour
$c'>0$, $\alpha'>0$ ind\'ependants de $\varphi$.
D'apr\`es la proposition 2.1,
il existe une constante $A_2>0$
ind\'ependante de $\varphi$ telle que 
$\int|\varphi|\d\sigma\leq A_2$. On en d\'eduit que
$\sigma(\varphi<-M)\leq A_2M^{-1}$ pour tout $M>0$. Fixons un
$M>0$ assez grand tel que $A_2M^{-1}<
\sigma(\Psi_n(\Ball_1))$ pour tout $n$.
La derni\`ere relation implique que
$(\varphi<-M)$ ne peut contenir $\Psi_n(\Ball_1)$.
On peut donc choisir un point $a_n\in\Ball_1$ tel que
$\varphi_n(a_n)=\varphi(\Phi_n(a_n))\geq -M$.
\par
Posons $\psi_n:=\varphi_n+A_1(\|z\|^2-16)$. C'est une fonction
p.s.h. dans $\Ball_3$ v\'erifiant $\psi_n\leq \varphi_n$.
Montrons que $\int_{\Ball(a_n,2)}\exp(-\alpha'\psi_n)\omega_0^k
\leq c'$ avec $c'>0$ et
$\alpha'>0$ ind\'ependants de $\varphi$. On a
$\ddc\psi_n\geq 0$, $\psi_n(a_n)\geq -M-16A_1$ et $\psi_n\leq A_2$ sur
$\Ball(a_n,2)\subset \Ball_4$. Il suffit d'appliquer un
th\'eor\`eme de
H\"ormander \cite[p.97]{Hormander2} 
qui affirme que
$\int_{\Ball_1}\exp(-\phi)\omega_0^k\leq c'$ pour toute
fonction p.s.h. $\phi$ sur $\Ball_2$ avec $\phi(0)=0$ et $\phi\leq 1$.
On peut prendre donc $\alpha'=(A_2+16A_1+M)^{-1}$.  
\end{preuve}
\begin{remarque} \rm Les estimations peuvent \^etre raffin\'ees
en utilisant les r\'esultats de Skoda \cite{Skoda2}.
\end{remarque}
{\bf A.2}. {\it Estimation des constantes pour $\P^k$.}
\\
\ \par
Notons $\S^k$ (resp. $\S^{2k+1}$)
la sph\`ere unit\'e de $\R^{k+1}$ (resp. de $\C^{k+1}$) et
$\sigma_k$ (resp. $\sigma_{2k+1}$)
la mesure invariante de masse $1$ sur $\S^k$ (resp. $\S^{2k+1}$).
Soit $\pi:\C^{k+1}\setminus\{0\}\longrightarrow
\P^k$ la projection canonique. Notons $z=(z_0,\ldots,z_k)$ les
coordonn\'ees de $\C^{k+1}$.
On dira qu'une fonction $\Phi$ sur
$\C^{k+1}$ est {\it $\log$-homog\`ene} si on a $\Phi(\lambda
z)=\log|\lambda|+\Phi(z)$ pour tout $\lambda\in\C^*$. La fonction
$\log\|z\|$ \'etant p.s.h.
$\log$-homog\`ene, il existe une $(1,1)$-forme positive ferm\'ee
$\omega_\FS$ sur $\P^k$ telle que 
$\pi^*(\omega_\FS)=\ddc\log\|z\|$. C'est {\it la forme de
Fubini-Study}, elle est invariante par le groupe
unitaire.
\par
Notons $\Omega_\FS$ la mesure sur $\P^k$
associ\'ee \`a la forme volume $\omega_\FS^k$.
C'est {\it la mesure invariante de masse $1$} sur $\P^k$.
Soit $\varphi$ une
fonction q.p.s.h. sur $\P^k$ v\'erifiant
$\ddc\varphi\geq -\omega_\FS$. Posons
$\Phi:=\varphi\circ\pi+\log\|z\|$ et $\Phi(0):=-\infty$.
C'est une fonction p.s.h. $\log$-homog\`ene sur
$\C^{k+1}$ v\'erifiant
$\max_{\S^{2k+1}}\Phi=\max_{\P^k}\varphi$ et $\int_{\S^{2k+1}}
\Phi\d\sigma_{2k+1}=\int_{\P^k}\varphi\omega_\FS^k$.
\par
Notons $\R\P^k$ l'image de $\R^{k+1}$ par $\pi$. C'est un
sous-espace projectif r\'eel de dimension $k$ de
$\P^k$. 
Posons
$m_\FS:=\pi_*(\sigma_k)$ o\`u $\sigma_k$ est la mesure
de probabilit\'e sur $\S^k$ invariante par le groupe orthogonal.
On a aussi
$\max_{\S^{k}}\Phi=\max_{\R\P^k}\varphi$ et 
$\int_{\S^{k}}\Phi\d\sigma_k=\int\varphi\d m_\FS$.
\begin{proposition} On a
$$\RR^*_1(\P^k,\omega_\FS)\leq
\frac{1}{2}(1+\log k)\ \ \mbox{ et }\ \
\RR^*_2(\P^k,\omega_\FS,1)\leq 1+\log
k.$$
\end{proposition}
\begin{preuve} Rappelons que $r(\P^k,\omega_\FS)=1$. 
Soit $\varphi$ une fonction q.p.s.h. v\'erifiant
$\ddc\varphi\geq -\omega_\FS$. La fonction 
$\Phi:=\varphi\circ\pi+\log\|z\|$ est p.s.h., $\log$-homog\`ene.
Elle v\'erifie
$\int \Phi\d\sigma_{2k+1}=\int\varphi\omega_\FS^k$ et
$\max_{\S^{2k+1}}\Phi=\max_{\P^k}\varphi$.
Soit $m:=\max_{\P^k}\varphi$. 
D'apr\`es Alexander \cite[Theorem 2.2]{Alexander}, on a
\begin{eqnarray*}
\int_{\S^{2k+1}}\Phi\d\sigma_{2k+1} & \geq &
m+\int_{\S^{2k+1}}\log|z_1|\d\sigma_{2k+1}\\
& = & m-\frac{1}{2}\sum_{n=1}^k
\frac{1}{n}\geq m-\frac{1}{2}(1+\log k).
\end{eqnarray*}
On en d\'eduit que $m\leq
\frac{1}{2}(1+\log k)$ si $\int_{\P^k}\varphi \omega_\FS^k=0$.
Par d\'efinition (2.2), on a $\RR^*_1(\P^k,\omega_\FS)\leq
\frac{1}{2}(1+\log k)$. D'apr\`es la proposition 2.5, on a
$\RR^*_2(\P^k,\omega_\FS,1)\leq 1+\log k$. 
\end{preuve}
\par
La proposition suivante permet d'estimer les int\'egrales sur
$\S^k$ et $\S^{2k+1}$ 
en fonction d'int\'egrales sur des sous-espaces lin\'eaires.
\begin{proposition}
Soit $h$ une fonction mesurable positive sur la sph\`ere
unit\'e $\S^k$ de $\R^{k+1}$ (resp. $\S^{2k+1}$ de $\C^{k+1}$).
Soit $F$ un sous-espace r\'eel (resp. complexe)
de dimension $m$ de
$\R^{k+1}$ (resp. de $\C^{k+1}$),
$1\leq m\leq k$. Supposons que pour tout
sous-espace r\'eel (resp. complexe) $E$ de dimension $m+1$
contenant $F$ on ait $\int_{\S^k\cap E} h\d \sigma_m\leq A$ (resp.
$\int_{\S^{2k+1}\cap E}h \d \sigma_{2m+1}\leq A$)
o\`u $A>0$ est une constante.
Alors il existe une constante $c>0$ ind\'ependante
de $k,A$ et $h$ telle que $\int_{\S^k}h\d\sigma_k\leq cAk^{m/2}$
(resp. $\int_{\S^{2k+1}}h \d\sigma_{2k+1}\leq cAk^m$).
\end{proposition}
\begin{preuve} Consid\'erons d'abord le cas r\'eel.
Soit $(x_1,\ldots,x_{k+1})$ le syst\`eme de coordonn\'ees
canoniques de $\R^{k+1}$. Notons $\OO x_i$
les axes de coordonn\'ees et
$\OO x$ la demi-droite issue de $0$ passant par $x$.
On peut supposer que $F$ est le sous-espaces engendr\'e par les
$m$ derniers axes $\OO x_{k-m+2}$, $\ldots$, $\OO x_{k+1}$.
\par
Utilisons les coordonn\'ees
polaires de $\S^k$.
Soit $\theta_n$ l'angle entre $\OO x_{n+1}$ et la demi-droite 
issue de $0$ passant par le point
$(x_1,\ldots,x_{n+1},0\ldots,0)$, $1\leq n\leq k$.
Les angles $\theta_n$ v\'erifient $0\leq\theta_n< 2\pi$.
Pour tout $x\in \S^k$, on a $x_n=\cos\theta_{n-1}
\sin\theta_n\ldots\sin\theta_{k-1}$
o\`u on a pos\'e $\theta_0:=0$. 
L'\'el\'ement d'aire de $\S^k$ en m\'etrique euclidienne
est \'egal \`a
$$|(\sin\theta_2)(\sin\theta_3)^2\ldots(\sin\theta_k)^{k-1}|
\d\theta_1\ldots\d\theta_k.$$
Puisque $\sigma_k(\S^k)=1$, l'\'el\'ement d'aire
de $\sigma_k$ est \'egal \`a
$$c_k|(\sin\theta_2)(\sin\theta_3)^2\ldots(\sin\theta_k)^{k-1}|
\d\theta_1\ldots\d\theta_k$$
o\`u $c_k^{-1}$ est l'aire de $\S^k$
en m\'etrique euclidienne.
On a donc \cite[3.2.13]{Federer}
$$c_k=2^{-1}\pi^{-(k+1)/2}\Gamma((k+1)/2).$$
Evaluons l'int\'egrale de $h$
$$\int_{\S^k}h\d\sigma_k(K)=c_k\int h
|(\sin\theta_2)(\sin\theta_3)^2\ldots(\sin\theta_k)^{k-1}|
\d\theta_1\ldots\d\theta_k.$$
Posons $\theta:=(\theta_1,\ldots,\theta_k)$ et
$\theta':=(\theta_1,\ldots,\theta_{k-m})$.
Soit $\pi$ la projection de $[0,2\pi[^k$ dans $[0,2\pi[^{k-m}$
avec $\pi(\theta):=\theta'$.  
Par hypoth\`ese, on a
$$\int_{\S^k\cap\pi^{-1}(\theta')}
h\d\sigma_m\leq A\ 
\mbox{ pour tout }\theta'\in [0,2\pi[^{k-m}.$$
Or cette int\'egrale 
est \'egale \`a
$$c_m\int_{\S^k\cap\pi^{-1}(\theta')}h 
|(\sin\theta_{k-m+2})\ldots(\sin\theta_k)^m|
\d\theta_{k-m+1}\ldots\d\theta_k.$$
Utilisant le th\'eor\`eme de Fubini, on obtient 
\begin{eqnarray*}
\int_{\S^k} h\d \sigma_k & \leq &
c_k c_m^{-1}A\int_{[0,2\pi[^{k-m}}
|(\sin\theta_2)\ldots(\sin\theta_{k-m+1})^{k-m}|
\d\theta_1\ldots\d\theta_{k-m}\\
& = & A c_m^{-1}c_kc_{k-m}^{-1}. 
\end{eqnarray*}
Il suffit d'observer que
$$c_m^{-1}c_kc_{k-m}^{-1}=
2\pi^{1/2}\Gamma((m+1)/2)^{-1}
\Gamma((k+1)/2)\Gamma
((k-m+1)/2)^{-1}\leq c k^{m/2}$$
o\`u $c>0$ est une constante.
\\
\par
La preuve du cas complexe est identique \`a celle
du cas r\'eel. Il suffit d'utiliser les coordonn\'ees polaires
suivantes pour $z=(z_1,\ldots,z_{k+1})\in\S^{2k+1}$:
$$z_n=\cos\theta_{n-1}\sin\theta_n\ldots\sin
\theta_{k-1}\exp(i\alpha_n)$$
avec $0\leq\theta_n<2\pi$ pour $1\leq n\leq k$,
$\theta_0=0$
et $0\leq\alpha_n< 2\pi$ pour $1\leq n\leq k+1$.
\end{preuve}
\begin{proposition} Il existe des constantes $c>0$ et $\alpha>0$
ind\'ependantes de $k$,
telles que $\int_{\P^k}\exp(-\alpha\varphi)\omega_\FS^k\leq ck$
pour toute fonction q.p.s.h. $\varphi$ v\'erifiant
$\ddc\varphi\geq -\omega_\FS$ et $\max_{\P^k}\varphi=0$. En
particulier, on a
 $\Delta(\P^k,\omega_\FS,\Omega_\FS,t)\leq
ck \exp(-\alpha t)$ pour tout $t\in\R$.
\end{proposition}
\begin{preuve} Soit $\varphi$ une fonction comme dans la
proposition. Soit $a\in\P^k$ tel que $\varphi(a)=0$. Posons
$\Phi:=\varphi\circ \pi$ et $\Phi(0)=-\infty$. Soit $\P^1_\xi$
une droite projective passant par $a$, $F:=\pi^{-1}(a)\cup\{0\}$
et $E_\xi:=\pi^{-1}(\P^1_\xi)\cup\{0\}$. Le plan complexe
$E_\xi$ contient la droite $F$.
Soient $c_1$ et $\alpha$ les constantes v\'erifiant la proposition
A.7 pour la droite projective $(\P^1,\omega_\FS)$ (\cad que pour
$k=1$).
\\
\par 
On a
pour tout $\P^1_\xi$
$$\int_{\S^{2k+1}\cap E_\xi}\exp(-\alpha\Phi)\d\sigma_3\leq c_1.$$
D'apr\`es la proposition A.6, on a
$$\int_{\P^k}\exp(-\alpha \varphi)\omega_\FS^k
=\int_{\S^{2k+1}}\exp(-\alpha\Phi)
\d\sigma_{2k+1}\leq ck$$
o\`u $c>0$ est une constante ind\'ependante de $k$.
\par
On a montr\'e que la mesure $\Omega_\FS=\omega_\FS^k$ est
$(ck,\alpha)$-mod\'er\'ee. Puisque $r(\P^k,\omega_\FS)=1$, 
d'apr\`es la proposition 2.5, on a
$$\Delta(\P^k,\omega_\FS,\Omega_\FS,t)\leq ck\exp(-\alpha t).$$
\end{preuve}
\par
La proposition suivante explique pourquoi les estimations 
pour $\R\P^k$ et $\P^k$ sont essentiellement les m\^emes.
Rappelons que dans $\P^k$ on a $\capacity(\P^k)=1$.
\begin{proposition}
Soit $\capacity(\R\P^k)$ la capacit\'e de $\R\P^k$
dans $\P^k$. Alors pour tout $k\geq 2$ et
toute fonction $\varphi$ q.p.s.h. sur
$\P^k$ v\'erifiant $\ddc\varphi\geq -\omega_\FS$ et
$\max_{\P^k}\varphi=0$, on a $\max_{\R\P^k}\varphi\geq \log
\capacity(\R\P^2)$.
En particulier, 
$\capacity(\R\P^k)=\capacity(\R\P^2)$ pour tout $k\geq 2$.
\end{proposition}
\begin{preuve} Soit $\varphi$ comme ci-dessus. Soient
$a\in\P^k$ et $b\in\R\P^k$ tels que $\varphi(a)=0$ et
$\varphi(b)=\max_{\R\P^k}\varphi$.
Observons que $\C^k$ est  r\'eunion des
sous-espaces complexes $F$ de
dimension $2$ v\'erifiant $\dim_\R(F\cap\R^k)=2$.
De m\^eme $\C^{k+1}$ est r\'eunion des sous-espaces complexes $E$ de
dimension $3$ contenant
une droite r\'eelle fixe de $\R^{k+1}$ et
v\'erifiant $\dim_\R(E\cap\R^{k+1})=3$.
Donc $\P^k$ est r\'eunion des plans
projectifs $P$ de $\P^k$ passant par
$b$ et v\'erifiant $\dim_\R(P\cap \R\P^k)=2$. Soit
$P$ un tel plan contenant $a$.
Par d\'efinition de la capacit\'e sur $P\simeq\P^2$ pour
$P\cap\R\P^k\simeq \R\P^2$, on a $\varphi(b)\geq \log \capacity
(\R\P^2)$. Donc $\capacity(\R\P^k)\geq \capacity(\R\P^2)$. 
On obtient l'autre in\'egalit\'e en observant que toute fonction
q.p.s.h. $\varphi$ sur $\P^2$ avec $\ddc\varphi\geq -\omega_\FS$
se prolonge en fonction q.p.s.h. $\widetilde\varphi$ sur $\P^k$ avec
$\ddc\widetilde \varphi\geq -\omega_\FS$ et
$\max_{\P^k}\widetilde\varphi=\max_{\P^2}\varphi$ (ceci se voit
ais\'ement sur le relev\'e de $\widetilde\varphi$ \`a $\C^{k+1}$).  
\end{preuve}
\begin{proposition} Il existe des constantes
$c>0$ et $\alpha>0$
ind\'ependantes de $k$ telles qu'on ait
$\int_{\P^k}\exp(-\alpha\varphi)\d m_\FS\leq c\sqrt{k}$,
$\RR_i(\P^k,\omega_\FS,m_\FS)\leq c(1+\log k)$
pour $i=1,2,3$ et
$\Delta(\P^k,\omega_\FS,m_\FS,t)\leq c\sqrt{k}\exp(-\alpha t)$
pour tout $t\in\R$.
\end{proposition}
\begin{preuve}
On peut reprendre la d\'emonstration de H\"ormander
\cite[p.98]{Hormander2} pour les fonctions $\psi$ sous-harmoniques
sur le disque unit\'e de $\C$ qui v\'erifient $\psi(0)=0$ et
$\psi\leq 1$. On remplace la mesure de Lebesgue par la mesure $m$
sur le cercle $|z|=1/2$. On trouve $\int\exp(-\psi/2)\d m\leq
c_1$ o\`u $c_1>0$ est une constante ind\'ependante de $\psi$.
Un argument semblable \`a celui de la
proposition A.3 permet de montrer
qu'il existe $c'>0$ et $\alpha>0$ telle que
dans $\P^1$ on a $\int\exp(-\alpha'\varphi)\d m_\FS\leq c'$ pour
toute fonction $\varphi$ q.p.s.h. v\'erifiant $\ddc\varphi\geq
-\omega_\FS$ et $\max_{\P^1}\varphi=0$. Le passage \`a $\P^k$ avec
l'estimation $\int_{\P^k}\exp(-\alpha\varphi)\d m_\FS\leq c\sqrt{k}$
est une simple application de la
proposition A.6 comme dans la proposition A.7.
On en d\'eduit que 
$\Delta(\P^k,\omega_\FS,m_\FS,t)\leq c\sqrt{k}\exp(-\alpha t)$
pour tout $t\in\R$. 
\\
\ \par
Montrons que $\RR_1(\P^k,\omega_\FS,m_k)\leq c(1+\log k)$.
Soit $\varphi$ une fonction q.p.s.h.
v\'erifiant $\ddc\varphi\geq -\omega_\FS$ et
$\max_{\P^k}\varphi=0$.
D'apr\`es (2.4), il suffit de 
v\'erifier que $-\int\varphi\d m_\FS\leq c(1+\log k)$ pour
un $c>0$ ind\'ependant de $k$ et de $\varphi$.
On peut supposer que $k\geq 2$. 
\par
D'apr\`es l'estimation
 $\int_{\P^k}\exp(-\alpha\varphi)\d m_\FS\leq c\sqrt{k}$
ci-dessus, on a pour tout $t_0\geq 0$
\begin{eqnarray*}
-\int\varphi\d m_\FS & \leq & t_0+\int_{t_0}^{+\infty}
m_\FS(\varphi<-t)\d t
 \leq  t_0+\int_{t_0}^{+\infty} c\sqrt{k} \exp(-\alpha t)\d t\\
& = & t_0+ c\sqrt{k}\alpha^{-1}\exp(-\alpha t_0).
\end{eqnarray*}
Pour $t_0=2^{-1}\alpha^{-1}\log k$, on obtient l'in\'egalit\'e
voulue.
\par
D'apr\`es la proposition 2.5 et la proposition A.5,
on a $\RR_i(\P^k,\omega_\FS,m_\FS)\leq
c(1+\log k)$ pour $i=2,3$ et pour une constante $c>0$
ind\'ependante de $k$.
\end{preuve}
\par
\
\\
{\bf A.3}. {\it Espaces produits et espaces multiprojectifs.}  
\\
\ \par
On associe
\`a la vari\'et\'e $X:=X_1\times X_2$ la forme de K\"ahler
$$\omega:=
c_{12}(\pi_1^*\omega_1+\pi_2^*\omega_2)$$
o\`u $\pi_1$, $\pi_2$ sont
les projections canoniques de $X$ sur $X_1$ et $X_2$
et o\`u
$c_{12}>0$ est telle que $\int_X\omega^{k_1+k_2}=1$.
La constante $c_{12}$ est calcul\'ee par la formule
$$c_{12}^{-k_1-k_2}= {k_1+k_2 \choose k_1}.$$
Consid\'erons deux mesures de probabilit\'e $\mu_1$ sur $X_1$ et
$\mu_2$ sur $X_2$. Notons $\mu$ le produit de $\mu_1$ et $\mu_2$.
C'est une mesure de probabilit\'e sur $X$. 
Posons $r:=r(X,\omega)$.
\begin{proposition} Soient $X$, $\omega$, $\Omega$, $r$ et $\mu$
comme ci-dessus. Soit $\RR$ une constante v\'erifiant $\RR\geq
\RR_1(X_i,\omega_i,\mu_i)$ pour $i=1,2$.
Supposons qu'il existe $c>0$ et
$\alpha>0$ telles que $\Delta(X_i,\omega_i,\mu_i,t)\leq
c\exp(-\alpha t)$ pour $t\in\R$ et $i=1,2$. Alors on a
$$\RR_1(X,\omega,\mu)\leq 2r\RR+2r\alpha^{-1}\log c+4\alpha^{-1}r$$
et
$$\Delta(X,\omega,\mu,t)\leq
2c\exp(\alpha\RR)\exp(-\alpha r^{-1} t/2)$$
pour tout $t\in\R$.
\end{proposition}
\begin{preuve} Fixons une fonction $\psi$ sur $X=X_1\times X_2$
telle que $\max_X \psi=0$ et $\ddc\psi\geq -r\omega$. Soit
$(a_1,a_2)\in X$ un point tel que $\psi(a_1,a_2)=0$. On veut
estimer la $\mu$ mesure de l'ensemble
$E:=(\psi<-t)$ pour $t\geq 0$. Posons
$$F:=\{x_2\in X_2,\ \psi(a_1,x_2)< -t/2\}$$
et 
$$E_{x_2}:=\{x_1\in X_1,\ \psi(x_1,x_2)< -t\}.$$
D\'efinissons
$$E':=\bigcup_{x_2\in X_2\setminus F}
(E_{x_2}\times\{x_2\}).$$  
On a $E\subset \pi_2^{-1}(F)\cup E'$.
\\
\ \par
Estimons la mesure
de $\pi_2^{-1}(F)$.
Si $\psi_1(x_2):=\psi(a_1,x_2)$, on a
$\max_{X_2}\psi_1=\psi_1(a_2)=0$. Posons
$\psi_2:=\psi_1-\int\psi_1\d\mu_2$. On a $\int\psi_2\d\mu_2=0$,
$\psi_2\geq \psi_1$ et $\ddc\psi_2\geq -r\omega_2$. Par
d\'efinition de $\RR_1(X_2,\omega_2,\mu_2)$,
puisque $r(X_2,\omega_2)\geq 1$, on a
$$-\int\psi_1\d\mu_2=\max_{X_2}\psi_2\leq
r\RR_1(X_2,\omega_2,\mu_2)\leq r\RR.$$
D'o\`u
\begin{eqnarray*}
\mu_2(F) &\leq & \mu_2(\psi_2\leq r\RR-t/2)
=\mu_2(r^{-1}\psi_2\leq\RR-r^{-1}t/2)\\
& \leq & \Delta(X_2,\omega_2,\mu_2,\RR-r^{-1}t/2)
\leq
c\exp(\alpha \RR-\alpha r^{-1}t/2)\\
& = & c\exp(\alpha
\RR)\exp(-\alpha r^{-1}t/2).
\end{eqnarray*}
Donc $\mu(\pi_2^{-1}(F))\leq c\exp(\alpha \RR)\exp(-\alpha
r^{-1}t/2)$.
\\
\ \par
Estimons la mesure de $E'$. Pour $x_2\in X_2\setminus
F$, posons $\psi_3(x_1):=\psi(x_1,x_2)$. On a $\psi_3\leq
0$, $\max_{X_1}\psi_3\geq \psi(a_1,x_2)\geq -t/2$ et
$\ddc\psi_3\geq -r\omega_1$. Posons
$\psi_4:=\psi_3-\int_{X_1}\psi_3\d\mu_1$. On a
$$-\int\psi_3\d\mu_1\leq\max_{X_2}\psi_4+t/2\leq
r\RR_1(X_2,\omega_2,\mu_2)+t/2\leq r\RR+t/2.$$
et
\begin{eqnarray*}
\mu_1(E_{x_2}) &\leq & \mu_1(\psi_4\leq r\RR-t/2)\leq
c\exp(\alpha \RR-\alpha r^{-1}t/2)\\
& = & c\exp(\alpha
\RR)\exp(-\alpha r^{-1}t/2).
\end{eqnarray*}
D'apr\`es le th\'eor\`eme de Fubini, on a
$\mu(E')\leq c\exp(\alpha \RR)\exp(-\alpha r^{-1}t/2)$.
\\
\ \par
On d\'eduit des estimations pr\'ec\'edentes que pour tout $t\geq 0$
$$\mu(\psi<-t)\leq 2c\exp(\alpha\RR)\exp(-\alpha r^{-1}t/2).$$
C'est aussi vrai pour tout $t\in\R$ car $\psi\leq 0$.  
Si $\varphi$ une fonction q.p.s.h. sur $X$ telle que
$\ddc\varphi\geq -\omega$ et $\int\varphi\d\mu=0$, on peut
appliquer l'estimation pr\'ec\'edente \`a la fonction
$\psi:=\varphi-\max_X\varphi$. Cette derni\`ere fonction est plus
petite que $\varphi$ et elle v\'erifie $\max_X\psi=0$.
On en d\'eduit que pour tout $t\geq 0$
$$\mu(\psi<-t)\leq 2c\exp(\alpha\RR)\exp(-\alpha r^{-1}t/2).$$
Donc $\Delta(X,\omega,t)\leq 2c\exp(\alpha
\RR)\exp(-\alpha r^{-1}t/2)$.
\\
\ \par
Estimons $\RR_1(X,\omega,\mu)$. Soit $\psi$ comme
ci-dessus, il faut montrer que
$$-\int\psi\d\mu\leq
2r\RR+2r\alpha^{-1}\log c +4\alpha^{-1}r.$$
On a pour tout $t_0\geq 0$
\begin{eqnarray*}
-\int\psi\d\mu & = & \int_0^{+\infty} \mu(\psi\leq -t)\d t\\
& = & \int_0^{t_0} \mu(\psi\leq -t)\d t
+\int_{t_0}^{+\infty} \mu(\psi\leq -t)\d t\\
& \leq & \int_0^{t_0} \d t
+\int_{t_0}^{+\infty} 2c\exp(\alpha\RR)\exp(-\alpha r^{-1}t/2)\d t\\
& = & t_0 + 4c\exp(\alpha\RR)\alpha^{-1} r \exp(-\alpha
r^{-1}t_0/2)
\end{eqnarray*}
En prenant $t_0=2r\RR+2r\alpha^{-1}\log c$,
on obtient
$$-\int\psi\d\mu\leq 2r\RR + 2r\alpha^{-1}\log c
+4\alpha^{-1}r.$$
\end{preuve}
\par
Consid\'erons l'espace multi-projectif
$\P^{k,l}:=\P^k\times\cdots\times\P^k$ ($l$ fois).
Notons $\pi_i$ la projection de $\P^{k,l}$ sur le $i$-\`eme
facteur. Posons $\omega_\MP:=c_{k,l}\sum\pi_i^*(\omega_\FS)$ et
soit $\Omega_\MP$ la mesure associ\'ee \`a la forme volume
$\omega_\MP^{kl}$.
La constante $c_{k,l}>0$ est choisie de
sorte que $\Omega_\MP$ soit une mesure de probabilit\'e. On a
\begin{eqnarray}
(c_{k,l})^{-kl} & = &
{kl\choose k}{kl-k \choose k}\ldots {2k \choose k}.
\end{eqnarray}
Observons que $c_{k,l}\leq 1$ et que
si $l$ est fix\'e, $c_{k,l}$ est
minor\'ee par une constante positive   
ind\'ependante de $k$. Pour cela, il suffit d'utiliser la formule d'\'equivalence
de Stirling $n!\simeq \sqrt{2\pi n} n^n \exp(-n)$.
Notons $m_\MP$ la mesure produit des $m_\FS$
sur $\R\P^{k,l}:=\R\P^k\times\cdots\times\R\P^k$ ($l$ fois).
\begin{proposition} Il existe des constantes $c>0$, $\alpha >0$
et $m>0$,
qui ne d\'ependent que de $l$, telles que pour tout $k\geq 1$ on ait
$\RR_i(\P^{k,l},\omega_\MP,\mu)\leq c(1+\log k)$,
$i=1,2,3$, 
$\Delta(\P^{k,l},\omega_\MP,\mu,t)\leq ck^m\exp(-\alpha t)$ pour
tout $t>0$ et pour $\mu=\Omega_\MP$ ou $\mu=m_\MP$.
\end{proposition}
\begin{preuve}
Observons que la constante $r(\P^{k,l},\omega_\MP)$ est major\'ee
par $l/c_{k,l}$. Il suffit d'appliquer les propositions A.5, A.7, A.9 et
A.10.  
\end{preuve}
\small

Tien-Cuong Dinh et Nessim Sibony,\\
Math\'ematique - B\^at. 425, UMR 8628, 
Universit\'e Paris-Sud, 91405 Orsay, France. \\
E-mails: Tiencuong.Dinh@math.u-psud.fr et
Nessim.Sibony@math.u-psud.fr.
\end{document}